\documentclass[11pt]{article}
\RequirePackage[OT1]{fontenc}
\RequirePackage{amsthm,amsmath}
\RequirePackage{hyperref}
\usepackage{amsfonts}
\usepackage{amssymb}
\usepackage{amstext}
\usepackage{parskip}
\usepackage{fullpage}
\usepackage{graphicx}
\usepackage{subcaption}

\RequirePackage[numbers]{natbib}
\usepackage{amsfonts,amssymb,amsthm,amsmath,enumerate}

\newcommand{\prob}[1]{\ensuremath{\mathbb P}\left(#1\right)}

\newcommand{\R}{\ensuremath{\mathbb R}}

\newcommand{\Z}{\ensuremath{\mathbb Z}}

\newcommand{\size}[1]{\ensuremath{\left|#1\right|}}

\newcommand{\argmin}{\operatorname{argmin}}

\newcommand{\silent}[1]{}

\newcommand{\Ball}{{\mathbf B}}

\newcommand{\MC}{{\mathcal C}}

\newcommand{\M}{{\mathcal M}}

\newcommand{\cov}{\textsf{Cov}}

\newcommand{\inv}[1]{\frac{1}{#1}}
\newcommand{\abs}[1]{\left\lvert#1\right\rvert}
\newcommand{\twonorm}[1]{\left\lVert#1\right\rVert_2}
\newcommand{\norm}[1]{\left\lVert#1\right\rVert}
\newcommand{\ldot}{\ldots}

\def\cp{{\mathcal P}}
\def\cpk{{\mathcal P}_k}
\def\Z{{\mathbb Z}}
\def\E{{\mathbb E}}

\def\M{{\mathcal M}}

\def\supp{\mathop{\text{supp}\kern.2ex}}
\def\argmin{\mathop{\text{arg\,min}\kern.2ex}}
\let\hat\widehat
\let\tilde\widetilde

\def\qed{\hskip1pt $\;\;\scriptstyle\Box$}

\newcommand{\vecone}{{\bf 1}}
\newcommand{\vc}{{\bf c}}

\newcommand{\BW}{\ensuremath{\mathbb{W}}}
\def\argmax{\mathop{\text{arg\,max}\kern.2ex}}

\def\fatnorm#1{|\kern-.2ex|\kern-.2ex| #1 |\kern-.2ex|\kern-.2ex|}

\newcommand{\fnorm}[1]{\left\lVert#1\right\rVert_F}

\newcommand{\onenorm}[1]{\left\lVert#1\right\rVert_{1}}

\def\SDP{\mathop{\text{SDP}\kern.2ex}}
\newcommand{\func}[1]{\ensuremath{\mathrm{#1}}}
\newcommand{\diag}{\func{diag}}

\newcommand{\offd}{\func{offd}}

\newcommand{\Sp}{\mathbb{S}}
\newcommand{\ip}[1]{\;\langle{\,#1\,}\rangle\;}

\newcommand{\beq}{\begin{equation}}
\newcommand{\eeq}{\end{equation}}
\newcommand{\ben}{\begin{eqnarray}}
\newcommand{\een}{\end{eqnarray}}
\newcommand{\bnum}{\begin{enumerate}}
\newcommand{\enum}{\end{enumerate}}
\newcommand{\bit}{\begin{itemize}}
\newcommand{\eit}{\end{itemize}}
\newcommand{\bens}{\begin{eqnarray*}}
\newcommand{\eens}{\end{eqnarray*}}
\newcommand{\G}{{\mathcal G}}

\DeclareMathOperator*{\Span}{span}

\newcommand{\tr}{{\rm tr}}

\newcommand{\mvec}[1]{\rm{vec}\left\{\,#1\,\right\}}
\newcommand{\opt}{\ensuremath{\text{opt}}}
\newcommand{\CC}{{\mathcal C}}

\newtheorem{theorem}{Theorem}[section]

\newtheorem{fact}[theorem]{Fact}
\newtheorem{lemma}[theorem]{Lemma}
\newtheorem{proposition}[theorem]{Proposition}

\newtheorem{definition}[theorem]{Definition}

\newtheorem{remark}[theorem]{Remark}
\newtheorem{corollary}[theorem]{Corollary}

\def\qed{\hskip1pt $\;\;\scriptstyle\Box$}

\newenvironment{proofof}[1]{\hspace*{20pt}{\it Proof}{ of #1}.\hskip10pt}{\qed\vskip5pt}
\newenvironment{proofof2}{\hskip10pt}{\qed\vskip5pt}

\begin{document}

\title{Semidefinite programming on population clustering: a local analysis}

  \author{Shuheng Zhou\\
    University of California, Riverside, CA 92521}

\date{}

\maketitle

\begin{abstract}

In this paper, we consider the problem of partitioning a small data sample of
  size $n$ drawn from a mixture of $2$ sub-gaussian distributions. In
  particular, we design and analyze two computational efficient algorithms to
  partition data into two groups approximately according to their population of
  origin given a small sample in a recent paper (Zhou 2023a).  Our work is
  motivated by the application of clustering individuals according to their
  population of origin using markers, when the divergence between any two of
  the populations is small.    Moreover, we are interested in the case that
  individual features are of low average quality $\gamma$, and we want to use
  as few of them as possible to correctly partition the sample.  Here we use $p
  \gamma$ to denote the $\ell_2^2$ distance between two population centers
  (mean vectors), namely, $\mu^{(1)}, \mu^{(2)} \in \R^p$.  We allow a full
  range of tradeoffs  between $n, p, \gamma$ in the sense that partial recovery
  (success  rate $< 100\%$) is feasible once the signal to noise ratio $s^2 :=
  \min\{np \gamma^2,  p \gamma\}$ is lower bounded by a constant.  Our work
  builds upon the semidefinite relaxation of an integer quadratic program that
  is formulated essentially as finding the  maximum cut on a graph, where edge
  weights in the cut represent dissimilarity  scores between two nodes based on
  their $p$ features in Zhou (2023a). More importantly, we prove that the
  misclassification error decays exponentially with respect to the SNR $s^2$ in
  the present  paper. The significance of such an exponentially decaying error
  bound is: when $s^2 =\Omega(\log n)$, perfect recovery of the cluster
  structure is  accomplished. This result was introduced in Zhou (2023a)
  without a proof. We  therefore present the full proof in the present work.
\end{abstract}

\section{Introduction}
\label{sec:intro}
We explore a type of classification problem that arises in the context
of computational biology.
The biological context for this problem is we are
given DNA information from $n$ individuals from $k$ populations 
of origin and we wish to classify each individual into the correct category.  
DNA contains a series of markers called SNPs  (Single Nucleotide Polymorphisms), each of which has two variants (alleles). We use bit $1$ and bit $0$ to denote them.
The problem is that we are given a small
sample of size $n$, e.g., DNA of $n$ individuals (think of $n$ in the hundreds or thousands),
each described by the values of $p$ {\em features} or {\em markers},
e.g., SNPs (think of $p$ as an order of magnitude larger than $n$).
Our goal is to use these features to classify the individuals according to
their population of origin.
Given the population of origin of an individual, the genotypes can be
reasonably assumed to be generated by drawing alleles independently
from the appropriate distribution.

Suppose we are given a data matrix $X \in \R^{n \times p}$ with
samples from two populations $\MC_1, \MC_2$, such  that 
\ben 
\label{eq::Xmean}
\forall i \in  \MC_g, \; \; \E (X_{i j}) =  \mu^{(g)}_{j}  \; \; \;  
g =1,  2, \forall j =1, \ldots, p. 
\een 
Our goal is to estimate the group membership
vector $u_2 \in \{-1,1\}^n$ such that
\ben
\label{eq::u2}
{u}_{2,j} = 1 \; \text{ for } \; j \in \MC_1 \;  \text{ and } \;
{u}_{2,j} =- 1 \; \text{ for } \; j \in \MC_2,
\een
where the sizes of clusters $\abs{\MC_{\ell}} =: n_{\ell}, \forall \ell$ may not be
the same. Our ultimate goal is to estimate the solution to the discrete optimization problem:
\ben
\label{eq::quadratic}
\text{maximize} \; \; x^T R x \quad \text{subject to} \; \; x 
\in  \{-1, 1\}^n,
\een
where $R$ is a static reference matrix to be specified.
It was previously shown that, in expectation, among all balanced cuts
in the complete graph formed among $n$ vertices (sample points), the
cut of maximum weight corresponds to the correct partition of the $n$
points according to their distributions in the balanced case ($n_1 = n_2  = n/2$).
Here the weight of a cut is 
the sum of weights across all edges in the cut, and the edge weight 
equals the Hamming distance between the bit vectors of the two
endpoints; See~\cite{CHRZ07,Zhou06}.
Under suitable conditions, the statement above also holds with high 
probability (w.h.p.); The analyses in~\cite{CHRZ07} and~\cite{Zhou06}
focused on the high dimensional setting, where $p \gg n$.

Features have slightly different probabilities 
depending on which population the individual belongs to.
Denote by $\Delta^2$ the $\ell_2^2$ distance between two  population
centers (mean vectors), namely, $\mu^{(1)}, \mu^{(2)} \in \R^p$.
We focus on the case where $p >n$, although it is not needed.
 Note that $\Delta$ measures the Euclidean distance between
 $\mu^{(1)}$ and $\mu^{(2)}$ and thus represents
 their separation.
The objective we consider is to minimize the total data size $D=n p$
needed to correctly classify the individuals in the sample as a
function of the ``average quality'' $\gamma$ of the features:
\ben
\label{eq::Delta}
 \gamma := \Delta^2/p, \; \text{ where} \; \; \Delta^2  :=
 \sum_{k=1}^p (\mu^{(1)}_k - \mu^{(2)}_k)^2 \;
\text{ and }  \mu^{(i)} = (\mu^{(i)}_1, \ldots, \mu^{(i)}_p) \in \R^p, i =1, 2.
\een
In other words, in the context of population
clustering, it has been previously shown one can use a random instance of the integer quadratic program: 
\ben
\label{eq::graphcut}
(Q) \quad \text{maximize} \; \; x^T A x \quad \text{ subject to} \quad x \in  \{-1, 1\}^n
\een
to identify the correct partition of nodes according to their population
of  origin w.h.p. so long as the data size $D$ is sufficiently large
and the separation metric is at the order of $\Delta^2 = \Omega(\log
n)$.  
Here $A =(a_{ij})$ is an $n \times n$ symmetric matrix,
where for $1\le i, j \le n$, $a_{ij} = a_{ji}$ denotes the edge weight
between nodes $i$ and $j$, computed from the individuals’ bit vectors.
This result is structural, rather than algorithmic.

The integer quadratic program~\eqref{eq::quadratic}
(or~\eqref{eq::graphcut}) is NP-hard~\cite{Karp72}.
In a groundbreaking paper~\cite{GW95}, Goemans and Williamson
show that one can use semidefinite program (SDP) as relaxation to solve
these approximately. See references therein for earlier works.
In the present work, we use semidefinite relaxation of 
the graph cut problem~\eqref{eq::graphcut}, which was originally formulated 
in~\cite{CHRZ07,Zhou06} in the context of population clustering.
More generally, one may consider semidefinite relaxations for
the following sub-gaussian mixture model with $k$ centers
(implicitly, with rank-$k$ mean matrix embedded), where
we have $n$ observations $X_1, X_2 \ldot,  X_n \in \R^{p}$:
\ben
\label{eq::model}
X_i = \mu^{(\psi_i)} + \Z_i,
\een
where $\Z_1, \ldot, \Z_n \in \R^{p}$ are independent, sub-gaussian, mean-zero, random vectors and
$\psi_i: i \to \{1, \ldots, k\}$ assigns node $i$ to a group $\MC_j$ with the 
mean $\mu^{(j)} \in \R^{p}$ for some $j \in [k]$.
Here we denote by $[k]$ the set of integers $\{1, \ldots, k\}$.
Hence, each row vector $X_i$ of data matrix $X$ is a $p$-dimensional 
sub-gaussian random vector and we assume rows are independent.

The proposed semidefinite relaxation framework in~\cite{Zhou23a} was
inspired by~\cite{GV15}. The important distinction of the present work from that as considered
in~\cite{GV15} is: to estimate the group membership  vector ${u}_2$,
we replace the random adjacency matrix stochastic block models with an
instance of symmetric matrix computed from the centered data matrix,
which we now introduce.
First, we define the {\it global center} of the data.
Let   $X$ be as in~\eqref{eq::Xmean} with  row vectors 
$X_1, \ldots, X_n \in \R^{p}$. 
\begin{definition}{\textnormal{ \bf (The global center)}}
  \label{def::estimators}
Denote by $\Z_i = X_i - \E (X_i)$.
Then the global center is the average over $n$ row vectors $X_1,
\ldots, X_n \in \R^{p}$, for $X$ as in~\eqref{eq::Xmean}, 
\ben
\label{eq::muvector}
\label{eq::muhat}
\hat{\mu}_n & = & \inv{n}\sum_{i=1}^n X_i, \quad\text{ and } \quad
\hat{\mu}_n - \E \hat{\mu}_n :=\inv{n}
\sum_{i=1}^n \Z_i = \inv{n} \sum_{i=1}^n  X_i - \E (X_i).
\een
\end{definition}
\noindent{\bf Estimators and convex relaxation.}
\label{sec::estimators}
Let $\vecone_n = [1, \ldots, 1] \in \R^n$ denote a vector of all $1$s
and  $E_n  := \vecone_n \vecone_n^T$.
As in many statistical problems,  one simple but crucial step is 
to first obtain the centered data $Y$. 
In the present work, to find a convex relaxation, we will first construct a matrix $Y$ by subtracting the sample mean $\hat{\mu}_n \in \R^{p}$ as computed from \eqref{eq::muhat} from each row vector $X_i$ of the data matrix $X$.
Denote by
\ben
\label{eq::defineY}
Y & = & X - P_1 X, \quad \text{ where}
\quad P_1 := \inv{n} \vecone_n \vecone_n^T = E_n/n
\een
is a projection matrix.
Loosely speaking, this procedure is 
called ``global centering'' in the statistical literature, for 
example, see~\cite{horns19}.
Given matrix $Y$, we construct:
\ben
\label{eq::defineAintro}
A & := & Y Y^T -\lambda (E_n - I_n), \;\; \text{ where} \;
\lambda  = \frac{2}{n(n-1)}\sum_{i<j} \ip{Y_i, Y_j}
\een
and $I_n$ denotes the  identity matrix.
To estimate the group membership vector ${u}_2 \in \{-1, 1\}^n$, 
we consider the following semidefinite optimization problem from~\cite{Zhou23a}:
\ben 
\label{eq::sdpmain}
\text{(SDP) } && \text{maximize}
\; \;  \ip{A, Z} \;\; \text{ subject to} \quad Z \succeq 0, 
\; I_n \succeq \diag(Z).
\een
Here $Z \succeq 0$ indicates that the matrix $Z$ is constrained to be 
positive semidefinite, $A \succeq B$ means that $A - B \succeq 0$, and
the inner product of matrices $\ip{A, B} = \tr(A^T B)$.
Here and in the sequel, denote by 
$$\M^{+}_{G} := \{Z: Z \succeq 0,  I_n \succeq \diag(Z)\} \subset [-1,
1]^{n \times n}$$
the set of positive semidefinite matrices whose entries 
are bounded by 1 in absolute value.
Since all possible solutions of (Q) are feasible for SDP, 
the optimal value of SDP is an upper bound on the optimal value of 
(Q).  Moreover, SDP~\eqref{eq::sdpmain} is shown to be equivalent to the optimization problem 
\ben 
\label{eq::hatZintro}
\text{(SDP1)} && \text{maximize}
\; \;  \ip{YY^T - \lambda E_n,  Z}  \quad   \text{ subject to} \quad 
Z \succeq 0,  \diag(Z) = I_n.
\een
See the supplementary Proposition~\ref{prop::optsol}.
The (global) analysis framework for the semidefinite relaxation
in~\cite{GV15} was set in the context of community detection in
sparse networks, where $A$ represents the symmetric adjacency matrix of a random
graph. In other words, they study the semidefinite relaxation of the
integer program~\eqref{eq::graphcut}, where an $n \times n$ 
random matrix $A$ (observed) is used to replace the hidden
static $R$ in the original problem~\eqref{eq::quadratic} such that $\E (A) = R$.
The innovative proof strategy of~\cite{GV15} is to apply the Grothendieck's
inequality for the random error $A - \E(A)$ rather than the original
matrix $A$ as considered in the earlier literature. We call this
approach the global analysis in~\cite{Zhou23a},
following~\cite{CCLN21}.
Throughout this paper, denote by $w_{\min} := \min_{j =1, 2} w_j$,
where $w_j = n_j/n$.  We use $n_{\min}  := n w_{\min}$ and $n_{\max} := n w_{\max}$ to 
represent the size of the smallest and the largest cluster 
respectively. For a symmetric matrix $M$, let $\lambda_{\max}(M)$ and 
$\lambda_{\min}(M)$ be the largest and the smallest eigenvalue of $M$ respectively. 
The operator norm $\twonorm{M}$ is defined to be 
$\sqrt{\lambda_{\max}(M^T M)}$. 

\subsection{Our approach and contributions}
Our approach of using the centered data matrix $Y$ in SDP1 
is novel to the best of our knowledge. 
Following~\cite{Zhou23a}, we use
$YY^T$ and the corresponding $A$ as the input to our optimization
algorithms, ensuring both computational efficiency and statistical
convergence, especially in the low signal-to-noise ratio (SNR) case when $s^2 \asymp 
(p \gamma) \wedge (n p \gamma^2)= o(\log n)$.
The concentration bounds on the operator norm $\twonorm{YY^T - \E 
  YY^T}$ imply that, up to a constant factor,  the same bounds also 
hold for $\twonorm{A - \E A}$, in view of 
Lemma~\ref{lemma::TLbounds} and~\eqref{eq::Adev}. This is a desirable
property of the SDP \eqref{eq::sdpmain} for its global analysis.
Denote by $C_0$ the $\psi_2$-constant of $\Z_i \in \R^{p}$, cf.~\eqref{eq::covZ1}.
When the independent sub-gaussian random vectors $\Z_i$ in 
\eqref{eq::model} are isotropic, we use the following notion of signal-to-noise ratio:
\ben
\label{eq::SNR}
\text{\bf SNR isotropic } \quad s^2 = ({p \gamma}/{C_0^2})
\wedge ({n p \gamma^2}/{C_0^4}).
\een
This quantity appears in~\cite{GV19}
explicitly and implicitly in~\cite{CHRZ07,Zhou06}, and is compatible
with the total sample size ($np$) lower bound in~\eqref{eq::kilo}, to ensure
partial recovery of the clusters when $s^2 =o(\log n)$. 

In the present work: (a) we will present a transparent and unified 
global and local analysis framework for the semidefinite  programming
optimization problem~\eqref{eq::sdpmain}; and (b) we prove that
the error decays exponentially with respect to the signal-to-noise
ratio (SNR) in Theorem~\ref{thm::exprate}, even for the sub-gaussian design with
dependent features so long as the lower bounds in~\eqref{eq::NKlower}
hold. In particular, the implication of such an exponentially decaying error 
bound is: when  $s^2 =\Omega(\log n)$, perfect recovery of the cluster 
structure is  accomplished.  See~\cite{Royer17,FC18,GV19} and references therein.
Although the result in Theorem~\ref{thm::exprate}
is in the same spirit as that in~\cite{GV19}, we prove these error
bounds for the SDP~\eqref{eq::sdpmain}, which is motivated by the 
graph partition problem~\eqref{eq::graphcut}, while they establish 
such error bounds for the semidefinite relaxation based on the $k$-means 
criterion~\eqref{eq::SSE} directly following~\cite{PW07}; 
See Section~\ref{sec::variations}, cf.~\eqref{eq::relax16} and 
\eqref{eq::relaxadjust}.

In other words, although our general result in
Theorem~\ref{thm::exprate} coincides with that of~\cite{GV19} for
$k=2$, we prove it for the much simpler SDP1 (and SDP) with the only
set of constraints being $\{Z_{ii} = 1, \forall i \in[n]\}$ besides
$Z$ being positive semidefinite. This improvement may be of theoretical interest.
In particular, unlike~\cite{GV19}, we do not enforce entries of $Z$ in
SDP and SDP1 to be nonnegative, which enables faster computation.
Moreover, we do not need an explicit de-biasing step so 
long as assumption (A2) holds, showing that SDP has an inherent 
tolerance on the variance discrepancy between the two component 
distributions, at least up to a certain threshold as specified 
in~\eqref{eq::bias}.
Finally, through the refined local analysis in the present work, we gained
new insights that for the balanced case  ($n_1 = n_2  = n/2$), the bias
may have been substantially reduced due to our  centering and
adjustment steps in the construction of $A$ as
in~\eqref{eq::defineAintro} and an OracleSDP, cf.~\eqref{eq::sdpB},
and hence we expect that assumption (A2) can be substantially relaxed
for  balanced partitions.

\noindent{\bf Construction of a reference matrix.}
Let $Y$ be a centered matrix with row vectors $Y_1, \ldots, Y_n$
as in~\eqref{eq::defineY}. The SDP1 estimator as
in~\eqref{eq::hatZintro} crucially exploits the geometric properties
of the two mean vectors in matrix $Y$, resulting in a natural choice
of $R$. By linearity of expectation,
\ben 
\nonumber
\E \hat{\mu}_n & := & w_1 \mu^{(1)}+
w_2 \mu^{(2)}, \quad \text{where} \quad w_i = \abs{\MC_i} /n, i =1, 2,
\\
\label{eq::EYpre}
\E (Y_i) & := &
\E X_i - \E \hat \mu_n =
\left\{
  \begin{array}{rl} w_2 (\mu^{(1)} - \mu^{(2)})  & \text{ if }  \; i \in \MC_1; \\
    w_1 (\mu^{(2)} - \mu^{(1)}) &  \text{ if }  \; i \in \MC_2. 
  \end{array}\right. 
\een
\begin{definition}{\textnormal{ \bf (The reference matrix)}}
\label{def::reference}
W.l.o.g., we assume that the first $n_1$ rows of 
$X$ belong to $\MC_1$ and the rest belong to $\MC_2$. 
Denote by $n_1 = w_1 n$ and $n_2 = w_2 n$. For $Y$ as defined in
\eqref{eq::defineY},  we construct
\ben
 R:= \E(Y) \E(Y)^T
 &= &
\label{eq::Rtilt}
p \gamma 
\left[
\begin{array}{cc}
 w_2^2 E_{n_1} &- w_1 w_2E_{n_1 \times n_2} \\
 - w_1 w_2 E_{n_2 \times n_1} & w_1^2 E_{n_2} 
\end{array}
\right].
\een
\end{definition}
Recall the overall goal of convex relaxation is to:
(a) estimate the solution of the 
integer quadratic problem~\eqref{eq::quadratic} with an 
appropriately chosen reference matrix $R$ such that 
solving the integer quadratic problem~\eqref{eq::quadratic}
will recover the cluster {\it exactly}; and (b)
choose the convex set $\M^{+}_{G}$ (resp. $\M_{\opt}$, cf. ~\eqref{eq::moptintro}) so that the semidefinite relaxation of the static problem~\eqref{eq::quadratic} is tight.
This means that when we replace $A$  (resp. $\tilde{A} := YY^T - \lambda E_n$)
with $R$ in  SDP~\eqref{eq::sdpmain} (resp. SDP1~\eqref{eq::hatZintro})
as done in Lemma~\ref{lemma::ZRnormintro}, we obtain a solution $Z^{*}$ which can
then be used to recover the clusters exactly.
Here, $Z^* =u_2 u_2^T$ will maximize the reference objective function
$\ip{R, Z}$ among all $Z \in [-1, 1]^{n   \times n}$, and naturally
among all $Z \in \M_{\opt} \subset [-1, 1]^{n \times n}$,
given that $Z^{*} \in \M_{\opt}$. 
\begin{lemma}{\textnormal{~\citep{Zhou23a}}}
 \label{lemma::ZRnormintro}
  Let  $\M_{\opt} \subseteq \M_G^{+} \subset [-1, 1]^{n \times n}$
  be as defined in~\eqref{eq::moptintro}:
\ben 
\label{eq::moptintro}
\M_{\opt}
& = & \{Z: Z \succeq 0,  \diag(Z) = I_n\} \subset \M^{+}_{G}. 
\een
Then for $R$ as in Definition~\ref{def::reference} and $u_2 \in \{-1,
1\}^n$ as in~\eqref{eq::u2},
$$Z^{*} = \argmax_{Z \in \M_{\opt}}\ip{R, Z} = u_2 u_2^T.$$
\end{lemma}

\noindent{\bf Construction of an OracleSDP.}
A straightforward calculation leads to the expression of the 
reference matrix $R:= \E(Y) \E(Y)^T$ as in Definition~\ref{def::reference}. 
However, unlike the settings of~\cite{GV15}, 
$\E A \not= R$, resulting in a large bias. 
It is crucial to bridge the gap between $YY^T$ and the 
reference matrix $R = \E(Y) \E(Y)^T$. 
A remedy was proposed in~\cite{Zhou23a} to
transform~\eqref{eq::sdpmain} into an equivalent {\bf OracleSDP} formulation:
\ben 
\label{eq::sdpB}
{\bf (OracleSDP) }    && 
\text{maximize} \quad \ip{B, Z} \quad \text{ subject to} \quad Z \in
\M_{\opt} \quad \text{where} \\
\label{eq::defineBintro}
&& B  :=  A - \E \tau I_n \quad \text{ for } \quad \tau =\inv{n}
\sum_{i=1}^n \ip{Y_i, Y_i},
\een
where the adjustment term $\E \tau I_n$ plays no role in 
optimization, since the extra trace term $\propto  \ip{I_n, Z} =
\tr(Z)$ is a constant function of $Z$ on the feasible set $\M_{\opt}$.
In other words, optimizing the original SDP~\eqref{eq::sdpmain}
over the larger constraint set $\M^{+}_{G}$ is equivalent to 
maximizing $\ip{B, Z}$ over $Z \in \M_{\opt}$ as in \eqref{eq::sdpB}.
We emphasize that our algorithm solves 
SDP1~\eqref{eq::hatZintro} rather than the oracle 
SDP~\eqref{eq::sdpB}. However, formulating the OracleSDP helps us with the global and local analyses, in controlling 
the operator norm of the bias term, namely, $\twonorm{\E B - R}$ and a
related quantity given in Lemma~\ref{lemma::optlocal}.

\noindent{\bf A unified framework for the local and global analyses.}
Given this OracleSDP formulation, the bias analysis on $\E B - R$ as in  
Lemma~\ref{lemma::EBRtilt}, and the concentration of measure 
bounds on $YY^T -\E YY^T$ have already enabled simultaneous analyses
of the SDP and the closely related spectral clustering method in~\cite{Zhou23a};
cf. Section~\ref{sec::kmeans}.
In the local analysis in the present paper, cf. Lemma~\ref{lemma::YYlocal},
we will obtain a high probability uniform control over $\abs{\ip{Y Y^T -
    \E(YY)^T,  \hat{Z} - Z^{*} }}$ in a {\em local neighborhood} of $Z^{*}$,
\ben 
\label{eq::signal2}
&& \sup_{\hat{Z} \in \M_{\opt} \cap (Z^{*} + r_1 B_1^{n \times n})}
\abs{\ip{YY^T -\E(YY^T), \hat{Z} -Z^{*}} }\le f(r_1) \quad (W_1), \\
\text{ where} \quad
\nonumber
&& \onenorm{\hat{Z} - Z^{*}} :=
\sum_{i,j} \abs{\hat{Z}_{ij}  - Z_{ij}^{*}} \le r_1, \quad \text{ for
} \hat{Z}, Z^{*} \in \M_{\opt},
\een
and $f(r_1)$ is a function that depends on the $\ell_1$ radius $r_1$,
$r_1 \le 2n(n-1)$. Lemma~\ref{lemma::YYlocal} is one of the main
technical contributions in this paper, which in turn depends on the geometry of the constraint
set $\M_{\opt}$ and  the sharp concentration of
measure bounds on $YY^T-\E YY^T$ in Theorem~\ref{thm::YYaniso}.
This approach also gives rise to the notion of a {\em local analysis},
following~\cite{CCLN21}.
Similar to the global analysis, the performance of
SDP also crucially depends on controlling the bias
term $\E B - R$; however, we now control $\ip{\E B - R, \hat{Z}-Z^{*}}$
uniformly over all $\hat{Z} \in \M_{\opt}$ in
Lemma~\ref{lemma::optlocal}, rather than the operator
norm $\twonorm{\E B - R}$. The bias and variance tradeoffs are
described in the sequel.
Combining these results leads to an exponentially decaying error bound 
with respect to the SNR $s^2$ to be presented in
Theorem~\ref{thm::exprate}.

We need to introduce some notation.
Recall for a random variable $X$, the $\psi_2$-norm of $X$, 
  denoted by $\norm{X}_{\psi_2}$, is $\norm{X}_{\psi_2} = \inf\{t >
  0\; : \; \E \exp(X^2/t^2) \le 2 \}$.
We use $\Z = (z_{ij})$ to denote the mean-zero random 
matrix with independent, mean-zero, sub-gaussian row vectors $\Z_1,
\ldots, \Z_n \in \R^p$ as considered in~\eqref{eq::model}, where for a
constant $C_0$ and $\cov(\Z_j) := \E (\Z_j \Z_j^T)$,
\ben
\label{eq::covZ1}
\forall j=1, \ldots, n, \; 
\norm{\ip{\Z_j, x}}_{\psi_2} & \le & C_0 \norm{\ip{\Z_j, x}}_{L_2}
\text{ for any } x \in  \R^{p},  \\
\label{eq::covZ2}
\text{ where } \; \norm{\ip{\Z_j, x}}^2_{L_2} &:=& x^T \E (\Z_j
\Z_j^T) x = x^T \cov(\Z_j) x.
\een
Throughout this paper, we use $V_i$ to denote the trace of
covariance $\cov(\Z_j)$, for all $j \in \MC_i$:
\ben
\label{eq::V2intro}
V_1 & = & \E \ip{\Z_j, \Z_j}  \quad \forall j \in \MC_1 \;\; \text{and}\;\; 
V_2  =  \E \ip{\Z_j, \Z_j} \quad \forall j \in \MC_2.
\een
\noindent{\bf (A1)}
Let $\Z = X - \E X$. Let $\Z_i = X_i - \E X_i, i \in [n]$ be 
  independent, mean-zero, sub-gaussian random vectors with independent
  coordinates such that for all $i,j$, $\norm{X_{ij} - \E 
    X_{ij}}_{\psi_2} \le C_0$. \\
  \noindent{\bf (A2)}
The two distributions have bounded discrepancy in their variance
profiles:
\ben
\label{eq::Varprofile}
\abs{V_1 - V_2}  & \le & \xi 
n p \gamma/3 \; \; \text{ for some } \quad w_{\min}^2/8 > 2 \xi =
\Omega(1/n_{\min}).
\een
Although (A1) assumes that the random matrix $\Z$ has independent
sub-gaussian entries, matching the separation (and SNR) condition~\eqref{eq::kilo},
we emphasize that the conclusions as stated in  
Theorems~\ref{thm::SDPmain} and~\ref{thm::exprate} hold for the general two-group model as considered in Lemma~\ref{lemma::twogroup} so long as (A2) holds, upon 
adjusting the lower bounds in~\eqref{eq::kilo}.
In particular, we allow each population to have distinct covariance 
structures, with diagonal matrices as special cases. 
\subsection{Prior results}
 For completeness, we first state in Theorem~\ref{thm::SDPmain} the 
 main result under assumptions (A1) and (A2),
 using the global analysis \citep{Zhou23a,GV15}, in order to set the
 context for local analysis in this paper.
 For a matrix $B= (b_{ij})$ of size $n \times n$, 
let $\mvec{B}$ be formed by concatenating columns of matrix 
$B$ into a vector of size $n^2$; 
we use $\onenorm{B} =  \sum_{i,j=1}^n |b_{ij}|$ to denote 
the $\ell_1$ norm of $\mvec{B}$
and $\fnorm{B} = (\sum_{i, j} b_{ij}^2)^{1/2}$ to denote 
the $\ell_2$ norm of $\mvec{B}$. 
\begin{theorem}{\textnormal{\citep{Zhou23a}}}
  \label{thm::SDPmain}
  Let $\MC_j \subset [n]$ denote the group membership, with $\abs{\MC_j} =
  n_j$ and $\sum_{j} n_j = n$.
 Suppose that for $j \in \MC_i$, $\E X_j = \mu^{(i)}$, where $i =1, 2$.
  Let $\hat{Z}$ be a solution of SDP1 (and SDP).
  Suppose that (A1) and (A2) hold, and for some  absolute constants
  $C, C'$,
 \ben 
  \label{eq::kilo}
 p \gamma = \Delta^2 \ge \frac{C' C_0^2}{\xi^2}  \; \;  \text{ and } \; 
 p n \ge \frac{C C_0^4}{\xi^2 \gamma^2}, \text{ where $\xi$ is the
   same as in~\eqref{eq::Varprofile}.}
 \een
Then with probability at least $1-2 \exp(-c n)$, we have for ${u}_2$ as in \eqref{eq::u2},
\ben
\label{eq::hatZFnorm}
\onenorm{\hat{Z} - u_2 u_2^T} /n^2 & = &   \delta \le {2 K_G
  \xi}/{w_{\min}^2}   \; \; \text{and} \; \;
\fnorm{\hat{Z} - u_2 u_2^T}/ n^2\le {4 K_G \xi}/{w_{\min}^2}.
\een
\end{theorem}
The parameter $\xi$ in (A2)  is understood to be set as a fraction of
the parameter $\delta$ appearing in Theorem~\ref{thm::SDPmain}, 
which in turn is assumed to be lower bounded: $1> \delta =
\Omega(1/n)$, since we focus on the partial recovery of the
group membership vector $u_2$ in the present work.
In particular, $\xi^2$ can be chosen conservatively to be inversely
proportional to $s^2$; cf. Section~\ref{sec::theory}.
This choice enables the bias and variance tradeoffs as elaborated in
Section~\ref{sec::theory}.

\subsection{Related work}
\label{sec::related}
In the theoretical computer science literature,
earlier work focused on learning from mixture of well-separated
Gaussians (component distributions), where one aims to 
classify each sample according to which component distribution it 
comes from; See for example~\cite{DS00,AK01,VW02,AM05,KMV10,KK10}.
In earlier works~\cite{DS00,AK01}, the separation requirement 
depends on the number of dimensions of each distribution; this has recently 
been reduced to be independent of $p$, the dimensionality of 
the distribution for certain classes of
distributions~\cite{AM05,KSV05}.
While our aim is different from those results, where $n > p$ is almost
universal and we focus on cases $p > n$ (a.k.a. high dimensional
setting), we do have one common axis for comparison, the $\ell_2$-distance between any two centers of
the distributions as stated in~\eqref{eq::NKlower}
(or~\eqref{eq::kilo}), which is essentially optimal.
Results in~\cite{Zhou06,CHRZ07} were among the first such 
results towards understanding rigorously and intuitively why their 
proposed algorithms and previous methods in~\cite{PattersonEtAl,PriceEtAl} work with low sample settings
when $p \gg n$ and $n p$ satisfies~\eqref{eq::kilo}.
However, such results were only known to exist for balanced max-cut 
algorithms~\citep{Zhou06,CHRZ07}, and hence these were structural as no polynomial 
time algorithms were given for finding the max-cut.

The main contribution of the present work and~\cite{Zhou23a} is: we use the 
proposed SDP~\eqref{eq::sdpmain} and the related spectral algorithms 
to find the partition, and prove quantitatively tighter bounds than those 
in~\cite{Zhou06} and~\cite{CHRZ07} by removing these logarithmic factors 
entirely. Recently, these barriers have also been broken down by a sequence of 
work~\cite{Royer17,FC18,GV19}.
For example,~\cite{FC18,FC21} have also established such bounds, but
they focus on balanced clusters and require an extra 
$\sqrt{\log n}$ factor in \eqref{eq::FZ18} in the second component: 
\ben 
\label{eq::FZ18}
\text{In~\cite{FC18}, cf. eq.(8): } \quad \Delta^2 = p \gamma &= & 
  \Omega\left(1 +  \sqrt{{p \log n}/{n}}  \right) \quad  \text{or } \\
  \label{eq::FZ21}
    \text{In~\cite{FC21},  cf. eq.(13): } \quad \Delta^2 = p \gamma &= & \Omega\left((1 \vee {p}/{n})  + \sqrt{{p \log n}/{n}}  \right). 
    \een
As a result, in~\eqref{eq::FZ21}, a lower bound on the sample size is
imposed: $n \ge 1/\gamma$ in case $p > n$, and moreover, the size of
the matrix $n p =\Omega(\log n/\gamma^2)$,
similar to the bounds in~\cite{BCFZ09}; cf. Theorem 1.2 therein.
Hence, these earlier results still need the SNR to be at the order of $s^2 =
\Omega(\log n)$. We compare with~\cite{GV19} in 
Section~\ref{sec::kmeans}. 
Importantly, because of these new concentration of measure bounds,
our theory works for the small sample and low SNR ($n < p$ and $s^2 =
o(\log n)$) cases for both SDP (in its local and global analyses) and
the spectral algorithm under (A2), which is perhaps the more
surprising case. For example, previously, the spectral algorithms
in~\cite{BCFZ09} partition sample points based on the top few
eigenvectors of the gram matrix $XX^T$, following an idea that goes
back at least to~\cite{Fie1973}.
In~\cite{BCFZ09}, the two parameters
$n, p$ are assumed to be roughly at the same order, hence not allowing a full 
range of tradeoffs between the two dimensions as considered in the 
present work. Such a lower bound on $n$ was deemed to be unnecessary given the 
empirical evidence~\cite{BCFZ09}. 
In contrast, the spectral analysis in~\cite{Zhou23a} uses the leading
eigenvector of  the gram matrix $YY^T$, based on centered data,
which will directly improve the results in~\cite{BCFZ09} as we can now
remove the lower bound on $n$. 
We also refer to~\cite{KMV10,RCY11,GV15,LZ16,Abbe16,BWY17,CY18,BMVV+18,GV19,LLLS+20,FC21,LZZ21,AFW22,nda22}
and references therein for related work on the Stochastic Block Models 
(SBM), mixture of (sub)Gaussians and clustering in more general metric
spaces.

\noindent{\bf Organization.}
The rest of the paper is organized as follows.
In Section~\ref{sec::theory},  we present the sub-gaussian mixture
models and our main theoretical results.
In Section~\ref{sec::kmeans}, we highlight the connections and key 
differences between our approach and convex relaxation algorithms based on 
the $k$-means criterion, and other related work. 
We present preliminary results for the local analysis (for proving 
Theorem~\ref{thm::exprate}) in Section~\ref{sec::exprate}.
We prove Theorem~\ref{thm::exprate} in
Section~\ref{sec::proofexprate}.
We provide a proof sketch for the key lemmas for
Theorem~\ref{thm::exprate} in Section~\ref{sec::keylemmas}.
 We conclude in Section~\ref{sec::conclude}.
We place all additional technical proofs
in the supplementary  material.

\noindent{\bf Notation.}
\label{sec::notation}
Let $\Ball_2^n$ and $\Sp^{n-1}$ be the unit Euclidean ball and the
unit sphere in $\R^n$ respectively.
Let $e_1, \ldots, e_n$ be the canonical basis of $\R^n$.
For a set $J \subset \{1, \ldots, n\}$, denote
$E_J = \Span\{e_j: j \in J\}$.
For a vector $v \in \R^n$, we use $v_{J}$ to denote the subvector
$(v_j)_{j \in J}$.
For a vector $x$, $\norm{x}_{\infty} := \max_{j} \abs{x_j}$, $\onenorm{x} := \sum_{j} \abs{x_j}$, 
and  $\twonorm{x} := \sqrt{\sum_{j} x^2_j}$; $\diag(x)$ denotes the diagonal matrix whose 
main diagonal entries are the entries of $x$.
For a matrix $B \in \R^{n \times n}$, $\tr(B) = \sum_{i=1}^n B_{ii}$.
Denote by $B_1^{n \times n} = \{(b_{ij}) \in \R^{n \times n}  : 
\sum_{i} \sum_{j} \abs{b_{ij}} \le 1\}$ the $\ell_1$ unit ball over
$\mvec{B}$. 
Let $\diag(A)$ and $\offd(A)$ be the diagonal and the off-diagonal 
part of matrix $A$ respectively. 
For a matrix $A$, let $\norm{A}_{\infty} = \max_{i} \sum_{j=1}^n
|a_{ij}|$ denote the maximum absolute row sum;
Let $\norm{A}_{\max} = \max_{i,j} |a_{ij}|$ denote the component-wise
max norm.
For two numbers $a, b$, $a \wedge b := \min(a, b)$, and 
$a \vee b := \max(a, b)$.
We write $a \asymp b$ if $ca \le b \le Ca$ for some positive absolute
constants $c,C$ which are independent of $n, p$, and $\gamma$.
We write $f = O(h)$ or $f \ll h$ if $\abs{f} \le C h$ for some absolute constant
$C< \infty$ and $f=\Omega(h)$ or $f \gg h$ if $h=O(f)$.
We write $f = o(h)$ if $f/h \to 0$ as $n \to \infty$, where the
parameter $n$ will be the size of the matrix under consideration.
In this paper, $C, C_1, C_2, C_4, c, c', c_1$, etc,
denote various absolute positive constants which may change line by line.

\section{Theory}
\label{sec::theory}
In this section, we first state the main result derived from the local
analysis, namely, Theorem~\ref{thm::exprate}, where we present an error bound 
\eqref{eq::hatZFnorm2} that decays 
exponentially in the SNR parameter $s^2$ as defined in~\eqref{eq::SNR2}. 
We then present Corollary~\ref{coro::misexp}, followed by discussions.
We prove Theorem~\ref{thm::exprate} in Section~\ref{sec::proofexprate}.
First, we have two definitions.
\begin{definition}
A random vector $W \in \R^m$ is called sub-gaussian if the one-dimensional marginals $\ip{W, h}$ are sub-gaussian random 
variables for all $h \in \R^{m}$:
(1) $W$ is called isotropic if for every $h \in \Sp^{m-1}$,
$\E \abs{\ip{W, h}}^2= 1$;
(2)  $W$ is $\psi_2$ with a constant $C_0$ if for every $h \in
\Sp^{m-1}$, $\norm{\ip{W, h}}_{\psi_2} \le C_0$.  The sub-gaussian norm of $W \in \R^m$ is denoted by
  \ben
  \label{eq::Wpsi}
  \norm{W}_{\psi_2}  := \sup_{h \in \Sp^{m-1}} \norm{\ip{W, h}}_{\psi_2}.
  \een
\end{definition}

\begin{definition}{\textnormal{\bf (Data generative process.)}}
  \label{def::WH}
   Let $H_1, \ldots, H_n$ be deterministic $p \times m$ matrices, 
 where we assume that $m \ge p$.
  Suppose that random matrix $\BW =(w_{jk})\in \R^{n \times m}$ has
  $W_1, \ldots, W_n \in   \R^{m}$ as independent row vectors, where
  $W_j$, for each $j$, is an isotropic sub-gaussian random vector with independent entries satisfying 
\ben 
\label{eq::Wpsi2}
\forall j \in [n], \quad \cov(W_j) :=\E (W_j W_j^T) = I_m, \; \; 
\E[w_{jk}] = 0, \quad \text{ and } \quad \max_{jk} \norm{w_{jk}}_{\psi_2} \le C_0, \quad \forall k.
\een
Suppose that we have for row vectors $\Z_1, \ldots, \Z_n \in \R^{p}$
of the noise matrix $\Z \in \R^{n \times p}$,   
 \bens 
\forall j =1, \ldots, n,\quad  \Z^T_j  =  W^T_j H_j^T, \quad \text{
  where} \quad H_j \in \R^{p \times    m}, \; \; 0< \twonorm{H_j} < \infty,
\eens
and $H_j$'s are allowed to repeat, for example, across rows from the
same cluster $\MC_i$ for some $i =1, 2$. 
Throughout this paper, we assume that $m \ge p$ to simplify our 
exposition, although this is not necessary.
\end{definition}

\noindent{\bf Signal-to-noise ratios.}
This notion of SNR \eqref{eq::SNR}
can be properly adjusted when coordinates in $\Z_i$ are
dependent in view of~\eqref{eq::covZ1} and \eqref{eq::covZ2}:
\ben
\label{eq::SNR2}
\text{\bf SNR anisotropic: } \quad s^2
=\frac{\Delta^2}
{C_0^2 \max_{j}\twonorm{\cov(\Z_j)}} \wedge \frac{n p \gamma^2}
{C_0^4 \max_{j}\twonorm{ \cov(\Z_j)}^2}.
\een
\begin{theorem}
  \label{thm::exprate}
  Suppose that for $j \in \MC_i \subset [n]$, $\E X_j = \mu^{(i)}$, where $i =1,
  2$.   Let $\hat{Z}$ be a solution of SDP1.
  Let $s^2$ be as defined in~\eqref{eq::SNR2}. Let $C_0$ be as defined
  in \eqref{eq::Wpsi2}.
Suppose the noise matrix $\Z = X-\E(X)$ is generated according to Definition~\ref{def::WH}:
\bens
\forall  j \in \MC_i, \quad \Z_j = H_i W_j,  \quad  \text{ where}
\quad H_i \in \R^{p \times m} \quad \text{ is deterministic and }
\quad 0< \twonorm{H_i} < \infty,
\eens
for $i \in \{1, 2\}$. 
Suppose that  for some absolute constants $C, C_{1}$ and some $0< \xi
\le w_{\min}^2/16$,
\ben
\label{eq::NKlower}
p \gamma \ge \frac{C C_0^2 \max_{j}\twonorm{\cov(\Z_j)} }{\xi^2}
\quad \text{ and } \quad np \ge \frac{C_{1} C_0^4
  \max_{j}\twonorm{\cov(\Z_j)}^2}{\gamma^2 \xi^2}.
\een
  Suppose (A2) holds, where the parameter $\xi$ in~\eqref{eq::Varprofile} is understood to  be the same as the parameter $\xi$ in \eqref{eq::NKlower}.
Then with probability at least $1-2 \exp(-c_1 n) - c_2/n^2$, for some
absolute constants $c_0, c_1, c_2$,
\ben 
\label{eq::hatZFnorm2}
\onenorm{\hat{Z} - {u}_2 {u}_2^T}/n^2 \le \exp\big(-c_0 s^2  w_{\min}^4\big).
\een
\end{theorem}

\begin{corollary}{\textnormal{{\bf (Exponential decay in $s^2$)}}}
  \label{coro::misexp}
    Let $\hat{x}$ denote the eigenvector of $\hat{Z}$
  corresponding to the largest eigenvalue, with $\twonorm{\hat{x}}
  = \sqrt{n}$.
  Denote by $\theta_{\SDP} = \angle(\hat{x}, u_2)$, the angle
  between $\hat{x}$ and $u_2$,
  where recall $u_{2j} = 1$ if $j 
  \in \MC_1$ and $u_{2j} = -1$ if  $j \in \MC_2$.
  In the settings of Theorem~\ref{thm::exprate},
  with probability at least $1-2 \exp(-c n) - 2/n^2$,
for some absolute constants $c, c_0, c_1$,
\ben
\nonumber
 \sin(\theta_{\SDP}) & \le & 2 \twonorm{\hat{Z} - u_2 u_2^T}/{n}  \le \exp(- c_1 s^2 w_{\min}^4) \; \;  \text{ and } \\
 \label{eq::eigenconv}
\min_{\alpha = \pm 1}  \twonorm{(\alpha \hat{x} - u_2)/\sqrt{n}}
  & \le &
  {2^{3/2} \twonorm{\hat{Z} - u_2 u_2^T}}/{n}  \le
 4 \exp(-c_0 s^2 w_{\min}^4/2).
 \een 
\end{corollary}
Corollary \ref{coro::misexp} follows from the Davis-Kahan Theorem and 
Theorem~\ref{thm::exprate}, and is given in~\cite{Zhou23a}; cf. Proof
of Corollary 2.8 therein.
Our local and global analyses also show the surprising 
result that Theorems~\ref{thm::SDPmain} and~\ref{thm::exprate} do not 
depend on the clusters being balanced,  nor do they require identical 
variance or covariance profiles, so long as (A2) holds.

\noindent{\bf The bias and variance tradeoffs.}
Roughly speaking, to ensure that $\twonorm{\E B - R}$ is bounded, we require
\ben
\label{eq::bias}
(A2') \quad \abs{V_1 - V_2} \le  w =
O\left(n \Delta  \vee \sqrt{n p} \right), \; \text{where}  \; \; \Delta =
\sqrt{p \gamma}.
\een
Notice that $\abs{V_1 - V_2} = O(p)$ by definition, and hence
\eqref{eq::bias} holds trivially in the large sample setting where $n > p$.
Given a fixed average quality parameter $\gamma$, the tolerance on
$\abs{V_1 - V_2}$ depends on the sample size $n$, the total data size
$D:= np$, and the separation parameter $\Delta$. In particular, the tolerance
parameter $w$ is chosen to be at the same order as the upper bound we
obtain on $\twonorm{YY^T -  \E YY^T}$ in Theorem~\ref{thm::YYaniso}.
In the supplementary Theorem~\ref{thm::reading}, we show that
the operator norm for $B - R$ is controlled essentially at the same order as
$\twonorm{YY^T - \E (Y Y^T)}$. To see this, first, we have for  $A$ as in~\eqref{eq::sdpmain},
\ben
\label{eq::YYop}
\twonorm{B -\E B} & \le & \twonorm{Y Y^T - \E Y Y^T}  + \abs{\lambda- \E
  \lambda}\twonorm{E_n - I_n},\quad \text{since } \\
\label{eq::Adev}
B - \E B & := &  A - \E A = Y Y^T - \E Y Y^T  - (\lambda- \E \lambda)(E_n - I_n).
\een
Combining Lemmas~\ref{lemma::EBRtilt},~\ref{lemma::TLbounds},
Theorem~\ref{thm::YYaniso}, with~\eqref{eq::YYop},
we have with probability at least $1 - 2\exp(-c n)$,
\ben
\label{eq::Bdev2}
\twonorm{B - R}
& \le & \twonorm{B- \E B} + \twonorm{\E B -R}
\le 2 \twonorm{YY^T -  \E (Y Y^T)}+ \twonorm{\E B  - R} \le  \xi n p \gamma.
\een
First, we state Lemma~\ref{lemma::twogroup}, which characterizes the
two-group design matrix covariance structures as considered in
Theorems~\ref{thm::SDPmain} and~\ref{thm::exprate}, with isotropic
design as the special case.
\begin{lemma}\textnormal{\bf (Two-group sub-gaussian mixture model)}
\label{lemma::twogroup}
Denote by $X$ the two-group design matrix as considered in
\eqref{eq::Xmean}.
Let $W_1, \ldots, W_n \in \R^{m}$ be independent, mean-zero, 
isotropic, sub-gaussian random vectors satisfying~\eqref{eq::Wpsi2}. 
Let $\Z_j = X_j - \E X_j = H_i W_j$, for all $j \in \MC_i$, where $H_i
\in \R^{p \times m}$ is deterministic and $0<
\twonorm{H_i} < \infty$, for each $i \in \{1, 2\}$.
Then $\Z_1, \ldots, \Z_n$ are independent sub-gaussian random vectors
with $\cov(\Z_i)$ satisfying \eqref{eq::covZ1} and \eqref{eq::covZ2}, where
\ben 
\label{eq::rowcov}
\forall j \in \MC_i, \quad
\cov(\Z_j) := \E (\Z_j \Z_j^T)  = \E (H_i W_j W_j^T H_i^T)  =  H_i  H_i^T \text{ and } \;
V_i := \fnorm{H_i}^2 = \tr(\Sigma_i).
\een
\end{lemma}
\begin{lemma}{\textnormal{\cite{Zhou23a}}}
  \label{lemma::EBRtilt}
  Let $X$ be as in Lemma~\ref{lemma::twogroup}. 
  Suppose (A2) holds.
  Suppose that $\xi \ge \inv{2n} (4 \vee \inv{w_{\min}})$ and $n \ge 4$.
  Then we have $\twonorm{\E B - R}  \le \frac{2}{3} \xi  n p \gamma.$
Finally, when $V_1 = V_2$, we have $\twonorm{\E B - R} \le  p \gamma/3$.
\end{lemma}
\begin{theorem}{\textnormal{{\bf (Anisotropic design matrix)}~\citep{Zhou23a}}}
\label{thm::YYaniso}
Let $X$ be as in Lemma~\ref{lemma::twogroup}. 
Let $Y$ be as in \eqref{eq::defineY}.
Let $\mu =\frac{\mu^{(1)}-\mu^{(2)}}{\sqrt{p \gamma}} \in \Sp^{p-1}$.
Then with probability at least $1 - 2\exp(-c_8 n)$, for $C_0$ as defined in \eqref{eq::Wpsi2},
  \bens
\twonorm{YY^T - \E (Y Y^T)}
& \le &
C_3  (C_0 \max_{i} \twonorm{H_i^T \mu}) n \sqrt{p \gamma} + C_4  (C_0
\max_{i} \twonorm{H_i})^2 (\sqrt{p n} \vee n) \le \inv{6}  \xi n p \gamma.
\eens
\end{theorem}
See the supplementary Section~\ref{sec::oracleSDP} for a proof sketch 
of Theorem~\ref{thm::YYaniso} and the proof of
Lemma~\ref{lemma::EBRtilt}, which we include for self-containment.
The error bound in \eqref{eq::Bdev2} was crucial in the global analysis 
leading to Theorem~\ref{thm::SDPmain}, as well as the simultaneous 
analysis of the spectral clustering method;
cf. Section~\ref{sec::kmeans}, Variation 2.
Indeed, our previous and current results show that even when $n$ is
small, by increasing $p$ such that the total sample size satisfies the second condition
in~\eqref{eq::kilo} and~\eqref{eq::NKlower}, we ensure
partial recovery of cluster structures using the
SDP~\eqref{eq::sdpmain} or the spectral algorithm as described in~\cite{Zhou23a}.
On the other hand, once there are enough features such that the
distance between the two centers are bounded below by a constant, adding more samples (individuals) to the clusters will reduce the number of features we
need to do partial recovery.
Now~\eqref{eq::NKlower} implies that 
\ben
\label{eq::invxi}
&&
256/ w_{\min}^4 \le \inv{\xi^2}  \le
\frac {p \gamma}{C C_0^2 \max_{j}\twonorm{\cov(\Z_j)} } \wedge
\frac {n p \gamma^2}{C_{1} C_0^4 \max_{j}\twonorm{\cov(\Z_j)}^2}
\asymp s^2,  \quad \text{ and hence }\\
\nonumber
&& \xi n p \gamma/2 \ge 
C_0^2 \max_{i} \twonorm{\cov(\Z_i)}
\left(\frac{n}{\xi} \vee \sqrt{np} \right) =: w, \quad \text{ given that }
\quad \inv{\xi} \le \frac{\sqrt{p \gamma}}{C C_0 
  \max_{j}\twonorm{\cov(\Z_j)}^{1/2}},
\een
resulting in \eqref{eq::bias} being a sufficient condition for (A2) to
hold. Given fixed $p$ and $\gamma$, increasing the sample size $n$
also leads to a {\em higher tolerance} on the variance discrepancy $\abs{V_1 - V_2}$
between the two component distributions as shown in~\eqref{eq::bias}.
To control the misclassification error using the
global analysis, the parameters $(\delta, \xi)$ in
Theorem~\ref{thm::SDPmain} can be chosen to
satisfy the following relations:
\ben
\label{eq::trend}
\xi^2 \asymp \frac{C_0^2}{p \gamma} \vee \frac{C_0^4}{ n p \gamma^2} = \frac{1}{s^2}
\quad \text{ and } \quad \delta =  \frac{2 K_G \xi}{w_{\min}^2},
\een
  in view of~\eqref{eq::kilo} and~\eqref{eq::hatZFnorm}.
Thus we obtain in Theorem~\ref{thm::SDPmain} that
the misclassification error is inversely proportional to 
the square root of the SNR parameter $s^2$ as in~\eqref{eq::SNR}
in view of \eqref{eq::trend},
 while in Theorem~\ref{thm::exprate} and Corollary~\ref{coro::misexp},
 the error decays exponentially in the SNR parameter $s^2$ as defined
 in~\eqref{eq::SNR2}. 

\section{The $k$-means criterion of a partition}
\label{sec::kmeans}
We now discuss the $k$-means criterion and its
semidefinite relaxations.
Denote by $X \in \R^{n \times p}$ the data matrix with row vectors 
$X_i$ as in~\eqref{eq::SSE} (see also~\eqref{eq::model}).
The $k$-means criterion of a partition $\MC = \{\MC_1, \ldots, \MC_k\}$ of sample points $\{1,\ldots, n\}$ is based on the total sum-of-squared Euclidean distances from each point
$X_i \in \R^{p}$ to its assigned cluster centroid $\vc_j$, namely,
\ben
\label{eq::SSE}
g(X, \MC, k) := \sum_{j=1}^k \sum_{i \in \MC_j} \twonorm{X_i - \vc_j}^2, \quad \text{
  where} \quad \vc_j := \inv{\abs{\MC_j}} \sum_{\ell \in \MC_{j}}
X_{\ell} \in \R^{p}.
\een
Getting a global solution to~\eqref{eq::SSE} through an integer
programming formulation as in~\cite{PX05} and~\cite{PW07} is NP-hard
and it is NP-hard for $k=2$~\citep{DFK+04,ADHP09}.
The partition $\CC$ can be represented by
a block diagonal matrix of size $n \times n$, defined as: $\forall i, j \in [k]^2, \forall
a, b \in \MC_{i} \times \MC_{j}$,
\bens
B_{ab}^* = 1/{\size{\MC_j}} \; \text{ if } \; i = j \text{ and }
B_{ab}^* = 0  \text{ otherwise. }
\eens
The collection of such matrices can be described by
\ben
\label{eq::cpk}
\mathcal{P}_k =\{B \in \R^{n \times n}: B = B^T,  B \ge 0, B^2 = B, \tr(B) = k,
B \vecone_n = \vecone_n\},
\een
where $B \ge 0$ means that all elements of $B$ are nonnegative.
Minimizing the $k$-means objective $g(X, \MC, k)$ is equivalent to~\citep{ZDGH+02,PX05,PW07}, 
\ben 
\label{eq::relax7}
&& \text{maximize} \quad \ip{X X^T, Z}  \quad \text{ subject to } \quad Z \in 
\cpk.
\een

\noindent{\bf Peng-Wei relaxations.}
Let $\Psi_n$ denote the linear space of real $n$ by $n$ symmetric matrices.
Now consider the semidefinite relaxation of the $k$-means objective~\eqref{eq::relax7},
\ben
\label{eq::relax16}
&& \text{maximize}
\quad \ip{ X X^T, Z}  \quad \text{ subject to }  \quad  Z \in \M_k,
\\
\nonumber
&& \text{  where } \; \M_k = \{Z \in \Psi_n:  Z   \ge 0, Z \succeq 0, \tr(Z) = k, 
Z \vecone_n = \vecone_n \}
\een
as considered in~\cite{PW07}.
  The key differences between this and the SDP~\eqref{eq::sdpmain} are:
  (a) In the convex set~$\M_{\opt}$~\eqref{eq::moptintro}, we do not enforce that all
  entries are nonnegative, namely, $Z_{ij} \ge 0, \forall i, j$.
This allows faster computation; (b) In order to derive concentration of
measure bounds that are sufficiently tight,
we make a natural, yet important data processing step, where we center
the data according to their column mean following
Definition~\ref{def::estimators} before computing $A$ as
in~\eqref{eq::defineAintro};
(c) We do not enforce $Z \vecone_n = \vecone_n$.
Note that for $\M_k$ as in \eqref{eq::relax16},  when $Z_{ij} \ge 0, 
\forall i, j$, row sum $\onenorm{Z_{j, \cdot}} =1, \forall j$ and 
hence  $\twonorm{Z} \le \norm{Z}_{\infty} = 1$. \\
\noindent{\bf Variation 2.}
To speed up computation, one can drop the nonnegative 
constraint on elements of $Z$
leading to  the following semidefinite relaxation~\citep{PW07}: 
\ben
\label{eq::relax17}
\text{maximize} \quad \ip{ X X^T, Z} \quad \text{ subject to }
I_n \succeq Z \succeq 0, \quad \tr(Z) = k,  \text{ and } \quad Z \vecone_n = \vecone_n.
\een
Now
when $Z$ is a feasible solution to~\eqref{eq::relax17}, define
\ben 
\label{eq::Z1}
Z_1 & := & Z - \vecone_n \vecone_n^T/n \; \; \text{ and hence 
}\; \; Z_1 :=  (I-P_1) Z = (I-P_1)Z (I-P_1),
\een
where recall $\vecone_n/\sqrt{n}$ is the unit-norm leading eigenvector
of $Z$.
Now~\cite{PW07} shows that the set of feasible
solutions to~\eqref{eq::relax17} have immediate connections to the SVD
of $YY^T$ through the following reduction:
 \ben 
 \label{eq::relax20}
 &&\text{maximize} \quad \ip{YY^T, Z_1} \quad \text{ subject to }
\quad  I_n \succeq Z_1 \succeq 0, \quad \tr(Z_1) = k-1.
\een
Then the algorithm for solving~\eqref{eq::relax20} and~\eqref{eq::relax17} is given as follows~\citep{PW07}: \\
{\bf Spectral Clustering: (a)}
Using singular value decomposition (SVD) method to compute the first $k-1$
largest eigenvalues of $YY^T$ and their corresponding eigenvectors
$v_1, \ldots, v_{k-1}$; 
\bens
\text{and} \quad \text{\bf (b) Set } \quad
Z_1 =\sum_{j=1}^{k-1} v_j v_j^T; \quad \text{and return} \quad Z =
\inv{n} \vecone_n \vecone_n^T + Z_1 \; \text{ as a solution to \eqref{eq::relax17}}.
\eens
SVD-based algorithms allow even faster computation.
In particular, for $k=2$, we have the optimal solution 
of~\eqref{eq::relax17} being $Z_1 = v_1 v_1^T$, where $v_1$ is the
leading eigenvector of $YY^T$.
The signs of the coefficients of $v_1$ correctly estimate the
partition of the vertices, up to $O(\xi^2 n)$
misclassified  vertices, where we set $\xi^2 \asymp 1/s^2$
in~\cite{Zhou23a}; cf.~\eqref{eq::trend}.
Hence the misclassification error is bounded to be inversely proportional to the SNR
parameter $s^2$; cf. \eqref{eq::trend}.
This should be compared with Theorem~\ref{thm::exprate}
and Corollary~\ref{coro::misexp}, where the misclassification error is  improved to $O( n \exp(-c_0 s^2 w_{\min}^4))$.
 Moreover, one can sort the values of $v_1$ and 
find the {\it nearly optimal } partition according to the $k$-means 
criterion; See  Algorithm 2 and numerical examples in~\cite{Zhou23a}.

\noindent{\bf Discussions.}
\label{sec::variations}
The main issue with the $k$-means relaxation is that the solutions
tend to put sample points into groups of the same sizes, and moreover,
the diagonal matrix $\Gamma$ can cause a bias, where
\bens
\Gamma = (\E [\ip{\Z_i, \Z_j}])_{i,j} = \diag([\tr(\cov(\Z_1)),
\ldots, \tr(\cov(\Z_n)) ]),
\eens
especially when $V_1, V_2$ differ from each other.
The authors of~\cite{Royer17},~\cite{GV19}, and~\cite{BGLR+20} propose
a preliminary estimator of $\Gamma$, denoted by $\hat{\Gamma}$, and consider
\ben
\label{eq::relaxadjust}
&& \hat{Z} \in \arg\max_{Z \in \M_k} \ip{X X^T - \hat{\Gamma}, Z},\quad
\text{ where} \; \M_k \; \text{ is as in~\eqref{eq::relax16} },
\een
instead of the original Peng-Wei SDP relaxation~\eqref{eq::relax16}.
Consequently, besides computation,
another main advantage of our SDP and spectral formulation in~\cite{Zhou23a} is that we do not need to have a separate estimator for $\tr(\Sigma_j)$, where $\Sigma_j, j=1, 2$ denotes the covariance matrices of sub-gaussian random vectors $\Z_j, j \in [n]$, so long as (A2) holds. In some sense, SDP1 has a certain inherent tolerance on 
the variance discrepancy between the two populations in the sense of \eqref{eq::bias}.
When (A2) is violated, we may adopt similar ideas 
in~\cite{Royer17,GV19} to make adjustments to SDP1 to correct the
bias. It is an arguably simpler task to deal with 
$\abs{V_1 - V_2}$ than actual values of $V_1, V_2$. See Remark~\ref{rm::target}.
\section{Proof sketch for Theorem~\ref{thm::exprate}}
\label{sec::exprate}
In this section, we need to state some intermediate results followed by the actual
proof of Theorem~\ref{thm::exprate} in
Section~\ref{sec::proofexprate}.
Let the reference matrix $R$ be constructed as in Definition~\ref{def::reference}.
Lemma \ref{lemma::signalgrad} follows from definitions of $B$ and $R$,
and the lower bound in Lemma~\ref{lemma::onenorm}.
\begin{lemma}{\textnormal{\bf (Elementary Inequality)}}
\label{lemma::signalgrad}
By optimality of $\hat{Z} \in \M_{\opt}$, we have
\ben
\nonumber
p \gamma w_{\min}^2 
\onenorm{\hat{Z} - Z^{*}}    & \le  &
\ip{R, Z^{*} - \hat{Z}} \le
\ip{Y Y^T - \E(YY^T), \hat{Z} - Z^{*} } + \ip{\E B - R,  \hat{Z} - Z^{*} } \\
&&
\label{eq::signal}
- \ip{ (\lambda - \E  \lambda ) (E_n - I_n),   \hat{Z} - Z^{*} } =:
W_1 + W_2 + W_3.
\een
\end{lemma}
We prove Lemma~\ref{lemma::signalgrad} in Section~\ref{sec::proofofsignal}.
To bound the excess risk $\ip{R, Z^{*} - \hat{Z}}$,
we present an upper bound on the bias term $W_2 := \abs{\ip{\E B - R,  \hat{Z} - Z^{*}
  }}$ as well as the large deviation bound on $W_3:= \abs{
  \ip{(\lambda - \E  \lambda ) (E_n - I_n),   \hat{Z} - Z^{*}}}$ over
all $\hat{Z} \in \M_{\opt}$. Lemma~\ref{lemma::optlocal} shows that
each term takes out at most a fraction of the total signal strength on
the LHS of \eqref{eq::signal} respectively.
We then obtain an upper bound on the variance term $W_1$ in
Lemma~\ref{lemma::YYlocal} in the local neighborhood of $Z^{*}$.
\begin{lemma}
  \label{lemma::optlocal}
  Suppose all conditions in Theorem~\ref{thm::exprate} hold.
  Let $\inv{4(n-1)}<\xi \le  w_{\min}^2/16$.
  Then
\ben
\label{eq::event0}
\forall \hat{Z} \in \M_{\opt}, \quad
\abs{\ip{\E B - R,  \hat{Z} - Z^{*} }} =: \abs{W_2}
& \le &
2 p \gamma \onenorm{\hat{Z} - Z^{*} } \big(\xi  + \inv{4(n-1)}\big);
\een
Under the settings of Theorem~\ref{thm::YYaniso},
we have with probability at least $1-\exp(cn)$, 
\bens
\forall \quad \hat{Z} \in \M_{\opt}, \quad
\abs{\ip{ (\lambda - \E  \lambda ) (E_n - I_n),   \hat{Z} - Z^{*} }} =: \abs{W_3}
& \le&
\inv{3} \xi p \gamma  \onenorm{\hat{Z} - Z^{*}}.
\eens
\end{lemma}

\begin{lemma}
  \label{lemma::YYlocal}
Suppose all conditions in Theorem~\ref{thm::YYaniso} hold.
Let $r_1 \le 2 q n$ for a positive integer $1 \le q < n$. 
Then on event $\G_{1}$, where $\prob{\G_{1}} \ge 1-{c}/{n^2}$, we
have for $\xi \le  w_{\min}^2/16$,
\bens
\lefteqn{
\sup_{\hat{Z} \in \M_{\opt} \cap (Z^{*} + r_1 B_1^{n \times n})}
  \abs{\ip{Y Y^T - \E(YY^T), \hat{Z} - Z^{*} }} \le \frac{5}{6} \xi p \gamma \onenorm{\hat{Z} - Z^{*}} }\\
& + & 
C' (C_0 \max_{i} \twonorm{H_i}) \left(n  \sqrt{p \gamma}  + (C_0
  \max_{i} \twonorm{H_i}) \sqrt{n p} \right) \left\lceil \frac{r_1}{2n}\right\rceil \sqrt{\log (2e n/  \lceil \frac{r_1}{2n}\rceil)}.
\eens
\end{lemma}

The final result we need for both the local and global analyses is to
verify a non-trivial global curvature of the excess risk $\ip{R, Z^{*}
  - \hat{Z}}$ for the feasible set $\M_{\opt}$ at the maximizer
$Z^{*}$.
\begin{lemma}{\textnormal{(Excess risk lower bound)~\citep{Zhou23a}}}
  \label{lemma::onenorm}
  Let $R$ be as in Definition~\ref{def::reference} and $Z^{*} =
  u_2 u_2^T$. Then
 \ben
 \label{eq::Rlower}
 \text{for every } \quad Z \in \M_{\opt}, \quad \ip{R, Z^{*} - Z}
 \ge p \gamma w_{\min}^2 \onenorm{Z- Z^{*}}.
 \een
\end{lemma}
This is already stated as the lower bound in~\eqref{eq::signal}.
Hence we can control $\xi \le w_{\min}^2 /16$ and the sample size
lower bound in~\eqref{eq::NKlower} so that each 
component $W_2, W_3$ is only a fraction of the signal at the level of
$p \gamma w_{\min}^2\onenorm{(Z^{*} -  Z)}$.
 We prove Lemma~\ref{lemma::optlocal} in  Section~\ref{sec::proofofoptlocal}.
 We prove Lemma~\ref{lemma::YYlocal} in 
Section~\ref{sec::proofofYYlocal}. 
In particular, the proof of Lemma~\ref{lemma::YYlocal} requires the 
large deviation bound in Theorem~\ref{thm::YYaniso}.

\subsection{Proof of  Theorem~\ref{thm::exprate}}
\label{sec::proofexprate}
\begin{proofof2}
Let $\lceil  r_1/(2n) \rceil =: q< n$ and $r_1 :=\onenorm{\hat{Z}  - Z^{*}}$.
Then, we have by Lemma~\ref{lemma::YYlocal}, on event $\G_1$,
  \ben 
\nonumber 
\lefteqn{\G_{1}:\quad \sup_{\hat{Z} \in \M_{\opt} \cap (Z^{*} + r_1
    B_1^{n \times n})}  \abs{\ip{Y Y^T - \E(YY^T), \hat{Z} - Z^{*} }}
\le \frac{5}{6}\xi p \gamma \onenorm{\hat{Z} - Z^{*}}  +
}\\
\label{eq::eventG1}
&&  
C' C_0 \max_{i} \twonorm{H_i} \left(n \sqrt{p \gamma} + (C_0 \max_{i} \twonorm{H_i}) \sqrt{n p }\right) \left\lceil 
  \frac{r_1}{2n}\right\rceil \sqrt{\log (2e n/ \lceil    \frac{r_1}{2n}\rceil)},
\een
and by Lemma~\ref{lemma::optlocal}, $\forall \hat{Z} \in \M_{\opt}$,
\ben 
\label{eq::eventG2}
\text{ on event } \; \G_{2} : &&
\abs{W_3} = \abs{\ip{ (\lambda - \E  \lambda ) (E_n - I_n),   \hat{Z} - Z^{*} }}
\le \inv{3} \xi p \gamma \onenorm{\hat{Z} - Z^{*}}.
\een
Let $C,  C', C_1, c, c_1, c_2,....$ be absolute positive constants.
For the rest of the proof, we assume  $\G_{1} \cap \G_{2}$ holds,
where $\prob{\G_{2}\cap \G_1} \ge 1- \exp(cn)-c/n^2$.
By Lemmas~\ref{lemma::signalgrad} and~\ref{lemma::onenorm}, 
\eqref{eq::event0}, \eqref{eq::eventG1}, and \eqref{eq::eventG2}, we have  for all $\hat{Z}
\in \M_{\opt} \cap (Z^{*} + r_1 B_1^{n\times n})$,
\ben
\nonumber
0 & \le & p \gamma w_{\min}^2 \onenorm{Z^{*} - \hat{Z}} \le
\ip{R, Z^{*} - \hat{Z}}  \le
\abs{\ip{Y Y^T - \E(YY^T), \hat{Z} - Z^{*} }} + \abs{W_2} + \abs{W_3} 
\\
\label{eq::Rleft}
& \le&  (10/3) \xi  p \gamma \onenorm{\hat{Z} - Z^{*}} + 
C' (C_0  \max_{i} \twonorm{H_i} ) q \sqrt{\log (2e n/ q)} \left(n
  \sqrt{p \gamma} + C_0 (\max_{i} \twonorm{H_i}) \sqrt{n p}\right). 
\een
By moving $(10/3)\xi  p \gamma \onenorm{\hat{Z} - Z^{*}} < \frac{5}{24}  p
\gamma  w_{\min}^2 r_1$, where $\xi \le  w_{\min}^2/16$, to the LHS of \eqref{eq::Rleft}, we have
\ben
\nonumber
\lefteqn{\frac{19}{24} p \gamma w_{\min}^2 r_1
\le  p \gamma w_{\min}^2 \onenorm{Z^{*} - \hat{Z}} -   (10/3) \xi  p \gamma \onenorm{\hat{Z} - Z^{*}} }\\
\label{eq::finalcond}
& \le  & C C_0 (\max_{i} \twonorm{H_i}) \sqrt{\log (2e n/q)} \left(2
  n q \sqrt{p \gamma}  + C_0 (\max_{i} \twonorm{H_i}) q \sqrt{n p }\right).
\een
Then $q$ must satisfy one of the following two conditions
in order to guarantee~\eqref{eq::finalcond}:  for $r_1 \le 2n q$,
\bit
\item[1.]
Under the assumption that 
  $\Delta^2 w_{\min}^4 = p \gamma w_{\min}^4 = \Omega((C_0 \max_{i} \twonorm{H_i})^2 )$,
  suppose
  \bens
  p \gamma w_{\min}^2 
  & \le  &  C_1 C_0 (\max_{i} \twonorm{H_i}) \sqrt{\log (2e n/q)}
  \sqrt{p \gamma},  \quad    \text{ which implies that } \\
\; \frac{c p \gamma w_{\min}^4}{(C_0 \max_{i} \twonorm{H_i})^2}
& \le &  \log (2e n/q)
\text{  and hence } \quad 
q \le  n \exp\big(-\frac{c_1 \Delta^2 w_{\min}^4}{(C_0 \max_{i} \twonorm{H_i})^2}\big). 
\eens
\item[2.]
Alternatively, under the assumption that  $n p \gamma^2
w_{\min}^4=\Omega( C_0^4)$, we require
\bens
p \gamma w_{\min}^2 r_1 \le
 2 qn p \gamma w_{\min}^2
\le  C_2 (C_0 \max_{i} \twonorm{H_i})^2 \sqrt{\log (2e n/q)} q \sqrt{n p },
\eens
and hence $n p \gamma^2 w_{\min}^4 \le  C_3 (C_0 \max_{i}
\twonorm{H_i})^4  \log (2e n/q)$,  which implies that 
\bens
\log (q/2e n) & \le & - \frac{C' n p 
\gamma^2 w_{\min}^4}{(C_0 \max_{i} \twonorm{H_i})^4}, \quad \text{
resulting in }\quad
q \le n \exp\big(- \frac{c_2 n p \gamma^2 w_{\min}^4}{(C_0 \max_{i}
  \twonorm{H_i})^4}\big).
\eens
\eit
Putting things together, we have
\ben
\label{eq::localSNR}
\onenorm{Z^{*} - \hat{Z}} & := & r_1 \le 2 q n \le  n^2 \exp(-C s^2
w_{\min}^4), \; \; \;\text{ where $s^2$ is as defined
  in~\eqref{eq::SNR2},}  \\
\nonumber
\text{ and } \quad
q & \le & n \exp(- C s^2 w_{\min}^4) \; \;\text{ for } \; \;
s^2 = \frac{\Delta^2}{(C_0 \max_{i} \twonorm{H_i})^2} \wedge
\frac{n p \gamma^2}{(C_0 \max_{i} \twonorm{H_i})^4}.
\een
\end{proofof2}

\section{Proof of key lemmas for Theorem~\ref{thm::exprate}}
\label{sec::keylemmas}
First we state two facts. We prove 
Fact~\ref{fact::trP2ZZ} in Section~\ref{sec::expskeleton}.
\begin{fact}
  \label{fact::diag}
By  the supplementary Proposition~\ref{prop::optsol}, we have for the optimal solution $\hat{Z}$ as in SDP~\eqref{eq::sdpmain},
\bens 
\diag(\hat{Z}) = \diag(Z^{*}) = I_n \quad  \text{ and} \quad 
\ip{I_n, \hat{Z} - Z^{*} }= \tr(\hat{Z}) - \tr(Z^{*}) = 0. 
\eens 
\end{fact}

\begin{fact}
  \label{fact::trP2ZZ}
  Recall $\abs{\MC_1} = n_1$ and  $\abs{\MC_2} = n_2$.
  Let $Z^{*}$ be as defined in Lemma~\ref{lemma::ZRnormintro}:
  \ben
  \label{eq::refoptsol}
Z^{*} =   u_2 u_2^T = \left[\begin{array}{cc}
  E_{n_1} &-E_{n_1 \times n_2} \\
 - E_{n_2 \times n_1} & E_{n_2} 
\end{array}
\right].
\een
Then, for  $P_2 = Z^{*}/n$, we have
$\tr(P_2 (Z^{*} - \hat{Z})) = 
\inv{n} \tr(Z^{*} (Z^{*} - \hat{Z})) =
\inv{n} \onenorm{Z^{*} -
    \hat{Z}} \le 2 (n-1)$.
\end{fact}

\subsection{Proof of Lemma~\ref{lemma::signalgrad}}
\label{sec::proofofsignal}
\begin{proofof2}
Denote by $\tilde{A} =Y Y^T -\lambda E_n$.
Now for $B$ as in~\eqref{eq::defineBintro} and $R = \E(Y)\E(Y)^T$, we have
\ben
\label{eq::redefineB}
\E B -R & := & \E (Y Y^T) -\E(Y)\E(Y)^T- \E \lambda (E_n - I_n) - \E \tau I_n.
\een
Clearly $\hat{Z}, Z^{*} \in 
\M_{\opt}$ by definition, we have $\ip{\hat{Z}, I_n}
=\ip{I_n, Z^{*}} = n$.
By optimality of $\hat{Z} \in \M_{\opt}$, we have by~\eqref{eq::defineBintro} and \eqref{eq::sdpmain},
\ben
\nonumber
\ip{B, \hat{Z} - Z^{*}  }
& := &  \ip{Y Y^T -\lambda (E_n - I_n) - \E \tau 
  I_n,  \hat{Z} - Z^{*}  } \\
\label{eq::posadd}
& := &  \ip{Y Y^T -\lambda E_n,  \hat{Z} - Z^{*}  } \ge 0,
\een
where \eqref{eq::posadd} holds since
$\hat{Z} := \arg\max_{Z  \in \M_{\opt}} \ip{\tilde{A}, Z}$ and $\ip{I_n, \hat{Z} -  Z^{*} } = 0$.
Now we have by \eqref{eq::redefineB},
\ben
\nonumber
\lefteqn{M :=  \E(YY^T)- \E(Y)\E(Y)^T -  \lambda (E_n - I_n) - \E \tau
  I_n }\\
\nonumber
& = &   \E(Y Y^T) - \E(Y)\E(Y)^T -  \E \lambda (E_n - I_n) - \E \tau I_n 
-  (\lambda - \E  \lambda ) (E_n - I_n) \\
\label{eq::Wbound}
& = &  \big( \E B - R\big) -  (\lambda - \E  \lambda ) (E_n - I_n).
\een
Then we have by \eqref{eq::posadd} and \eqref{eq::Wbound},
\bens
\lefteqn{
\ip{R, Z^{*} - \hat{Z}} := \ip{-\E(Y)\E(Y)^T, \hat{Z} - Z^{*} }}\\
  & \le &
  \ip{Y Y^T  -  \lambda (E_n - I_n) - \E \tau 
    I_n,   \hat{Z} - Z^{*} } - \ip{\E(Y)\E(Y)^T, \hat{Z} - Z^{*} }\\
  & = &
  \ip{Y Y^T - \E(YY^T), \hat{Z} - Z^{*} } +
  \ip{\E B -R,  \hat{Z} - Z^{*} } -    \ip{(\lambda - \E  \lambda)
    (E_n - I_n), \hat{Z} - Z^{*} }.
\eens
Thus we have the RHS of \eqref{eq::signal} holds.
Now the LHS of  \eqref{eq::signal} holds by
Lemma~\ref{lemma::onenorm}.
\end{proofof2}

\subsection{Proof of Lemma~\ref{lemma::optlocal}}
\label{sec::proofofoptlocal}
 \begin{proofof2}
To be fully transparent, we now decompose the bias term into three components.
Intuitively, $W_0, W_2$ and $\mathbb{W}$ arise due to the imbalance in 
variance profiles. Hence (A2) is needed to control this bias.
\begin{proposition}{\textnormal{{\bf (Bias decomposition)}~\citep{Zhou23a}}}
  \label{prop::biasfinal}
  Denote by 
\ben
\label{eq::defineW2}              
 W_0   & =&
(V_1 -V_2) \left[\begin{array}{cc} w_2 I_{n_1} & 0 \\
0 & -w_1 I_{n_2} \end{array}    \right]\quad \text {and} \quad
 W_2  := 
\frac{V_1-V_2}{n}
\left[
\begin{array}{cc}
E_{n_1} &  0 \\
0  & - E_{n_2} 
\end{array}
\right].
\een 
Then for $W_0, W_2$ as defined in~\eqref{eq::defineW2},
\ben
\label{eq::wow}
\E B - R 
& =&  W_0 - \mathbb{W} -\frac{\tr(R)}{(n-1)} (I_n - \frac{E_n}{n}),\quad
\text{ where} \; \; \\
\label{eq::wow2}
\mathbb{W}
& := & W_2 
+
\frac{ (V_1 - V_2)(w_2 - w_1)}{n} E_n, \quad \text{and when } V_1 = V_2, W_0 = \mathbb{W}=0.
\een
\end{proposition}

\begin{lemma}\textnormal{(Deterministic bounds)}
  \label{lemma::TLbounds}
  Let $Y$ be as specified in Definition~\ref{def::estimators}. 
By definition of $\tau$ and $\lambda$,
\ben
 \label{eq::lambdabound}
(n-1) \abs{\lambda - \E \lambda}
& = &  \abs{\tau - \E \tau} \le \twonorm{YY^T - \E (Y Y^T)}.
\een
\end{lemma}
By Lemma~\ref{lemma::TLbounds} and Theorem~\ref{thm::YYaniso},
we have with probability at least $1 - 2 \exp(- c_6 n)$,
\bens
\abs{(\lambda - \E  \lambda ) \ip{(E_n - I_n),   \hat{Z} - Z^{*} }}
& \le&
\frac{\twonorm{YY^T - \E (YY^T)}}{n-1}
\abs{\sum_{i\not=j } (\hat{Z} -
Z^{*} )_{ij} } \le
\inv{3} \xi p \gamma \onenorm{\hat{Z} - Z^{*}}.
\eens
Now we have by  Proposition~\ref{prop::biasfinal},
\bens
\ip{\E B - R,  \hat{Z} - Z^{*} } 
& = &
\ip{W_0 - \BW, \hat{Z} - Z^{*} }  - \frac{\tr(R)}{n-1} \ip{I_n -
  E_n/n, \hat{Z} - Z^{*} },
\eens
where $W_0$ is a diagonal matrix as in~\eqref{eq::defineW2}.
Hence by Fact~\ref{fact::diag},  (A2), we have for $\BW$ as in~\eqref{eq::wow2},
\ben
\label{eq::BRdecomp}
\ip{W_0, \hat{Z} - Z^{*} }  & = & \ip{W_0, \diag(\hat{Z} - Z^{*}) }  =
0; \\
\label{eq::W2sum}
\abs{\ip{\BW, \hat{Z} - Z^{*} }} & \le &
\norm{\BW}_{\max} \onenorm{ \hat{Z} - Z^{*} }  \le
\frac{2\abs{V_1 - V_2}}{n} \onenorm{(Z^{*} - Z)}   \le
\xi p \gamma \onenorm{(Z^{*} - Z)}.
\een
Hence, by the triangle inequality and ${\tr(R)}/{n}  =  p \gamma w_1 w_2$,
\bens 
\ip{\E B - R,  \hat{Z} - Z^{*} } 
& \le & 
\abs{\ip{\BW, \hat{Z} - Z^{*} } } + \abs{\frac{\tr(R)}{n(n-1)} \ip{E_n, \hat{Z} - Z^{*} }} \\
& \le & \frac{2 \abs{V_1 - V_2}}{n} \onenorm{\hat{Z} - Z^{*} } + \frac{p
 \gamma w_1 w_2}{n-1} \onenorm{\hat{Z} - Z^{*} } \\
& \le &
p \gamma \onenorm{\hat{Z} - Z^{*} } \big(\xi  +  \inv{4(n-1)} \big)
\le 2 \xi p \gamma \onenorm{\hat{Z} - Z^{*} }.
\eens
The lemma thus holds.
\end{proofof2}

\begin{remark}{\textnormal{\bf Balanced Partitions}}
  \label{rm::target}
Suppose that $w_1 = w_2$, then we may hope to obtain a tighter  bound
for $\ip{\E B - R,  \hat{Z} - Z^{*} }$.
Indeed, it is our conjecture that~\eqref{eq::W2sum} can be substantially 
tightened for the balanced case, where $\BW = W_2$.
For $w_1 = w_2$, we have (a) $\inv{n}{\vecone_n^TZ^* u_2} =\vecone_n^T
u_2=0$, and (b) for $\BW = W_2$, 
\ben
\label{eq::W2pro}
\inv{\abs{V_1 - V_2}}
\abs{\ip{\BW, \hat{Z} - Z^{*}} } & = &
\inv{n}  \abs{\sum_{i, j \in \MC_1, j\not= i} (Z^{*}_{ij} -\hat{Z}_{ij}) -
  \sum_{i, j \in \MC_2, j \not= i}  (Z^{*}_{ij} - \hat{Z}_{ij})}\\
& = &
\label{eq::target}
\inv{n}
\abs{\ip{\left[
\begin{array}{cc}
E_{n_1} &  0 \\
0  & - E_{n_2} 
\end{array}
\right], Z^{*} - \hat{Z}  }} =\inv{n} \abs{\vecone_n^T( \hat{Z} -Z^{*}) u_2}   = 
\inv{n} \abs{u_2^T \hat{Z} \vecone_n} \\
& = &
\label{eq::rowsum}
\inv{n}  \abs{\sum_{i \in \MC_1} \sum_{j \in [n], j \not= i} \hat{Z}_{ij} -
    \sum_{i \in \MC_2} \sum_{j \in [n], j \not= i}  \hat{Z}_{ij}},
  \een
  where recall $\diag(\hat{Z}) =\diag(Z^{*}) = I_n$ and
  \eqref{eq::target}
  holds by symmetry and (a).
Suppose that we add an additional
 constraint such that the row sums of $\hat{Z}$ are equal, namely,
 for all $i \in  [n]$, $\sum_{j \in \{[n]\setminus i\}} \hat{Z}_{ij} =
 r$ for some chosen $r$, then~\eqref{eq::rowsum} becomes 0, since the
 average row sum for each cluster will be the same. Without this constraint, a bound tighter
 than $\onenorm{Z^{*} - \hat{Z}} /n$, cf.~\eqref{eq::W2sum} and the supplementary Fact~\ref{fact::Zu2}, is expected for \eqref{eq::target} due to cancellations, since $\forall (i, j) \in \MC_k, k=1, 2$, we have $(Z^{*}_{ij} - \hat{Z}_{ij}) \ge 0$ for the two block-wise sums in~\eqref{eq::W2pro}.
Although we do not have a proof, it is our conjecture that for the
balanced case where $n_1  = n_2$, the bias term essentially becomes a
small order term compared with the variance; cf. the supplementary
Facts~\ref{fact::Rtrace},~\ref{fact::Zu2} and~\eqref{eq::BRdecomp}.
See also the related Fact~\ref{fact::trP2ZZ}, where the bound is tight.
\end{remark}

\subsection{Proof of Lemma~\ref{lemma::YYlocal}}
\label{sec::proofofYYlocal}
\begin{proofof2}
Let $\Lambda = YY^T - \E YY^T$ and  $\Psi := \Z \Z^T - \E \Z \Z^T$. 
It remains to obtain an upper bound for the quantity in~\eqref{eq::signal2}.
The local analysis on~\eqref{eq::signal2} relies on the operator norm
bound we obtained in Theorem~\ref{thm::YYaniso} as already shown in
Lemma~\ref{lemma::optlocal}.
The key distinction from the global analysis in~\cite{Zhou23a}
is, we need to analyze the projection
operators in the following sense.
Following the proof idea of~\cite{FC18},
we define the following projection operator: 
\bens
\cp: M \in \R^{n \times n} &\rightarrow& P_2 M + M P_2 - P_2 M P_2
\quad \text{for } \; P_2 = Z^{*} /{n}, \\
\cp^{\perp}: M \in \R^{n \times n} & \rightarrow &
M-\cp(M) = (I_n-P_2) M (I_n -P_2).
\eens
Now suppose that we use the following decomposition to bound
for all $\hat{Z} \in \M_{\opt} \cap (Z^{*} + r_1 B_1^{n \times n})$,
where $\M_{\opt} := \{Z: Z \succeq 0, \diag(Z) = I_n\} \subset [-1, 1]^{n \times n}$:
\bens
\ip{\Lambda, \hat{Z} -Z^{*}}
& = &
\ip{\Lambda, \cp(\hat{Z} -Z^{*})}  +
\ip{\Lambda, \cp^{\perp}(\hat{Z} -Z^{*})} =: S_1(\hat{Z}) +
S_2(\hat{Z}), \\
\text{where }  \quad
S_1(\hat{Z})  & = & \ip{\Lambda, \cp(\hat{Z} -Z^{*})}   =
\ip{\cp(\Lambda), \hat{Z} -Z^{*}}, \quad \text{and } \\
S_2(\hat{Z}) & := & \ip{\Lambda, \cp^{\perp}(\hat{Z} -Z^{*})} = \ip{\Lambda, \cp^{\perp}(\hat{Z} )} \le  \twonorm{\Lambda} \onenorm{Z^{*} -  \hat{Z}}/n. 
\eens 
First, we control $S_1(\hat{Z})$ uniformly for all $\hat{Z} \in 
\M_{\opt} \cap (Z^{*} + r_1 B_1^{n \times n})$ in
Lemma~\ref{lemma::S1}.
We then obtain an upper bound for $S_2(\hat{Z})$ in
Lemma~\ref{lemma::S2}, uniformly for all $\hat{Z} \in \M_{\opt}$.
Lemmas~\ref{lemma::S1} and~\ref{lemma::S2} show that 
the total {\it noise} contributed by the second term, namely, $S_2(\hat{Z})$ and
a portion of that by $S_1(\hat{Z})$ is only a fraction of the total
signal strength on the LHS of~\eqref{eq::signal}. This allows us to
prove~\eqref{eq::localSNR}, in combination with the deterministic and
high probability bounds in Lemma~\ref{lemma::optlocal}.
\begin{lemma}
   \label{lemma::S1}
Suppose all conditions in Theorem~\ref{thm::YYaniso} hold.
Then, with  probability at least $1-\frac{c'}{n^2}$, for some
absolute constants $c', C$,
\bens
\lefteqn{
\sup_{\hat{Z} \in \M_{\opt} \cap (Z^{*} + r_1 B_1^{n \times n})}
S_1(\hat{Z})
\le  \frac{2}{3} \xi p \gamma \onenorm{\hat{Z} -Z^{*}} + }\\
&&
C C_0 (\max_{i} \twonorm{H^T_i}) \lceil \frac{r_1}{2n}\rceil
\sqrt{\log (2e n/\lceil \frac{r_1}{2n}\rceil)} \left( n \Delta + C_0
  (\max_{i} \twonorm{H^T_i})\sqrt{np} \right).
\eens
\end{lemma}
\begin{lemma}
  \label{lemma::S2}
Suppose all conditions in Theorem~\ref{thm::YYaniso} hold. 
  Let $\xi \le  w_{\min}^2/16$.
Then, with probability at least $1-\exp(c n)$, for all $\hat{Z} \in \M_{\opt}$,
\bens
S_2(\hat{Z})
& \le &
\twonorm{YY^T -\E(YY^T)} \onenorm{Z^{*} -  \hat{Z}}/n \le \frac{1}{6} \xi p \gamma \onenorm{\hat{Z} -Z^{*}}.
 \eens
\end{lemma}
We give an outline for Lemma~\ref{lemma::S1} in Section~\ref{sec::proofofS1main}.
We prove Lemma~\ref{lemma::S2} in the current section.
Lemma~\ref{lemma::YYlocal} follows immediately from
Lemmas~\ref{lemma::S1} and~\ref{lemma::S2}.
\end{proofof2}
\begin{proofof}{Lemma~\ref{lemma::S2}}
  The proof of  Lemma~\ref{lemma::S2}
follows from arguments in~\cite{FC18}.
Denote the noise matrix by $\Lambda:=YY^T -\E(YY^T)$.
  Since $\hat{Z} \succeq 0$, we have by properties of projection,
  \bens
  \cp^{\perp}(\hat{Z}) = (I-P_2) (\hat{Z}) (I-P_2) \succeq 0 \quad
  \text{ and}   \quad \cp^{\perp}(Z^{*}) =  (I_n-P_2) Z^{*} (I_n -P_2) =  0,
  \eens
where $P_2 = Z^{*}/n$. Thus we have for $\tr(\hat{Z}) = \tr(Z^{*}) = n$,
\bens
\norm{ \cp^{\perp}(\hat{Z}) }_{*}
& := &
\norm{ (I-P_2) (\hat{Z}) (I-P_2)}_{*} = 
\tr( (I-P_2) (\hat{Z}) (I-P_2)) \\
&= & \tr(Z^{*} (Z^{*} - \hat{Z}))/n = \onenorm{Z^{*} -  \hat{Z}}/n,
  \eens
  where $\norm{\cdot}_{*}$ denotes the nuclear norm, the last inequality holds
  by Fact~\ref{fact::trP2ZZ}. Thus we have
  \bens
S_2(Z) = 
\ip{\Lambda, \cp^{\perp}(\hat{Z} )}
& \le & \twonorm{\Lambda} \norm{ \cp^{\perp}(\hat{Z}) }_{*}  = 
\twonorm{YY^T-\E YY^T} \onenorm{Z^{*} -  \hat{Z}}/n.
\eens
The lemma thus holds in view of Theorem~\ref{thm::YYaniso}.
\end{proofof}

\subsection{Proof outline for Lemma~\ref{lemma::S1}}
\label{sec::proofofS1main}
\begin{proofof2}
Throughout this section, we consider $\hat{Z} \in \M_{\opt} \cap 
(Z^{*} + r_1 B_1^{n \times   n})$, where $r_1 \le 2 qn$ for some positive 
integer $q \in [n]$.
As predicted, the decomposition and reduction ideas for the global 
analysis~\cite{Zhou23a}, cf. the supplementary
Section~\ref{sec::reduction},
are useful for the local analysis as well. We now briefly introduce them.
First,  we decompose
\ben
\nonumber
\lefteqn{YY^T - \E (Y Y^T) 
  = YY^T - \E (Y) \E(Y)^T  + \E (Y) \E(Y)^T -\E (Y Y^T) } \\
\label{eq::projection}
& = &\E(Y)(Y-\E(Y))^T + (Y-\E(Y))(\E(Y))^T  +
\hat{\Sigma}_Y - \Sigma_Y,
\een
where $\hat\Sigma_Y =  (Y-\E(Y))(Y-\E(Y))^T$ and
\ben
\nonumber
\hat{\Sigma}_Y - \Sigma_Y
& = &  (Y-\E(Y)) (Y-\E(Y))^T + \E (Y) \E(Y)^T -\E (Y Y^T) \\
\label{eq::covproj}
& = &  (I-P_1) (\Z \Z^T -\E (\Z \Z^T))(I-P_1).
\een
Then we have  by the triangle inequality,
\bens
\lefteqn{S_1 := \ip{\cp(\Lambda), \hat{Z} -Z^{*}} 
\le    \abs{\ip{\cp((Y -\E(Y)) \E(Y)^T), \hat{Z} -Z^{*}} }} \\
&&
\abs{\ip{\cp(\E Y (Y - \E(Y))^T), \hat{Z} -Z^{*}}} +
\abs{\ip{\cp(\hat\Sigma_Y- \Sigma_Y), \hat{Z} -Z^{*}} }.
\eens
We need to first present the following results.

\begin{lemma}
  \label{lemma::ensemble}
  Suppose all conditions in Theorem~\ref{thm::YYaniso} hold.
Then, with probability at least $1-\frac{c}{n^2}$,
  \bens
\lefteqn{   \sup_{\hat{Z} \in \M_{\opt} \cap (Z^{*} + r_1 B_1^{n \times     n})}
   \abs{\ip{\cp(\E(Y)(Y -\E(Y))^T), \hat{Z} -Z^{*}}}}\\
   & =&
   \sup_{\hat{Z} \in \M_{\opt} \cap (Z^{*} + r_1 B_1^{n \times     n})}
   \abs{\ip{\cp((Y -\E(Y)) \E(Y)^T), \hat{Z} -Z^{*}}} \\
  & \le  &
 C_5 C_0 (\max_{i} \twonorm{H^T_i \mu}) \sqrt{p \gamma} r_1 \sqrt{\log (2e n/\lceil \frac{r_1}{2n}\rceil)}.
\eens 
\end{lemma}

\begin{lemma}
  \label{lemma::P2cov}
  Under the conditions in Theorem~\ref{thm::YYaniso},
  we have with probability at least $1-\frac{c'}{n^2}$,
\bens
\lefteqn{
\sup_{\hat{Z} \in \M_{\opt} \cap (Z^{*} + r_1 B_1^{n \times  n})}
\abs{\ip{\cp(\hat\Sigma_Y- \Sigma_Y), \hat{Z} -Z^{*}} }}\\
& \le &
\frac{2}{3} \xi  p \gamma  \onenorm{\hat{Z} -Z^{*}} 
+ C_{10} (C_0 \max_{i} \twonorm{H_i})^2 \lceil \frac{r_1}{2n}\rceil 
\left(\sqrt{n p   \log (2e n/\lceil  \frac{r_1}{2n}\rceil)} +
  \sqrt{r_1} \log(2e n/\lceil \frac{r_1}{2n}\rceil ) \right).
\eens
\end{lemma}
Combining Lemmas~\ref{lemma::ensemble} and~\ref{lemma::P2cov},
We have
with probability at least $1-\frac{c}{n^2}$, for $\Delta^2 =p \gamma$,
\bens
\lefteqn{
\sup_{\hat{Z} \in \M_{\opt} \cap (Z^{*} + r_1 B_1^{n \times  n})}
S_1(\hat{Z}) \le 
\sup_{\hat{Z} \in \M_{\opt} \cap (Z^{*} + r_1 B_1^{n \times  n})}
\abs{\ip{\cp(\hat\Sigma_Y- \Sigma_Y), \hat{Z} -Z^{*}} } }\\
&& +
\sup_{\hat{Z} \in \M_{\opt} \cap (Z^{*} + r_1 B_1^{n \times  n})}
2 \abs{\ip{\cp((Y -\E(Y)) \E(Y)^T), \hat{Z} -Z^{*}} }\\
& \le &
\frac{2}{3} \xi  p \gamma  \onenorm{\hat{Z} -Z^{*}} +
2 C_5 C_0 (\max_{i} \twonorm{H^T_i \mu}) 
\lceil \frac{r_1}{2n}\rceil \sqrt{\log (2e n/\lceil
  \frac{r_1}{2n}\rceil)} n \Delta + \\
&&  C_{10}  (C_0 \max_{i} \twonorm{H_i})^2 \lceil \frac{r_1}{2n}\rceil \sqrt{  \log (2e n/\lceil  \frac{r_1}{2n}\rceil)} \left(\sqrt{n p } + \sqrt{r_1} \sqrt{\log(2e n/ \lceil \frac{r_1}{2n}\rceil)} \right) \\
& \le &
\frac{2}{3} \xi  p \gamma  \onenorm{\hat{Z} -Z^{*}} 
+ C' C_0 (\max_{i} \twonorm{H^T_i}) \lceil \frac{r_1}{2n}\rceil \sqrt{\log (2e n/\lceil \frac{r_1}{2n}\rceil)}
\left( n \Delta + C_0 (\max_{i} \twonorm{H^T_i})\sqrt{np} \right),
\eens
where in the last  inequality, we  have by~\eqref{eq::NKlower}, for $r_1 \le 2qn$ and $q = \lceil \frac{r_1}{2n}\rceil$,
\bens
C_0 (\max_{i} \twonorm{H^T_i})\sqrt{r_1} \sqrt{\log(2e n/ \lceil \frac{r_1}{2n}\rceil) } & \le &
C_0 (\max_{i} \twonorm{H^T_i}) \sqrt{2 n q \log (2e n/q)} \\
& \le &
C_0 (\max_{i} \twonorm{H^T_i}) n \sqrt{2 \log (2e)} =O(n \sqrt{p \gamma}),
\eens 
where the first inequality holds since $\max_{q \in [n]}
q \log (2e n/q)  = n \log(2e)$, given that $q \log (2e n/q)$  is
a monotonically increasing function of $q$,  for $1\le q < n$,
and the second inequality holds since $p\gamma = \Omega((C_0 \max_{i}
\twonorm{H^T_i})^2)$.  Lemma~\ref{lemma::S1} thus holds.
We prove Lemmas~\ref{lemma::ensemble} and~\ref{lemma::P2cov} in 
Sections~\ref{sec::proofens} and~\ref{sec::proofP2cov} respectively.
\end{proofof2}

\subsection{Proof of  Lemma~\ref{lemma::ensemble}}
\label{sec::proofens}
\begin{proofof2}
First, we have
\ben
\label{eq::rotation2}
\ip{(Y -\E(Y)) \E(Y)^T P_2,  \hat{Z} -Z^{*}}
&  = &
\ip{P_2 \E(Y) (Y -\E(Y))^T, \hat{Z} -Z^{*}}, \\
\label{eq::rotation}
\ip{P_2 (Y -\E(Y)) \E(Y)^T,  \hat{Z} -Z^{*}}
&  = &
\ip{\E(Y) (Y -\E(Y))^T P_2, \hat{Z} -Z^{*}}.
\een
Then we have  by the triangle inequality,~\eqref{eq::rotation} and~\eqref{eq::rotation2},
  \ben
  \label{eq::rotation3}
  \lefteqn{
\abs{\ip{\cp((Y -\E(Y)) \E(Y)^T), \hat{Z} -Z^{*}} } =
\abs{\ip{\cp(\E(Y)(Y -\E(Y))^T), \hat{Z} -Z^{*}} } }\\
&  \le &
\nonumber
\abs{ \ip{(\E(Y) (Y -\E(Y))^T)P_2, \hat{Z} -Z^{*}}
  - \ip{P_2 \E(Y) (Y -\E(Y))^T P_2, \hat{Z} -Z^{*}}} \\
&&
\nonumber
+ \abs{\ip{P_2 \E(Y) (Y -\E(Y))^T, \hat{Z} -Z^{*}} } =: \abs{T_1}+ \abs{T_2}.
\een
We will bound the two terms $T_1, T_2$ in Lemmas~\ref{lemma::P2minor}
and~\ref{lemma::P2final} respectively.
We prove Lemma~\ref{lemma::P2final} in 
Section~\ref{sec::P2finalproof}. 
We prove Lemma~\ref{lemma::P2minor} in 
the supplementary Section~\ref{sec::P2minorproof}, where we also show that
\eqref{eq::rotation2} and \eqref{eq::rotation} hold.
\begin{lemma}
  \label{lemma::P2minor}
 Under the conditions in Theorem~\ref{thm::YYaniso},
 we have for all $\hat{Z} \in \M_{\opt}$, with probability at least $1-2\exp(-cn)$,
\bens 
 \abs{T_1} := \abs{\ip{(I-P_2) \E(Y) (Y -\E(Y))^T P_2, \hat{Z} -Z^{*}}}
\le C_4 C_0 (\max_{i} \twonorm{H^T_i \mu}) \sqrt{p \gamma} \onenorm{\hat{Z} - Z^*}. 
 \eens 
 Moreover, when $w_1 = w_2$, $T_1 =0$. 
\end{lemma} 

\begin{lemma}
 \label{lemma::P2final}
 Suppose all conditions in Theorem~\ref{thm::YYaniso} hold.
 Suppose $r_1 \le 2 qn$   for some positive integer $1\le q < n$.
 Let $\Delta = \sqrt{p \gamma}$. 
 Then, with probability at least $1-\frac{c}{n^2}$,
 \bens 
\sup_{\hat{Z} \in \M_{\opt} \cap (Z^{*} + r_1 B_1^{n \times n})}
   \ip{P_2 \E(Y) (Y -\E(Y))^T, \hat{Z} -Z^{*}}
\le C_3 C_0 (\max_{i} \twonorm{H^T_i \mu}) \Delta r_1 \sqrt{\log (2e  n/ \left\lceil \frac{r_1}{2n}\right\rceil)}.
  \eens
\end{lemma}
Putting things together, we have 
an uniform upper bound on $\ip{\cp(\E(Y)(Y -\E(Y))^T), \hat{Z} -Z^{*}}$.
\end{proofof2}

\subsection{Proof of Lemma~\ref{lemma::P2cov}}
\label{sec::proofP2cov}
\begin{proofof2}
Upon obtaining the bounds in Lemma~\ref{lemma::ensemble}, 
there are only two unique terms left, 
which we bound in Lemmas~\ref{lemma::finalrate}
and~\ref{lemma::P2PsiP2} respectively.
Notice that by the triangle inequality,  we have
for symmetric matrices $\hat\Sigma_Y- \Sigma_Y = (I-P_1)\Psi(I-P_1)$
and $\hat{Z} -Z^{*}$, 
\ben
\nonumber
\abs{\ip{\cp(\hat\Sigma_Y- \Sigma_Y), \hat{Z} -Z^{*}} }
& \le & 2\abs{\ip{P_2 (\hat\Sigma_Y- \Sigma_Y), \hat{Z} -Z^{*}} }+
\abs{ \ip{P_2 (\hat\Sigma_Y- \Sigma_Y) P_2, \hat{Z} -Z^{*}}  }.
\een
As the first step, we first obtain in Lemma~\ref{lemma::P2PsiP2} a deterministic bound on $\ip{P_2(\hat\Sigma_Y- \Sigma_Y) P_2, \hat{Z} -Z^{*}}$.
Then, in combination with Theorem~\ref{thm::YYaniso},
we obtain the high probability bound as well.
We prove Lemma~\ref{lemma::P2PsiP2} in the supplementary Section~\ref{sec::proofP2ZZP2}.
Next we state Lemma~\ref{lemma::finalrate}, 
which  we prove in Section~\ref{sec::prooffinalrate}. 
\begin{lemma}
  \label{lemma::P2PsiP2}
  Denote by $\Psi := (\Z \Z^T -\E (\Z \Z^T))$.
  Then
\bens
\abs{\ip{P_2(\hat\Sigma_Y- \Sigma_Y) P_2, \hat{Z} -Z^{*}}}
& \le & (1+  \abs{w_1 -   w_2})^2 \twonorm{\Psi} \onenorm{\hat{Z} -Z^{*}} /n.
\eens
Moreover,  under the conditions in Theorem~\ref{thm::YYaniso},
we have with probability at   least $1-\exp(cn)$,
 \bens
\sup_{\hat{Z} \in 
\M_{\opt}} \abs{\ip{P_2(\hat\Sigma_Y- \Sigma_Y) P_2, \hat{Z} -Z^{*}}}
 & \le &
 \frac{1}{3} \xi  p \gamma  \onenorm{\hat{Z} -Z^{*}}. 
 \eens
\end{lemma}

\begin{lemma}
  \label{lemma::finalrate}
 Suppose $r_1 \le 2 qn$   for some positive integer $1\le q < n$.
Suppose all conditions in Theorem~\ref{thm::exprate} hold.
Then, with probability at least $1-\frac{c}{n^2}$,
\bens
\lefteqn{\sup_{\hat{Z} \in \M_{\opt} \cap (Z^{*} + r_1 B_1^{n \times  n})}
  2 \abs{\ip{P_2(\hat\Sigma_Y- \Sigma_Y), \hat{Z} -Z^{*}}}
  \le \frac{1}{3} \xi  p \gamma  \onenorm{\hat{Z} -Z^{*}} + }\\
&&
C_{10}  (C_0 \max_{i} \twonorm{H_i})^2\lceil \frac{r_1}{2n} \rceil \left(\sqrt{np \log(2e n/\lceil \frac{r_1}{2n} \rceil) }+ \sqrt{r_1} \log(2e n/\lceil \frac{r_1}{2n} \rceil) \right).
\eens
\end{lemma}

Lemma~\ref{lemma::P2cov} follows from Lemmas~\ref{lemma::P2PsiP2}
and~\ref{lemma::finalrate}, since with probability at least $1-\frac{c}{n^2} - \exp(cn)$,
\bens
\lefteqn{
\sup_{\hat{Z} \in \M_{\opt} \cap (Z^{*} + r_1 B_1^{n \times  n})}
\abs{\ip{\cp(\hat\Sigma_Y- \Sigma_Y), \hat{Z} -Z^{*}} }}\\
& \le & \sup_{\hat{Z} \in \M_{\opt} \cap (Z^{*} + r_1 B_1^{n \times  n})}
\left(2\abs{\ip{P_2 (\hat\Sigma_Y- \Sigma_Y), \hat{Z} -Z^{*}} }+
  \abs{ \ip{P_2 (\hat\Sigma_Y- \Sigma_Y) P_2, \hat{Z} -Z^{*}} }\right)\\
& \le &
\frac{2}{3} \xi  p \gamma \onenorm{\hat{Z} -Z^{*}}
+ C_{10} (C_0 \max_{i} \twonorm{H_i})^2\lceil \frac{r_1}{2n}\rceil 
\left(\sqrt{n p   \log (2e n/\lceil  \frac{r_1}{2n}\rceil)} +
  \sqrt{r_1} \log(2e n/\lceil \frac{r_1}{2n}\rceil ) \right).
\eens
\end{proofof2}

\subsection{Proof of Lemma~\ref{lemma::P2final}}
\label{sec::P2finalproof}
\begin{proofof2}
We state in Lemmas~\ref{lemma::P2major} and~\ref{lemma::P2majorstar}
two reduction steps.
Denote by $L_i :=\ip{\mu^{(1)} -\mu^{(2)}, \Z_i}$ throughout this
section. 
For all $\hat{Z} \in  M_{\opt}
\cap (Z^* + r_1 B_1^{n \times n})$, we have $\hat{w}  \in B_{\infty}^n
\cap \lceil r_1 /(2n) \rceil B_1^n$, where
\ben
\label{eq::weightnorm}
\hat{w}_j & := & 
\inv{2n}  \onenorm{ (\hat{Z}   -Z^{*})_{\cdot j}}  \le 1 
 \text{ and } \; \sum_{j} \hat{w}_j \le 
 \lceil r_1 /(2n) \rceil,
 \een
 where $B_{\infty}^n$ and $B_1^n$  denote the unit $\ell_{\infty}$ ball and 
 $\ell_1$ ball respectively; cf. Fact~\ref{fact::weights}.
 Denote by $L_1^{*} \ge L_2^{*} \ge \ldots \ge L_n^{*}$
 (resp. $\hat{w}_1^{*} \ge \hat{w}_2^{*} \ge \ldots \ge  \hat{w}_n^{*}$) 
the non-decreasing arrangement of $\abs{L_j}$ (resp. $\hat{w}_j$).
Lemma~\ref{lemma::P2major} is analogous to the supplementary
Lemma~\ref{lemma::tiltproject} for the global analysis.
We then have the reduction as in Lemma~\ref{lemma::P2majorstar}.
We prove Lemmas~\ref{lemma::P2major} and~\ref{lemma::P2majorstar}, in
the supplementary sections~\ref{sec::P2majorproof}
and~\ref{sec::proofP2majorstar} respectively.
\begin{lemma}{\bf (Deterministic bounds)}
\label{lemma::P2major}
Let $\hat{w}_j$ be as in~\eqref{eq::weightnorm}. Then
\ben
\label{eq::2L}
  \abs{\ip{P_2(\E(Y) (Y -\E(Y))^T), \hat{Z} -Z^{*}}}
  \le 4 w_1 w_2 n \sum_{j =1}^{n} \abs{L_i}\hat{w}_j + 
4 w_1 w_2  \abs{\sum_{i =1}^{n} L_i }\onenorm{\hat{Z} -Z^{*}}/(2n).
\een
\end{lemma}

\begin{lemma}{\bf(Reduction to order statistics)}
\label{lemma::P2majorstar}
Under the settings of Lemma~\ref{lemma::P2major}, we have
\ben
\label{eq::L1}
\sup_{\hat{Z} \in \M_{\opt} \cap (Z^{*} +  r_1 B_1^{n \times n})} \abs{  \ip{P_2(\E(Y) (Y -\E(Y))^T), \hat{Z} -Z^{*}}}  \le 3 n \sum_{j =1}^{\lceil {r_1}/{2n}\rceil} L_j^{*}.
\een
\end{lemma}
  It remains to bound the sum on the RHS of \eqref{eq::L1}, namely,
  $\sum_{j =1}^{\lceil {r_1}/{(2n)}\rceil} L_j^{*}$.
  To do so,  we state in Proposition~\ref{prop::orderL} a high probability bound on
  $\sum_{j=1}^q L_j^*$, which holds simultaneously for all $q \in
  [n]$. 
\begin{proposition}
  \label{prop::orderL}
  Suppose all conditions in Lemma~\ref{lemma::twogroup} hold.
  Let random matrix $\Z$ satisfy the conditions as stated  therein.
Let $1 \le q \le n$ denote a positive integer. Then for some absolute constants $C_5, c$,
we have with probability at least $1-\frac{c}{n^2}$, simultaneously for all $q \in [n]$,
\bens 
\sum_{j=1}^q L_j^* \le C_5 C_0 (\max_{i} \twonorm{H^T_i \mu}) \sqrt{p \gamma} q \sqrt{\log (2e n/q)},
\eens
where  $H_i$ and  $\mu = {(\mu^{(1)}-\mu^{(2)})}/\sqrt{p \gamma}$ are as
defined in Theorem~\ref{thm::YYaniso}.
\end{proposition}

We prove Proposition~\ref{prop::orderL} in
Section~\ref{sec::prooforderL}.
It remains to prove Lemma~\ref{lemma::P2final}.
Indeed, we have with probability at least $1-\frac{c}{n^2}$, by
\eqref{eq::2L} and \eqref{eq::L1}, and Proposition~\ref{prop::orderL},
for $ w_1 w_2 \le 1/4$,
\bens 
\lefteqn{
  \sup_{\hat{Z} \in \M_{\opt} \cap (Z^{*} +  r_1 B_1^{n \times n})} \abs{  \ip{P_2(\E(Y) (Y -\E(Y))^T), \hat{Z} -Z^{*}}} }\\
\nonumber
  & \le &
  3 n \sum_{j =1}^{\lceil {r_1}/{(2n)}\rceil} L_j^{*}
  \le C_3 C_0 (\max_{i} \twonorm{H^T_i \mu}) \Delta r_1 \sqrt{\log (2e
    n/\lceil \frac{r_1}{2n}\rceil)},
  \eens
  where it is understood that we set $q = \lceil \frac{r_1}{2n}\rceil 
  \in [n]$ for some $r_1 > 0$.
 The  lemma is thus proved.
\end{proofof2}

\subsection{Proof of Lemma~\ref{lemma::finalrate}}
\label{sec::prooffinalrate}
First we note that by \eqref{eq::covproj},
\ben
\label{eq::P2SigmaY}
\ip{(\hat\Sigma_Y- \Sigma_Y) P_2, \hat{Z} -Z^{*}}
& = &
\ip{P_2(\hat\Sigma_Y- \Sigma_Y),  \hat{Z} -Z^{*}} \\
& = &
\nonumber
\ip{P_2 \Psi, \hat{Z} -Z^{*}} - \ip{P_2 P_1 \Psi , \hat{Z} -Z^{*}} \\
&&    \nonumber
+ \ip{P_2 P_1 \Psi P_1, \hat{Z} -Z^{*}}
  -\ip{P_2 \Psi P_1 , \hat{Z} -Z^{*}}.
\een
We state the following deterministic bounds, which
we need in the proof of Lemmas~\ref{lemma::P2PsiP2}
and~\ref{lemma::P2SigmaY}. Moreover,
we need Lemma~\ref{lemma::P2SigmaY}.
We prove Lemmas~\ref{lemma::P2P1items} and~\ref{lemma::P2SigmaY} 
in the supplementary Sections~\ref{sec::proofofLemmaP2P1}
and~\ref{sec::proofofP2SigmaY} respectively.
\begin{lemma}{\bf(Deterministic bounds)}
\label{lemma::P2P1items}
Let $\Psi = \Z \Z^T -\E (\Z \Z^T)$. Then
\ben
\label{eq::P2major}
\abs{  \ip{P_2 \Psi P_1, \hat{Z} -Z^{*}}}
& \le &   \twonorm{\Psi} \onenorm{\hat{Z} -Z^{*}}/{n}, \\
\label{eq::P2P1P1}
\abs{\ip{P_2 P_1 \Psi P_1, \hat{Z} -Z^{*}}}
& \le &
\abs{w_1 - w_2} 
\twonorm{\Psi } \onenorm{\hat{Z} -Z^{*}} /n, \\
\text{ and } \; 
  \abs{\ip{P_2 P_1 \Psi P_1 P_2, \hat{Z} -Z^{*}}}
& \le &
\label{eq::P2P1P1P2}
\abs{w_1 - w_2}^2 \twonorm{\Psi} \onenorm{\hat{Z} -Z^{*}}/n.
\een
Finally, we will show the following bounds:
   \ben
   \label{eq::P22}
 \ip{P_2 \Psi P_2, \hat{Z} -Z^{*}}
 & \le & \twonorm{\Psi} \onenorm{\hat{Z} -Z^{*}} /n, \\
 \label{eq::P212}
 \abs{\ip{P_2 \Psi P_1 P_2, \hat{Z} -Z^{*}}} &  = &
 \abs{\ip{P_2 P_1 \Psi P_2, \hat{Z} -Z^{*}}}
 \le  \twonorm{\Psi} \abs{w_1 -w_2} \onenorm{\hat{Z} -Z^{*}}/n.
 \een
 When $w_1 = w_2$, the RHS of~\eqref{eq::P2P1P1}, \eqref{eq::P2P1P1P2}  and~\eqref{eq::P212} all become 0, since $P_2 P_1 =0$ in that case.
\end{lemma}

\begin{lemma}{\textnormal{(Deterministic bounds)}}
  \label{lemma::P2SigmaY}
Denote by $\Psi=  \Z \Z^T -\E (\Z \Z^T)$ and 
$\offd(\Psi)$ its off-diagonal component, where the diagonal elements are set to be $0$. Denote by 
\bens 
Q_{j} & := &
\sum_{t  \in \MC_1, t \not=j}  \Psi_{j t} - \sum_{t  \in \MC_2, t \not=j} \Psi_{j t}
= \sum_{t  \in \MC_1}  \offd(\Psi)_{jt} - \sum_{t  \in \MC_2}  \offd(\Psi)_{jt}, \\
S_{j} & := &
\sum_{t  \in [n], t \not=j}  \Psi_{j t} =\sum_{t  \in [n]}  \offd(\Psi)_{j t}.
\eens 
Let $\hat{w}_i$ be the same as in \eqref{eq::weightnorm}.
Then for $\norm{\diag(\Psi)}_{\max} := \max_{j} \abs{\Psi_{jj}}$,
\ben
\label{eq::P2ZZdiag}
\abs{\ip{P_2 \diag(\Psi), \hat{Z} -Z^{*}}}
& \le &
\norm{\diag(\Psi)}_{\max} \onenorm{\hat{Z} -Z^{*}}/n \\
\label{eq::Sdiag}
\abs{\ip{P_2 P_1 \diag(\Psi), \hat{Z} -Z^{*}}}
& \le &
2 \abs{w_1 - w_2} \norm{\diag(\Psi)}_{\max}\onenorm{\hat{Z} -Z^{*}}/n \\
\label{eq::P2ZZ}
\abs{\ip{P_2 \offd(\Psi), \hat{Z} -Z^{*}}}
  &  \le &
  2 \sum_{j \in [n]} \abs{Q_{j}}  \hat{w}_j, \quad \text{ and } \\
\label{eq::doublesum}
\abs{\ip{P_2 P_1 \offd(\Psi), \hat{Z} -Z^{*}}}
& \le &
2 \abs{w_1 - w_2} \sum_{k} \abs{S_{k} }    \cdot \hat{w}_k.
\een
\end{lemma}

Denote by $Q_1^{*} \ge Q_2^{*} \ge \ldots \ge Q_n^{*}$ (resp. $\hat{w}_1^{*} \ge \hat{w}_2^{*} \ge \ldots \ge  \hat{w}_n^{*}$) 
 the non-decreasing arrangement of $\abs{Q_i}$ (resp. $\hat{w}_j$). 
 Denote by $S_1^{*}\ge S_2^{*} \ge \ldots \ge S_n^{*}$
 the non-decreasing arrangement of $\abs{S_i}$  (resp. $\hat{w}_j$).
The RHS of  \eqref{eq::doublesum} and~\eqref{eq::P2ZZ} are bounded 
to be at the same order; cf. Proposition~\ref{prop::orderQ};
Moreover, the RHS of~\eqref{eq::Sdiag} and~\eqref{eq::doublesum} are
both 0 in case we have balanced clusters.
We finish the proof of Lemma~\ref{lemma::finalrate} in  the
supplementary Section~\ref{sec::suppfinalrate}.
Finally, we state Proposition~\ref{prop::orderQ}, which we prove 
in the supplementary Section~\ref{sec::prooforderQ}. 
 \begin{proposition}
   \label{prop::orderQ}
Let $q$ denote a positive integer. Let $C_4, C_5, c$ be absolute  constants.
Under the conditions in Theorem~\ref{thm::YYaniso},
we have with probability at least $1-\frac{c}{n^2}$, simultaneously
for all positive $q \in [n]$, 
\ben
\label{eq::Q2}
\sum_{j =1}^{q} Q_j^{*} 
& \le & C_4 (C_0 \max_{i} \twonorm{H_i})^2 q \big(\sqrt{np \log(2e
  n/q)} + \sqrt{nq} \log(2e n/q) \big), \\
\label{eq::Q3}
\sum_{j =1}^{q} S_j^{*}  & \le &
C_5 (C_0 \max_{i} \twonorm{H_i})^2 q \big(\sqrt{np \log(2e n/q)}
+ \sqrt{nq} \log(2e n/q) \big)
\een
\end{proposition}

\section{Conclusion}
\label{sec::conclude}
Our work is motivated by the two threads of work in combinatorial 
optimization and in community detection, and particularly 
by~\cite{GV15} to revisit the max-cut problem~\eqref{eq::graphcut} and
its convex relaxations.
Centering the data matrix $X$ plays a key role in the statistical
analysis and in understanding the roles of sample size lower bounds
for partial recovery of the clusters~\cite{Zhou23a}.
As in~\cite{Zhou23a},  we focus on the sample size lower bound and
show that a full range of tradeoffs between the sample size and the number of features are 
feasible so long as $s^2$ is lower bounded in the sense 
of~\eqref{eq::kilo} and~\eqref{eq::NKlower}.
In the present work, we further elaborate upon the roles of 
Theorem~\ref{thm::YYaniso} in the local analysis of 
the SDP in Section~\ref{sec::proofofYYlocal}, and the bias and
variance tradeoffs.
More importantly, we prove that the misclassification error decays 
exponentially with respect to the SNR $s^2$ in the present  paper.
We elaborated upon the connections and differences between 
our work and those in~\cite{FC18,GV19}.

\appendix

\section{Preliminary results}

We present preliminary results for the global and local analyses in this section.
These results are either proved in~\cite{Zhou23a}, or follow from
results therein.

\noindent{\bf  Remarks on covariance structures.}
Lemma~\ref{lemma::twogroup} characterizes the two-group design matrix
variance and covariance structures to be considered in Theorems~\ref{thm::SDPmain}
and~\ref{thm::exprate}.
It is understood that when $H_i$ is a symmetric square matrix, it can be
taken as the unique square root of the corresponding covariance
matrix. Expressions in~\eqref{eq::Wpsi}, \eqref{eq::covZ1}, and
\eqref{eq::covZ2} are compatible with each other~\cite{Zhou23a}.
More explicitly, we use the construction in
Definition~\ref{def::WH} to generate $\Z_j \in \R^{p}$ with certain
covariance structures.
When we allow each population to have distinct
covariance structures following Theorem~\ref{thm::exprate},
we have for some  universal constant $C$, and for all $j \in \MC_i$, 
\ben
\label{eq::Zpsi2}
\norm{\Z_j}_{\psi_2}  :=  \sup_{h \in \Sp^{p-1}}\norm{\ip{\Z_j,
    h}}_{\psi_2}  &\le & \norm{W_j}_{\psi_2} \twonorm{H_i} \le C C_0 \max_{i} \twonorm{H_i} \; \text{ since } \\
\forall h \in \Sp^{p-1}, \quad \norm{\ip{\Z_j, h}}_{\psi_2} & = &
 \norm{\ip{H_i W_j,  h}}_{\psi_2} \le  \norm{W_j}_{\psi_2}  \twonorm{H_i^T h},
\een
where $\norm{W_j}_{\psi_2}  \le   C C_0$  by definition
of~\eqref{eq::Wpsi2}.
 Without loss of generality (w.l.o.g.), one may assume that $C_0 = 1$, as one can 
 adjust $H_i$  to control the upper
 bound in~\eqref{eq::Zpsi2} through $\twonorm{H_i}$.

\begin{theorem}{\textnormal{{\bf (Hanson-Wright inequality for
        anisotropic random vectors)}~\cite{Zhou23a}}}
\label{thm::ZHW}
 Let $H_1, \ldots, H_n$ be deterministic $p \times m$ matrices, 
 where we assume that $m \ge p$.
Let $\Z^T_1, \ldots, \Z^T_n \in \R^{p}$ be row vectors of $\Z$.
 We generate $\Z$ according to
 Definition~\ref{def::WH}.
Then we have for $t > 0$, for any $A = (a_{ij}) \in \R^{n \times n}$,
\ben
\nonumber
\lefteqn{\prob{
  \abs{\sum_{i=1}^n  \sum_{j \not=i}^n \ip{\Z_{i}, \Z_{j}} a_{ij}} > t}}\\
\label{eq::genDH}
& \le &
2 \exp \left(- c\min\left(\frac{t^2}{(C_0 \max_{i} \twonorm{H_i})^4  p
      \fnorm{A}^2}, \frac{t}{(C_0 \max_{i} \twonorm{H_i})^2
      \twonorm{A}} \right)\right),
\een
where $\max_{i} \norm{\Z_i}_{\psi_2} \le C C_0 \max_{i} \twonorm{H_i}$
in the sense of \eqref{eq::Zpsi2}.
\end{theorem}

Denote by 
  \ben 
\label{eq::definemu}
\mu:= \frac{\mu^{(1)}-\mu^{(2)}}{\twonorm{\mu^{(1)}-\mu^{(2)}}}
=\frac{\mu^{(1)}-\mu^{(2)}}{\sqrt{p \gamma}} \in \Sp^{p-1}. 
\een 
\begin{lemma}
  \textnormal{{\bf (Projection for anisotropic sub-gaussian random
    vectors)}~\citep{Zhou23a}}
  \label{lemma::anisoproj}
  Suppose all conditions in Theorem~\ref{thm::exprate} hold. Let $R_i = H_i^T$.
  Let $\mu$  be as defined in~\eqref{eq::definemu}. 
 Then
\ben
\label{eq::Zcovdef}      
\norm{\ip{\Z_j, \mu}}^2_{\psi_2}
& \le & C^2_0 \norm{\ip{\Z_j, \mu}}^2_{L_2} :=  C^2_0 \mu^T \cov(\Z_j)  \mu \\
\nonumber
\text{where } \; \mu^T \cov(\Z_j)  \mu
& = &  \mu^T H_i H_i^T \mu = \twonorm{R_i \mu}^2 \; \; \text{ for
  each} \; j \in \MC_i, i =1, 2.
\een 
Thus for any $t > 0$, for some absolute constants $c, c'$,
we have for each $r =(r_1, \ldots, r_n) \in \Sp^{n-1}$ and $u =(u_1,\ldots, u_n) \in \{-1, 1\}^n$ the following tail bounds:
\ben
\label{eq::proanis}
\prob{\abs{\sum_{i=1}^n r_i \ip{\Z_i, \mu}} \ge t} & \le &
2 \exp\left(- \frac{c t^2}{(C_0 \max_{i} \twonorm{R_i \mu})^2}\right),
\; 
\text{ and} \; \\
\prob{\abs{\sum_{i =1}^{n} u_i \ip{\Z_i, \mu }} \ge t}
\label{eq::proanisum}
& \le &
2 \exp\left(-\frac{c' t^2}{n (C_0 \max_{i} \twonorm{R_i \mu})^2}\right).
\een
\end{lemma}

First, we state Lemma~\ref{lemma::tiltproject}.
Upon obtaining \eqref{eq::EYpre}, Lemma~\ref{lemma::tiltproject} is
deterministic and does not depend on the covariance structure of $\Z$.
See Lemma 5.4 in~\cite{Zhou23a}.
\begin{lemma}{\textnormal{{\bf (Reduction: a deterministic comparison
        lemma)}~\citep{Zhou23a}}}
  \label{lemma::tiltproject}
Let $X$ be as in Lemma~\ref{lemma::twogroup}. 
Let $Y$ be as in Definition~\ref{def::estimators}.
Let $\Z_j, j \in [n]$ be row vectors of $X - \E X$ and $\hat{\mu}_n$
be as defined in~\eqref{eq::muhat}. For $x_i \in \{-1,1\}$,
\ben
\sum_{i=1}^n x_i \ip{Y_i -\E Y_i, \mu^{(1)} -\mu^{(2)}}
& \le &
\label{eq::pairwise2}
\frac{2(n-1)}{n} \sum_{i=1}^n \abs{\ip{\Z_i, \mu^{(1)}  -\mu^{(2)}} }. 
\een
\end{lemma}

 \begin{corollary}
   \label{coro::sumY}
Let $X$ be as in Lemma~\ref{lemma::twogroup}. 
Let $Y$ be as in Definition~\ref{def::estimators}.
Under the conditions in Theorem~\ref{thm::YYaniso},
we have with probability at least $1-\exp(-c'n)$,
 \bens
 \lefteqn{
\abs{  \ip{\mu^{(1)} - \mu^{(2)}/\sqrt{p \gamma}, \sum_{j \in C_1} (Y_j - \E(Y_j))  -
    \sum_{j \in C_2} (Y_j - \E(Y_j))}}}\\
& \le  &
2\sum_{i=1}^n \abs{\ip{\frac{\mu^{(1)} -   \mu^{(2)}}{\sqrt{p \gamma}}, \Z_i}} \le C C_0 n \max_{i} \twonorm{R_i \mu}.
\eens
\end{corollary} 
\begin{proof}
Now for $t =  n C C_0 \max_{i} \twonorm{R_i \mu}/2$
\bens
\prob{\max_{ u \in \{-1, 1\}^n } \sum_{j=1}^n u_i \ip{\mu, \Z_j} \ge t}
& \le &  2^n 2 \exp\left(-\frac{c' t^2}{n (C_0 \max_{i} \twonorm{R_i 
      \mu})^2}\right) \\
& \le &  2 \exp(- c n),
\eens
where the last step follows from \eqref{eq::proanisum}.
And hence with probability at least $1-\exp(-c'n)$, we have by Lemma~\ref{lemma::tiltproject},
 \bens
 \lefteqn{
\abs{\ip{\mu^{(1)} - \mu^{(2)}/\sqrt{p \gamma}, \sum_{j \in C_1} (Y_j - \E(Y_j))  -
    \sum_{j \in C_2} (Y_j - \E(Y_j))}}}\\
& \le  &
\max_{x \in \{-1, 1\}^n} \sum_{i=1}^n x_i  \ip{\mu, Y_j - \E(Y_j)} \\
& \le  & 2 \max_{u \in \{-1, 1\}^n} \sum_{j=1}^n  u_j \ip{\mu,  \Z_j}
\le n C C_0 \max_{i} \twonorm{R_i \mu}.
\eens 
\end{proof}

\section{Additional proofs for Theorem~\ref{thm::exprate}}
\label{sec::expskeleton}
First we prove Fact~\ref{fact::trP2ZZ}.

\begin{proofof}{Fact~\ref{fact::trP2ZZ}}
  One can check that $Z^{*}_{ij} \cdot ((Z^{*} -
  \hat{Z})_{ij}) \ge 0$ for all $i, j$.
On the diagonal blocks in \eqref{eq::refoptsol}, we have  $(Z^{*} - \hat{Z})_{ij} \ge 0$, since $Z^{*}_{ij} =1$;  on the off-diagonal blocks, we have $(Z^{*} - \hat{Z})_{ij} \le 0$,
  since $Z^{*}_{ij} =-1$.  Now $\onenorm{\hat{Z} -Z^{*}} \le 2n(n-1)$ by Fact~\ref{fact::diag}.
\end{proofof}

Next we state Facts~\ref{fact::weights} and~\ref{fact::Zu2}.
\begin{fact}
  \label{fact::weights}
  Let $Z_{\cdot j}$ represent the $j^{th}$ column of $Z$.
 Let $\hat{Z} \in \M_{\opt}$ and $\onenorm{\hat{Z} - Z^*} \le r_1 \le 2 qn$,
 where $q = \lceil \frac{r_1}{2n}\rceil < n$ is a positive integer
 and  $Z^{*}$ is as defined in \eqref{eq::refoptsol}.
Denote by
\ben
\label{eq::weights}
   \hat{w}_j =   \inv{2n}\sum_{i \in [n]}  \abs{(\hat{Z} -Z^{*})_{ij}}
   &  := &    \onenorm{(\hat{Z} -Z^{*})_{\cdot j}}/(2n) < 1\\
   \label{eq::weightsum}
   \text{ thus we have} \quad
    \sum_{j =1}^{n} \hat{w}_j &  := & \onenorm{\hat{Z} -Z^{*}}/(2n) \le q  < n.
    \een
 Hence the weight vector $\hat{w} =(\hat{w}_1, 
   \ldots,\hat{w}_n) \in B_{\infty}^{n} \cap q B_1^{n}$, where 
   $B_{\infty}^{n} = \{y \in \R^n: \abs{y_j} \le 1\}$
   and $B_1^{n} = \{y \in \R^n: \sum_{j} \abs{y_j} \le 1\}$
 denote the unit $\ell_{\infty}$ ball and   $\ell_1$ ball respectively.
Moreover,  we have for $\hat{w}_j$ as in~\eqref{eq::weights}, 
\bens
\forall j \in \MC_1, \; \; 
\hat{w}_{j} & = &
-\inv{2n}\big(\sum_{i \in \MC_1}  (\hat{Z} -Z^{*})_{ij} - \sum_{i \in \MC_2}
  (\hat{Z} -Z^{*})_{ij} \big) \quad \text{ and} \\
\forall j \in \MC_2, \; \;
\hat{w}_{j}
& = &
\inv{2n} \big(\sum_{i \in \MC_1}  (\hat{Z} -Z^{*})_{ij} - \sum_{i \in \MC_2}
  (\hat{Z} -Z^{*})_{ij} \big).
\eens
\end{fact}

\begin{proof}
By definition of $Z^{*}$, we have $\forall j \in \MC_1$,
\bens
\forall  i \in \MC_2, &&  (\hat{Z} -Z^{*})_{ij}  \ge 0 \quad \text{ and } \quad
\forall i \in \MC_1, \; \; (\hat{Z} -Z^{*})_{ij}  \le 0,
\eens
and hence
\bens
\hat{w}_{j} & = &
\inv{2n}\big(\sum_{i \in \MC_2} (\hat{Z} -Z^{*})_{ij} - \sum_{i \in
  \MC_1}  (\hat{Z} -Z^{*})_{ij} \big) =\onenorm{(\hat{Z} -Z^{*})_{\cdot j}}/(2n) < 1.
\eens
Similarly, we have $\forall j \in \MC_2$,
\bens
\forall  i \in \MC_1, && (\hat{Z} -Z^{*})_{ij}  \ge 0
\text{ and } \quad \forall  i \in \MC_2,  (\hat{Z} -Z^{*})_{ij}  \le
0,
\eens
and  hence
\bens
\hat{w}_{j}
& = &
\inv{2n} \big(\sum_{i \in \MC_1}  (\hat{Z} -Z^{*})_{ij} - \sum_{i \in \MC_2}
  (\hat{Z} -Z^{*})_{ij} \big)= \onenorm{(\hat{Z} -Z^{*})_{\cdot j}}/(2n) < 1.
\eens
\end{proof}

\begin{fact}
 \label{fact::Zu2}
We have by symmetry of $\hat{Z} -Z^{*}$, under the settings in Fact~\ref{fact::weights},
\bens
  \inv{2n} \abs{\vecone_n^T( \hat{Z} -Z^{*}) u_2}
& = &
\abs{\inv{2n} \sum_{i=1}^n 
  \big(\sum_{j\in\MC_1} (\hat{Z} -Z^{*})_{ij} -\sum_{j \in \MC_2} (\hat{Z}
  -Z^{*})_{ij}\big)} \\
& \le &
\sum_{i=1}^n \abs{\hat{w}_i}=\inv{2n}  \onenorm{ (\hat{Z} -Z^{*})}
= \inv{2n} \sum_{j} \onenorm{(\hat{Z} -Z^{*})_{j\cdot} }.
\eens
\end{fact}

\section{Technical proofs for Section~\ref{sec::keylemmas}}
\label{sec::expdetail}

\subsection{Proof outline of Lemma~\ref{lemma::P2minor}}
\label{sec::P2minorproof}
\begin{proofof2}
Recall by~\eqref{eq::EYpre},
  \bens
\E (Y_i) & = & \left\{ \begin{array}{rl} w_2 (\mu^{(1)} -  \mu^{(2)})    &
\text{ if }  \; i \in \MC_1, \\    w_1 (\mu^{(2)} - \mu^{(1)})&\text{ if }
\; i \in \MC_2. \\
 \end{array}\right.
\eens                                             
 First, we show that \eqref{eq::rotation2} holds.
To see this, we have
     \bens
   \lefteqn{
     \ip{(Y -\E(Y)) \E(Y)^T P_2,  \hat{Z} -Z^{*}}
     = \tr((\hat{Z} -Z^{*}) (Y -\E(Y)) \E(Y)^T P_2)}\\
   &  = &
   \ip{P_2 \E(Y) (Y -\E(Y))^T, \hat{Z}-Z^{*}}.
   \eens
Similarly, \eqref{eq::rotation} holds since
\ben
\nonumber
     \ip{P_2 (Y -\E(Y)) \E(Y)^T,  \hat{Z} -Z^{*}}
 &  = &   \tr((\hat{Z} -Z^{*}) P_2 (Y -\E(Y)) \E(Y)^T) \\
   \label{eq::roationlocal}
   &  = &
   \ip{\E(Y) (Y -\E(Y))^T P_2, \hat{Z}-Z^{*}}.
   \een
 Denote by $V = (\E(Y) (Y -\E(Y))^T) u_2$, which is a vector with coordinates 
$V_1, \ldots, V_n$ such that 
\bens
\forall i \in \MC_1, \quad 
V_i &  = & - w_2
\ip{\mu^{(2)} - \mu^{(1)}, \sum_{j \in C_1} (Y_j - \E(Y_j))  -
  \sum_{j \in C_2} (Y_j - \E(Y_j))},  \\
\forall i \in \MC_2, \quad 
V_i &  = &  w_1 \ip{\mu^{(2)} - \mu^{(1)}, \sum_{j \in C_1} (Y_j - \E(Y_j)) 
  -  \sum_{j \in C_2} (Y_j - \E(Y_j))},
\eens
where $Y_j \in \R^{p}$ represents the $j^{th}$ row vector of matrix $Y$.
 Similarly, we have by \eqref{eq::roationlocal},
   \ben
\nonumber
    \lefteqn{\ip{P_2 (Y -\E(Y)) \E(Y)^T,  \hat{Z} -Z^{*}}    = 
      \tr((\hat{Z}-Z^{*}) \E(Y) (Y -\E(Y))^T P_2) } \\
\nonumber
 &  = & \tr(P_2 (\hat{Z}-Z^{*}) \E(Y) (Y -\E(Y))^T) 
  =  \tr(\frac{u_2 u_2^T}{n} (\hat{Z} -Z^{*}) (\E(Y) (Y -\E(Y))^T)) \\
\nonumber
&  = &
\tr(\frac{u_2^T}{n} (\hat{Z} -Z^{*}) (\E(Y) (Y -\E(Y))^T) u_2) 
= \ip{\inv{n}\big(\sum_{i \in \MC_1}  (\hat{Z} -Z^{*})_{i\cdot} - \sum_{i \in \MC_2}
  (\hat{Z} -Z^{*})_{i \cdot} \big), V}  \\
\nonumber
&  = &
2 \sum_{k}  V_k \cdot
\inv{2n} \big(\sum_{i \in \MC_1}
(\hat{Z} -Z^{*})_{i k} - \sum_{i \in \MC_2}
(\hat{Z} -Z^{*})_{i k} \big) =: S.
\een
Now 
\bens
\lefteqn{\forall k \in \MC_1, \quad
  \inv{2n}
  \big(\sum_{i \in \MC_1}  (\hat{Z} -Z^{*})_{i k} - \sum_{i \in \MC_2}
(\hat{Z} -Z^{*})_{i k} \big)  V_k } \\
&  = &
w_2  \hat{w}_k\ip{\mu^{(2)} - \mu^{(1)}, \sum_{j \in C_1} (Y_j - \E(Y_j)) 
  -  \sum_{j \in C_2} (Y_j - \E(Y_j))}, \\
\lefteqn{
\text{ and }  \; \forall k \in \MC_2, \quad
\inv{2n}\big(\sum_{i \in \MC_1}  (\hat{Z} -Z^{*})_{i k} - \sum_{i \in \MC_2}
(\hat{Z} -Z^{*})_{i k} \big)  V_k }\\
&  = &  w_1 \hat{w}_k \ip{\mu^{(2)} - \mu^{(1)}, \sum_{j \in C_1} (Y_j - \E(Y_j)) 
  -  \sum_{j \in C_2} (Y_j - \E(Y_j))}.
\eens
Hence by~\eqref{eq::weights} and Fact~\ref{fact::weights},
\ben 
\label{eq::finalYY}
S & = &
2 C_S
\cdot \ip{\mu^{(2)} - \mu^{(1)}, \sum_{j \in C_1} (Y_j - \E(Y_j))  -
  \sum_{j \in C_2} (Y_j - \E(Y_j))},  \\
\nonumber
\text{ where } && 
C_S := w_2 \sum_{k \in \MC_1} \hat{w}_k + w_1 \sum_{k \in \MC_2}
\hat{w}_k \le  \sum_{k \in [n]} \hat{w}_k = \onenorm{\hat{Z} -Z^{*}}/(2n).
\een
Next we define
\ben
\label{eq::P2sum}
W & := & 2 u_2^T \E(Y) (Y-\E(Y))^T u_2 /n = 2 u_2^T V/n =
\frac{2}{n}(\sum_{i \in \MC_1} V_i -\sum_{i \in \MC_2} V_i ) \\
& = &
\nonumber
\ip{\mu^{(1)}-\mu^{(2)}, 
  \sum_{j \in \MC_1} (Y_j - \E Y_j) -   \sum_{j \in \MC_2} (Y_j - \E
  Y_j)} 4 w_1 w_2.
\een
Then
\ben
\nonumber
\lefteqn{\ip{P_2 \E(Y) (Y -\E(Y))^T P_2, \hat{Z} -Z^{*}} =
  \tr(P_2 (\hat{Z} -Z^{*})  P_2 \E(Y) (Y -\E(Y))^T) }\\
&  = &
\label{eq::P2YEYP2}
\inv{2n} u_2^T (\hat{Z} -Z^{*}) u_2 \big(2u^T_2 \E(Y) (Y -\E(Y))^T u_2 /n\big)
=: -W \sum_{k} \hat{w}_k,
\een
where \eqref{eq::P2YEYP2} holds
by \eqref{eq::P2sum} and Fact~\ref{fact::trP2ZZ},
since for $P_2 = {u_2 u_2^T}/{n} = Z^*/n$,
\bens
\frac{u_2^T}{2n} (\hat{Z} -Z^{*}) u_2 & = & 
\tr(P_2 (\hat{Z} -Z^{*}))/2 = -\onenorm{\hat{Z} -
  Z^*}/(2n) = -\sum_{k} \hat{w}_k.
\eens
Putting things together, we have by \eqref{eq::weightsum}, \eqref{eq::finalYY} and \eqref{eq::P2YEYP2},
 \bens
 \lefteqn{T_1 :=
\ip{\E(Y) (Y -\E(Y))^T P_2, \hat{Z} -Z^{*}}
- \ip{P_2 \E(Y) (Y   -\E(Y))^T P_2, \hat{Z} -Z^{*}}} \\
& = &
S +W \sum_{k} \hat{w}_k = 
 2\ip{\mu^{(2)} - \mu^{(1)}, \sum_{j \in C_1} (Y_j - \E(Y_j))  -
   \sum_{j \in C_2} (Y_j - \E(Y_j))} \cdot \\
 &&
\big(w_2\sum_{k \in   \MC_1} \hat{w}_k 
 + w_1 \sum_{k \in \MC_2} \hat{w}_k - 2 w_1 w_2 \sum_{k \in [n]}
 \hat{w}_k \big),
 \eens
where 
\bens
\abs{w_2\sum_{k \in   \MC_1} \hat{w}_k 
 + w_1 \sum_{k \in \MC_2} \hat{w}_k - 2 w_1 w_2 \sum_{k \in [n]}
 \hat{w}_k}
& \le &
\abs{w_2\sum_{k \in   \MC_1} \hat{w}_k 
  + w_1 \sum_{k \in \MC_2} \hat{w}_k } \vee \abs{2 w_1 w_2
  \sum_{k \in [n]}
 \hat{w}_k} \\
& \le & \onenorm{\hat{Z} -Z^{*}}/(2n) \le q < n.
\eens
Now, for $w_1= w_2 = 1/2$, the above sum is 0 and hence $T_1 =0$.
Let $r_1 \le 2 qn$.
Then, by Corollary~\ref{coro::sumY}, we have with probability at least $1-2\exp(-cn)$,
for all $\hat{Z} \in \M_{\opt} \cap (Z^{*} + r_1 B_1^{n \times   n})$,
\bens
\lefteqn{
 \abs{T_1}= \abs{ \ip{(\E(Y) (Y -\E(Y))^T)P_2, \hat{Z} -Z^{*}}
    - \ip{P_2 \E(Y) (Y -\E(Y))^T P_2, \hat{Z} -Z^{*}}}} \\
 & \le &
\abs{\ip{\mu^{(2)} - \mu^{(1)}, \sum_{j \in C_1} (Y_j - \E(Y_j))  -
     \sum_{j \in C_2} (Y_j - \E(Y_j))}} \onenorm{\hat{Z} -Z^{*}}/n \\
 &  \le & C_4 C_0 (\max_{i} \twonorm{R_i \mu}) \sqrt{p \gamma} \onenorm{\hat{Z} - Z^*}
 \le  C_4 C_0 (\max_{i} \twonorm{R_i \mu}) r_1 \sqrt{p \gamma}.
 \eens
\end{proofof2}

\subsection{Proof of Lemma~\ref{lemma::P2PsiP2}}
\label{sec::proofP2ZZP2}

\begin{theorem}{\textnormal{~\cite{Zhou23a}}}
\label{thm::YYcovcorr}
  In the settings of Theorem~\ref{thm::YYaniso}, we have
  with probability at least $1 - 2\exp(-c_6 n)$, 
\ben
\label{eq::ZZHop}
\twonorm{\hat\Sigma_Y - \E \hat\Sigma_Y}
& \le &
\twonorm{\Z \Z^T - \E  \Z \Z^T} \le  C_2   C_0^2 \max_{j} \twonorm{\cov(Z_j)} (\sqrt{n p} \vee n).
\een
\end{theorem}

\begin{proofof}{Lemma~\ref{lemma::P2PsiP2}}
First,
\ben
\label{eq::P2PsiP2}
\ip{P_2(\hat\Sigma_Y- \Sigma_Y) P_2, \hat{Z} -Z^{*}}
& = &
\ip{P_2 (I-P_1) \Psi (I-P_1) P_2, \hat{Z} -Z^{*}} \\
& = &
\nonumber
\ip{P_2 \Psi P_2, \hat{Z} -Z^{*}} +
\ip{P_2 P_1 \Psi P_1 P_2, \hat{Z} -Z^{*}} \\
& - &
\nonumber
\ip{P_2 P_1 \Psi P_2, \hat{Z} -Z^{*}} - \ip{P_2 \Psi P_1 P_2, \hat{Z}
  -Z^{*}}.
\een
Hence the  deterministic statement in Lemma~\ref{lemma::P2PsiP2} holds in view of \eqref{eq::P2PsiP2}, \eqref{eq::P2P1P1P2},
~\eqref{eq::P22}, and~\eqref{eq::P212},
\bens
\abs{\ip{P_2(\hat\Sigma_Y- \Sigma_Y) P_2, \hat{Z} -Z^{*}}}
& \le &
\abs{\ip{P_2  \Psi P_2, \hat{Z} -Z^{*}}  }+
\abs{\ip{P_2 P_1 \Psi P_1 P_2, \hat{Z} -Z^{*}} } +  \\
&&
2 \abs{\ip{P_2 P_1 \Psi P_2, \hat{Z} -Z^{*}}} \\
& \le &
\inv{n}\onenorm{\hat{Z} -Z^{*}  } \twonorm{\Psi}
\left(1 + (w_1 -w_2)^2 + 2 \abs{w_2 - w_1} \right)\\
& \le &
\inv{n}\onenorm{\hat{Z} -Z^{*}  } \twonorm{\Psi} (1 + \abs{w_1 -w_2})^2.
\eens
Finally, we have by Theorem~\ref{thm::YYcovcorr}, with probability at least $1-\exp(c' n)$,
\bens
\twonorm{\Psi} = \twonorm{\Z \Z^T -\E \Z \Z^T} &\le&
C_2 C_0^2 \max_{j} \twonorm{\cov(\Z_j)} (\sqrt{n p} \vee n)  \le
\inv{12} \xi n p \gamma,
\eens
where the last step holds so long as~\eqref{eq::NKlower} holds
for sufficiently large constants $C, C_1$ and $\xi \le \frac{w_{\min}^2}{16}$.
Thus, we have with probability at least $1-\exp(c' n)$, for all $\hat{Z} \in \M_{\opt}$,
\bens
\abs{\ip{P_2(\hat\Sigma_Y- \Sigma_Y) P_2, \hat{Z} -Z^{*}}}
& \le &
\label{eq::lastfac}
8 \twonorm{\Psi} \onenorm{\hat{Z} -Z^{*}} /(2n) \le
\frac{1}{3} \xi  p \gamma   \onenorm{\hat{Z} -Z^{*}}.
\eens
The lemma is thus proved.
\end{proofof}

\subsection{Proof of Lemma~\ref{lemma::finalrate}}
\label{sec::suppfinalrate}

\begin{corollary}
\label{coro::orderQ}
Under the settings in Proposition~\ref{prop::orderQ},
we have with probability at least $1 - \frac{c}{n^2}$,
 \ben
 \label{eq::duet}
   \lefteqn{\sup_{\hat{Z} \in \M_{\opt} \cap (Z^{*} + r_1 B_1^{n \times n})}
   \abs{\ip{P_2 \offd(\Psi), \hat{Z} -Z^{*}}}
   \le  4 \sum_{j =1}^{\lceil r_1 /(2n) \rceil} Q_j^{*} }\\
 &\le &
 \nonumber
4 C_4 (C_0 \max_{i} \twonorm{H_i})^2
\lceil \frac{r_1}{2n} \rceil \left(\sqrt{np \log(2e n/\lceil
    \frac{r_1}{2n} \rceil) }+ \sqrt{r_1} \log(2e n/\lceil
  \frac{r_1}{2n} \rceil) \right),
\een
and
\ben
\label{eq::trio}
\lefteqn{
  \sup_{\hat{Z} \in \M_{\opt} \cap (Z^{*} + r_1 B_1^{n \times n})}
\abs{\ip{P_2 P_1 \offd(\Psi), \hat{Z} -Z^{*}}}  \le
4 \abs{w_1 - w_2}
\sum_{j =1}^{\lceil {r_1}/{(2n)} \rceil} S_j^{*} }\\
& \le &
\nonumber
4 \abs{w_1 - w_2}  C_5   (C_0 \max_{i} \twonorm{H_i})^2 \lceil \frac{r_1}{2n} \rceil \left(\sqrt{np \log(2e n/\lceil \frac{r_1}{2n} \rceil) }+ \sqrt{r_1} \log(2e n/\lceil \frac{r_1}{2n} \rceil) \right).
\een 
\end{corollary}

\begin{proofof}{Lemma~\ref{lemma::finalrate}}
  Let $\hat{Z} \in \M_{\opt} \cap ( Z^{*} + r_1 B_1^{n \times n})$, 
  where $r_1 \le 2 q n$  for some positive integer $1\le q \le n$. 
  With probability at least $1-\exp(c_4 n)$, by
  Lemma~\ref{lemma::P2P1items} (cf.  \eqref{eq::P2major} and
  \eqref{eq::P2P1P1}), we have for all $\hat{Z} \in \M_{\opt}$,
\bens
\lefteqn{W_1 := \abs{  \ip{P_2 \Psi P_1, \hat{Z} -Z^{*}}} + 
\abs{\ip{P_2 P_1 \Psi P_1, \hat{Z} -Z^{*}}} } \\
& \le &
(1+\abs{w_1 - w_2} ) \twonorm{\Psi} \inv{n} \onenorm{\hat{Z} -Z^{*}}.
\eens
Thus we have
by the triangle inequality, by~\eqref{eq::P2SigmaY},
   \bens
\lefteqn{\abs{\ip{P_2(\hat\Sigma_Y- \Sigma_Y), \hat{Z} -Z^{*}}}
  \le      \abs{ \ip{P_2 \Psi P_1 ,  \hat{Z} -Z^{*}}} + 
  \abs{\ip{P_2 P_1 \Psi P_1, \hat{Z}   -Z^{*}}}  +}\\
& &   \abs{\ip{P_2 \Psi, \hat{Z} -Z^{*}} } +  \abs{\ip{P_2 P_1 \Psi ,
    \hat{Z} -Z^{*}} } =: W_1 + W_2,
  \eens
  where we decompose
  \bens
\lefteqn{W_2 := \abs{\ip{P_2 \Psi, \hat{Z} -Z^{*}}}  + \abs{\ip{P_2 P_1 \Psi , \hat{Z} -Z^{*}} }   \le} \\
&& \abs{\ip{P_2 \diag(\Psi), \hat{Z} -Z^{*}} } + \abs{\ip{P_2 P_1
    \diag(\Psi) , \hat{Z} -Z^{*}} }  \quad (W_2(\diag)) \\
&& + \abs{\ip{P_2 \offd(\Psi), \hat{Z} -Z^{*}} } +
\abs{\ip{P_2 P_1 \offd(\Psi) , \hat{Z} -Z^{*}} } \quad (W_2(\offd)).
\eens
Then, we have for  all $\hat{Z} \in Z^{*} + r_1 B_1^{n \times n}$,
by \eqref{eq::P2P1P1}, \eqref{eq::P2major}, \eqref{eq::P2ZZdiag},
and \eqref{eq::Sdiag},
\ben
 \nonumber
W_2(\offd) + W_1 &= & \abs{\ip{P_2 \diag(\Psi), \hat{Z} -Z^{*}}} + 
\abs{\ip{P_2 P_1 \diag(\Psi), \hat{Z} -Z^{*}}} +\\
 \nonumber
&& \abs{  \ip{P_2 \Psi P_1, \hat{Z} -Z^{*}}} +
\abs{\ip{P_2 P_1\Psi P_1, \hat{Z} -Z^{*}}} \\
 \nonumber
& \le & (1+\abs{w_1 - w_2}) \max_{k} (\abs{\Psi_{k k}} +\twonorm{\Psi})
\onenorm{\hat{Z} -Z^{*}}/n \\
\label{eq::Psibounds}
& \le & 4 \twonorm{\Psi} \onenorm{\hat{Z} -Z^{*}}/n,
\een
and moreover, we have by~\eqref{eq::duet}, and~\eqref{eq::trio},
\ben
\nonumber
\lefteqn{   \sup_{\hat{Z} \in \M_{\opt} \cap (Z^{*} + r_1 B_1^{n \times n})}  
  W_2(\offd) } \\
& := &
\nonumber
\sup_{\hat{Z} \in \M_{\opt} \cap (Z^{*} + r_1 B_1^{n \times n})}
\big(  \abs{\ip{P_2 \offd(\Psi), \hat{Z} -Z^{*}} } +
\abs{\ip{P_2 P_1 \offd(\Psi), \hat{Z}   -Z^{*}}} \big)\\
& \le &
\label{eq::W2offd}
4 \sup_{\hat{w} \in B_{\infty}^{n} \cap \lceil r_1 /(2n) \rceil B_1^{n}}
\big(\sum_{j =1}^{\lceil r_1 /(2n) \rceil} Q_j^{*} +\abs{w_1 - w_2}\sum_{j =1}^{\lceil {r_1}/{(2n)} \rceil} S_j^{*}\big).
\een
Thus we have for $\Psi = \Z \Z^T -\E(\Z \Z^T)$,
by~\eqref{eq::Psibounds} and \eqref{eq::W2offd}, with probability at
least $1- c_2 /n^2$, and $q =\lceil r_1 /(2n) \rceil$,
\bens
\lefteqn{
\sup_{\hat{Z} \in \M_{\opt} \cap (Z^{*} + r_1 B_1^{n \times n})}
\abs{\ip{P_2(\hat\Sigma_Y- \Sigma_Y), \hat{Z} -Z^{*}}}
\le   \sup_{\hat{Z} \in \M_{\opt} \cap (Z^{*} + r_1 B_1^{n \times n})} 
    4 \twonorm{\Psi} \onenorm{\hat{Z} -Z^{*}}/n
  }\\
  & + &
4 \sup_{\hat{w} \in B_{\infty}^{n} \cap \lceil r_1 /(2n) \rceil B_1^{n}}
\big(\sum_{j =1}^{\lceil r_1 /(2n) \rceil} Q_j^{*} +\abs{w_1 - w_2}\sum_{j =1}^{\lceil {r_1}/{(2n)} \rceil} S_j^{*}\big) \\
& \le &
\inv{3} \xi p \gamma \onenorm{\hat{Z} -Z^{*}} + C_{10}   (C_0 \max_{i}
\twonorm{H_i})^2 \lceil \frac{r_1}{2n} \rceil \left(\sqrt{np \log(2e
    n/\lceil \frac{r_1}{2n} \rceil) }+ \sqrt{r_1} \log(2e n/\lceil
  \frac{r_1}{2n} \rceil) \right).
\eens 
The lemma thus holds. 
\end{proofof}

\begin{proofof}{Corollary~\ref{coro::orderQ}}
For any positive integer $q \in [n]$ and $\hat{w}$ as in~\eqref{eq::weightnorm},
\ben
  \nonumber
  \sup_{\hat{w} \in B_{\infty}^{n} \cap q B_1^{n}}
  \sum_{j \in [n]} \abs{Q_j} \hat{w}_j &\le &
  \sup_{\hat{w} \in B_{\infty}^{n} \cap q B_1^{n}}
  \sum_{j \in [n]} Q^*_j \hat{w}^*_j 
 =   \sup_{\hat{w} \in B_{\infty}^{n} \cap q B_1^{n}}
\big(\sum_{j =1}^{q} Q_j^{*}  \hat{w}^{*}_j +
\sum_{j =q+1}^{n} Q_j^{*}  \hat{w}^{*}_j \big) \\
 \label{eq::Qsum}
&\le &
\sum_{j =1}^{q} Q_j^{*}  + q Q_{q+1}^{*}  \le 2 \sum_{j =1}^{q}
Q_j^{*}.
\een
 Following the same arguments in Lemma~\ref{lemma::P2final}, we have
  by~\eqref{eq::P2ZZ} and \eqref{eq::weightnorm}, 
\ben
\nonumber
\sup_{\hat{Z} \in  M_{\opt}
   \cap Z^* + r_1 B_1^{n \times n}} \abs{\ip{P_2 \offd(\Psi), \hat{Z} - Z^{*}}}
  &  \le &   2  \sup_{\hat{w} \in B_{\infty}^{n} \cap \lceil r_1 /(2n) \rceil B_1^{n}}
  \sum_{j \in [n]} \abs{Q_{j}} \hat{w}_j
  \een
Now  \eqref{eq::duet} follows from~\eqref{eq::P2ZZ}, 
~\eqref{eq::Qsum}, and~\eqref{eq::Q2}, where we set  $q =\lceil r_1 
/(2n) \rceil$.
Notice that for any positive integer $q  \in [n]$,  following the same
argument, we have
 \ben
 \label{eq::Ssum}
 \sup_{\hat{w} \in B_{\infty}^{n} \cap q B_1^{n}}
\sum_{j \in [n]} \abs{S_j} \hat{w}_j 
& \le & \sup_{\hat{w} \in B_{\infty}^{n} \cap q B_1^{n}}
\sum_{j \in [n]} S^*_j \hat{w}^*_j  \le 2 \sum_{j =1}^{q} S_j^{*}.
\een
Hence \eqref{eq::trio}  follows from~\eqref{eq::doublesum}, 
   \eqref{eq::Ssum},~\eqref{eq::Q3},  and   Proposition~\ref{prop::orderQ}. 
\end{proofof}

\subsection{Proof of Lemma~\ref{lemma::P2major}}
\label{sec::P2majorproof}
\begin{proofof2}
  The proof of  Lemma~\ref{lemma::P2major} is entirely deterministic.
    Recall by~\eqref{eq::EYpre},
    we have
\bens
\inv{n} \E(Y)^T u_2 = \inv{n}(\sum_{i\in \MC_1} \E(Y_i) -\sum_{i\in \MC_2} \E(Y_i) )
= 2 w_1 w_2 (\mu^{(1)} -  \mu^{(2)}).
\eens
  Throughout this proof, we denote by
  $\mu := {(\mu^{(1)}-\mu^{(2)})}/\sqrt{p \gamma}$.
  Let $Z^{*}$ be as defined in \eqref{eq::refoptsol}.
  Recall $P_2 = u_2 u_2^T/n = Z^{*}/n$.
Thus we have by symmetry of $\hat{Z} -Z^{*}$ and definition of $\hat{w}_j$ in \eqref{eq::weights},
  \bens
  \lefteqn{  \ip{P_2(\E(Y) (Y -\E(Y))^T), \hat{Z} -Z^{*}} = 
  \inv{n} u_2^T (\E(Y) (Y -\E(Y))^T)( \hat{Z} -Z^{*}) u_2 }\\
  &  = &
  2 w_1 w_2 \sqrt{p \gamma} \sum_{j \in [n]}
  \ip{\mu, Y_j - \E (Y_j)} \cdot
  \big(  \sum_{i \in \MC_1}  (\hat{Z} -Z^{*})_{ji} - \sum_{i \in \MC_2}
  (\hat{Z} -Z^{*})_{ji} \big) \\
   &  = &
 4 n w_1 w_2 \sqrt{p \gamma}
\big(\sum_{j \in \MC_2}  \hat{w}_j \ip{\mu, Y_j - \E (Y_j)} 
  -   \sum_{j \in \MC_1} \hat{w}_j \ip{\mu, Y_j - \E (Y_j)} \big),
\eens
where $\hat{w}_j \ge0$ is as defined in Fact~\ref{fact::weights},
\bens
\forall j \in \MC_2, \quad 
\hat{w}_j 
& := &\inv{2n} \big(\sum_{i \in \MC_1}  (\hat{Z} -Z^{*})_{i j} - \sum_{i \in \MC_2}  (\hat{Z} -Z^{*})_{i j} \big), \\
\forall j \in \MC_1, \quad 
\hat{w}_j  & := & -\inv{2n}
\big(\sum_{i \in \MC_1}  (\hat{Z} -Z^{*})_{i j} - \sum_{i \in \MC_2}
(\hat{Z} -Z^{*})_{i j} \big).
\eens
Moreover, we have for each $j \in [n]$,
\bens
n \ip{\mu, Y_j - \E (Y_j)}
& =&  \sum_{i =1, i \not=j}^{n}  \ip{\mu, \Z_j - \Z_i} = 
\sum_{i =1, i \not=j}^{n}  \ip{\mu, -\Z_i} +
\sum_{i =1, i \not=j}^{n}  \ip{\mu, \Z_j} \\
   & = &
   \sum_{i =1}^{n}  \ip{\mu, -\Z_i} +
   \ip{\mu, \Z_j} +  (n-1) \ip{\mu, \Z_j} \\
   & = & n \ip{\mu, \Z_j} - \sum_{i =1}^{n}  \ip{\mu, \Z_i}.
\eens
Hence we have by  \eqref{eq::weightsum}, and $L_j  = \sqrt{p  \gamma}\ip{\mu, Z_j}$,
\bens
\ip{P_2(\E(Y) (Y -\E(Y))^T), \hat{Z} -Z^{*}}
&  = &
4  w_1 w_2 \big(\sum_{j \in \MC_2} \hat{w}_j n L_j
-   \sum_{j \in \MC_1}  \hat{w}_j  n L_j \big)\\
& - &
 4  w_1 w_2  \big(\sum_{i =1}^{n}  L_i\big)
\big(\sum_{j \in \MC_2} \hat{w}_j  -  \sum_{j \in \MC_1}
  \hat{w}_j \big)\\
& \le &
4 w_1 w_2 n \sum_{j \in [n]} \abs{L_i} \hat{w}_j 
+   4 w_1 w_2 \abs{\sum_{i =1}^{n}  L_i}  \onenorm{(\hat{Z}    -Z^{*})}/(2n).
     \eens
\end{proofof2}

\subsection{Proof of Lemma~\ref{lemma::P2majorstar}}
  \label{sec::proofP2majorstar}
  \begin{proofof2}
Let $1 \le q < n$ denote a positive integer.
Then for the first component on the RHS of~\eqref{eq::2L} in
Lemma~\ref{lemma::P2major}, we have
\ben
  \nonumber 
4 w_1 w_2\sum_{j \in [n]} \abs{L_j} \hat{w}_j
& = &  4 w_1 w_2\sum_{j \in [n]} \abs{\ip{\mu^{(1)} -\mu^{(2)}, Z_i}} \hat{w}_j \\
  & \le & \sum_{j \in [n]} L_j^{*}  \hat{w}^{*}_j =
\sum_{j =1}^{q} L_j^{*}  \hat{w}^{*}_j +
\sum_{j =q+1}^{n} L_j^{*}  \hat{w}^{*}_j, \\
\text{ and }  \sup_{\hat{w} \in B_{\infty}^{n} \cap q B_1^{n}}
  \sum_{j \in [n]} \abs{L_j} \hat{w}_j 
& \le &
\nonumber 
\sup_{\hat{w} \in B_{\infty}^{n} }
\sum_{j =1}^{q} L_j^{*}  \hat{w}^{*}_j +
  \sup_{\hat{w} \in q B_1^{n}}
  \sum_{j =q+1}^{n} L_j^{*}  \hat{w}^{*}_j \\
  & \le &
  \label{eq::L5}
\sum_{j =1}^{q} L_j^{*}  + q  L_{q+1}^{*}
\le 2 \sum_{j =1}^{q} L_j^{*}.
\een
By Fact~\ref{fact::weights}, we have 
$\hat{w} =(\hat{w}_1, \ldots,\hat{w}_n) \in B_{\infty}^{n} \cap \lceil \frac{r_1}{2n}\rceil  B_1^{n}$, for all $\hat{Z}  \in \M_{\opt} \cap (Z^{*} +  r_1 B_1^{n \times n})$.
By setting $q = \lceil \frac{r_1}{2n}\rceil \in [n]$ in\eqref{eq::L5},  we obtain
\ben
 \sup_{\hat{w} \in B_{\infty}^{n} \cap \lceil \frac{r_1}{2n}\rceil B_1^{n}}  \sum_{j \in [n]} \abs{L_j} \hat{w}_j 
   \le 
\label{eq::L6}
2 \sum_{j =1}^{\lceil {r_1}/{(2n)}\rceil} L_j^{*}.
\een
Moreover, for the second component on the RHS of~\eqref{eq::2L}, we have by~\eqref{eq::weightsum}, 
\ben 
\nonumber 
\lefteqn{
 \sup_{\hat{Z} \in \M_{\opt} \cap Z^{*} + r_1 B_1^{n \times n}}
 \abs{\inv{n} \sum_{j =1}^{n}  L_j}  \onenorm{(\hat{Z}
   -Z^{*})}/(2n)}\\
\label{eq::L2}
&  \le &
  \sup_{\hat{w} \in B_{\infty}^{n} \cap \lceil \frac{r_1}{2n}\rceil B_1^{n}}
  \inv{n} \sum_{j =1}^{n}  \abs{L_j} \abs{\sum_{j} \hat{w}_j }
  \le  \frac{\lceil {r_1}/{2n}\rceil}{n}\sum_{j=1}^n   L^*_j
  \le \sum_{j =1}^{\lceil {r_1}/{(2n)}\rceil} L_j^{*},
\een
where \eqref{eq::L2} holds since in the second sum $\sum_{j =1}^{q} L_j^{*}$,
where $q = \lceil \frac{r_1}{2n}\rceil$,
we give a unit weight to each of the $q$ largest components $L_j^{*},
j \in [q]$, while in $\frac{q}{n} \sum_{i=1}^n L^*_j$, we give the
same weight $q/n \le 1$ to each of $n$ components $L_j^{*}, j \in [n]$.
Putting things together, we have by~\eqref{eq::2L},~\eqref{eq::L6}, and~\eqref{eq::L2},
  \bens
  \nonumber
\lefteqn{
  \sup_{\hat{Z} \in \M_{\opt} \cap (Z^{*} +  r_1 B_1^{n \times n})} \abs{  \ip{P_2(\E(Y) (Y -\E(Y))^T), \hat{Z} -Z^{*}}}} \\
  \label{eq::L1proof}
& \le &
\sup_{\hat{w} \in B_{\infty}^{n} \cap \lceil \frac{r_1}{2n}\rceil B_1^{n}}  n \sum_{j \in [n]} \abs{L_j} \hat{w}_j + 
 \sup_{\hat{w} \in B_{\infty}^{n} \cap \lceil \frac{r_1}{2n}\rceil B_1^{n}}
\abs{\sum_{j =1}^{n}  L_j } \sum_{j} \hat{w}_j \le  3 \sum_{j =1}^{\lceil \frac{r_1}{2n}\rceil} L_j^{*}.
\eens
The lemma thus holds.
\end{proofof2}

\subsection{Proof of Proposition~\ref{prop::orderL}}
\label{sec::prooforderL}
\begin{proofof2}
Recall by~\eqref{eq::definemu}, 
\bens 
L_j :=\ip{\mu^{(1)} -\mu^{(2)}, \Z_j} = \sqrt{p \gamma} \ip{\mu, \Z_j}. 
\eens 

  First, we bound for any integer $q \in [n]$, and $\tau>0$,
\ben
\nonumber 
\prob{\sum_{j=1}^q L_j^* > \tau}
& \le &
\nonumber 
\prob{\exists J \in [n], \abs{J} = q, \sum_{j \in J} \abs{L_j} > \tau}
\\
& = &
\nonumber 
\prob{\max_{ J \in [n], \abs{J} = q} \max_{u =(u_1, \ldots, u_q) \in
    \{-1, 1\}^q}  \sum_{j \in J}   u_j L_j > \tau} \\
& \le &
\nonumber 
\sum_{J: \abs{J} = q} \sum_{u \in \{-1, 1\}^q} \prob{\sum_{j \in J}   u_j L_j/\Delta > \tau/\Delta} \\
\label{eq::L4}
& \le & {n \choose q} 2^{q}
\exp\left(-\frac{(\tau/\Delta)^2}{C_1 C^2_0 (\max_{j} \twonorm{R_j \mu}^2) q} \right),
\een
where \eqref{eq::L4} holds by the union bound and the sub-gaussian tail
bound; To see this, notice that for each fixed index set $J \in [n]$
and vector $u \in \{-1, 1\}^q$, following the proof of
Lemma~\ref{lemma::anisoproj}, we obtain for $\Delta^2 = p  \gamma$ and
$\mu \in \Sp^{p-1}$,
\ben
\nonumber 
\norm{\sum_{j \in J, \abs{J} =q}   u_j L_j/ \sqrt{p \gamma}}_{\psi_2}^2
& \le & 
\sum_{j \in J, \abs{J} =q}  \norm{L_j/ \sqrt{p \gamma}}_{\psi_2}^2 =
\sum_{j \in J, \abs{J} =q}  \norm{\ip{\mu^{(1)} -\mu^{(2)}/\sqrt{p \gamma} ,    \Z_j}}_{\psi_2}^2 \\
\label{eq::qsum}
& = &
\sum_{j \in J, \abs{J} =q}  \norm{\ip{\mu, \Z_j}}_{\psi_2}^2
\le C_1 q  (C_0 \max_{j} \twonorm{R_j \mu})^2.
\een
Hence by the sub-gaussian tail bound,
cf. Lemma~\ref{lemma::anisoproj}, upon replacing $n$ with $1\le q <
n$, we have for each fixed index set $J \subset [n]$,
  \ben
\label{eq::qtail}
\prob{\sum_{j \in J: \abs{J} =q}   u_j L_j/ \Delta > \tau/ \Delta}
\le \exp\left(-\frac{( \tau/ \Delta)^2}{C_1 (C_0 \max_{j} \twonorm{R_j \mu})^2 q} \right).
\een
Set 
$\tau_q = C_3 C_0 (\max_{j} \twonorm{R_j \mu})\Delta q \sqrt{\log (2e
  n/q)}$ for some  absolute constant $C_3 > c \sqrt{C_1}$ for $c \ge 2$.
Then by~\eqref{eq::L4} and \eqref{eq::qtail}, we have
\bens
\nonumber 
\prob{\sum_{j=1}^q L_j^* > \tau_q}
& \le &
\nonumber
{n \choose q} 2^{q} \exp\left(-\frac{C^2_3 (C_0 \max_{j} \twonorm{R_j \mu})^2
    q^2 \log (2e n/q)}{C_1 (C_0 \max_{j} \twonorm{R_j \mu})^2 q} \right) \\
& \le &
\nonumber 
(2en/q)^q \exp\left(-c q \log (2e n/q)\right)  \\
& \le &
\nonumber 
e^{q \log (2 en/q)} e^{- c q \log (2e n/q)}
\le e^{- c' q \log (2e n/q)},
\eens
where ${n \choose q} 2^{q}   \le ( en/q)^q 2^q = (2 en/q)^q = 
\exp(q \log (2en/q))$.
Following the calculation as done in~\cite{GV19}, and the inequality immediately above, we have
  \bens 
  \sum_{q=1}^n  {n \choose q} 2^{q} e^{-\tau_q^2/(C_1 C_0^2 q \Delta^2)}
\le \sum_{q=1}^n e^{- c' q \log (2e n/q)}  \le \frac{c}{n^2}.
\eens
The proposition thus holds.
\end{proofof2}

\subsection{Proof of Lemma~\ref{lemma::P2P1items}}
\label{sec::proofofLemmaP2P1}
\begin{proofof2}
  Throughout this proof, denote by $\Psi= \Z \Z^T -\E (\Z \Z^T)$.
  Let $u_2$ be   the group membership vector as in~\eqref{eq::u2}.
  Let $\vecone_n =(1, \ldots, 1)$ and $\vecone_n/\sqrt{n} \in
  \Sp^{n-1}$.  Hence we have $u_2^T \vecone_n = \abs{C_1} - \abs{C_2}$.   
First, we have by Fact~\ref{fact::Zu2}, \eqref{eq::P2major} holds since
\bens
\abs{  \ip{P_2 \Psi P_1, \hat{Z} -Z^{*}}}
& = & \abs{\tr((\hat{Z} -Z^{*}) P_2 \Psi P_1) } = 
\abs{\tr(P_2 \Psi P_1 (\hat{Z} -Z^{*}) )} \\
& = &  \abs{\tr(\frac{u_2 u_2^T}{n} \Psi \inv{n}\vecone_n    \vecone_n^T (\hat{Z} -Z^{*}))} \\
& \le &
\abs{\frac{2 u_2^T \Psi \vecone_n}{n} }
\abs{\frac{\vecone_n^T (\hat{Z} -Z^{*}) u_2}{2n}}
\le 
\twonorm{\Psi} \onenorm{\hat{Z} -Z^{*}}/{n}.
\eens
Second,~\eqref{eq::P2P1P1} holds
since ${\vecone_n^T  \Psi \vecone_n}/{n} \le \twonorm{\Psi}$ and
\bens
\ip{P_2 P_1\Psi P_1, \hat{Z} -Z^{*}}
& = &
\ip{\frac{u_2  u_2^T}{n}  \inv{n}\vecone_n \vecone_n^T 
  \Psi \inv{n}\vecone_n \vecone_n^T, \hat{Z} -Z^{*}} \\
& = &
\ip{\frac{u_2}{n} \big( u_2^T  \inv{n}\vecone_n\big)
  \inv{n}\big(\vecone_n^T \Psi \vecone_n \big)\vecone_n^T, \hat{Z} -Z^{*}} \\
& = &
2 \frac{\abs{C_1} - \abs{C_2}}{n}\frac{\vecone_n^T \Psi \vecone_n}{n}
\inv{2n} \vecone_n^T( \hat{Z} -Z^{*}) u_2.
\eens
Hence by Fact~\ref{fact::Zu2},
\bens
\abs{\ip{P_2 P_1 \Psi P_1, \hat{Z} -Z^{*}}}
\nonumber 
& \le &
\abs{w_1 - w_2} \twonorm{\Psi} 
\onenorm{\hat{Z} -Z^{*}}/{n}.
\eens
To bound \eqref{eq::P2P1P1P2}, we have by Fact~\ref{fact::trP2ZZ},
\bens
\ip{P_2 P_1\Psi P_1 P_2, \hat{Z} -Z^{*}}
& = &
\ip{\frac{u_2  u_2^T}{n}  \inv{n}\vecone_n \vecone_n^T 
  \Psi \inv{n}\vecone_n \vecone_n^T \frac{u_2  u_2^T}{n}, \hat{Z} -Z^{*}} \\
& = & \ip{\frac{u_2}{n} \big( u_2^T \vecone_n/n\big) 
  \frac{\big(\vecone_n^T  \Psi  \vecone_n \big) }{n} (\vecone_n^T u_2/n) u_2^T, \hat{Z} -Z^{*}} \\
& = &
(w_1 -w_2)^2\frac{\vecone_n^T \Psi \vecone_n}{n} \ip{\frac{u_2
    u_2^T}{n}, \hat{Z} -Z^{*}},
\eens
and hence
\bens
\abs{\ip{P_2 P_1\Psi P_1 P_2, \hat{Z} -Z^{*}}}
& \le  & (w_1 -w_2)^2 \twonorm{\Psi} \onenorm{\hat{Z} -Z^{*}}/n.
\eens
Finally, it remains to show ~\eqref{eq::P22} and~\eqref{eq::P212}.
First, we rewrite the inner product as follows:
 \bens
\ip{\hat{Z} -Z^{*},  P_2 \Psi P_2}
 &  = &
 \tr((\hat{Z} -Z^{*})^T \frac{u_2 u_2^T}{n} \Psi  \frac{u_2 u_2^T}{n}
) \\
&  = &
\frac{u_2^T (\hat{Z} -Z^{*})^T u_2}{n}
\frac{u_2^T \Psi  u_2}{n} =
- \inv{n}\onenorm{\hat{Z} -Z^{*}  } \frac{u_2^T \Psi  u_2}{n};
\eens
Then~\eqref{eq::P22} holds by Fact~\ref{fact::trP2ZZ}.
Finally, we compute
\bens
\ip{P_2 P_1 \Psi P_2, \hat{Z} -Z^{*}}
&  = &
\tr((\hat{Z} -Z^{*}) \frac{u_2 u_2^T}{n} \frac{\vecone_n 
  \vecone_n^T}{n} \Psi  \frac{u_2 u_2^T}{n})\\
&  = & \tr(\frac{u_2^T}{n} (\hat{Z} -Z^{*}) u_2) \frac{u_2^T \vecone_n }{n}
\frac{\vecone_n^T \Psi  u_2}{n},\\
\; \; \text{and} \quad \ip{P_2 \Psi P_1 P_2, \hat{Z} -Z^{*}}
&  = &
\tr((\hat{Z} -Z^{*}) \frac{u_2 u_2^T}{n} \Psi  \frac{\vecone_n 
    \vecone_n^T}{n} \frac{u_2 u_2^T}{n})\\
&  = &
\tr(P_2 (\hat{Z} -Z^{*}) )
\frac{u_2^T  \Psi  \vecone_n}{n}
\frac{\vecone_n^T u_2}{n}.
\eens
Then
\bens
\abs{  \ip{P_2 P_1 \Psi P_2, \hat{Z} -Z^{*}}}
  &  = &
\abs{\ip{P_2 \Psi P_1 P_2, \hat{Z} -Z^{*}}}\\
&  \le &
\abs{w_1 - w_2} \inv{n} \onenorm{\hat{Z} -Z^{*}} \twonorm{ \Psi}.
\eens
Hence~\eqref{eq::P212} holds and the lemma is proved.
\end{proofof2}

\subsection{Proof of Lemma~\ref{lemma::P2SigmaY}}
\label{sec::proofofP2SigmaY}
\begin{proofof2}
First, we have
\bens
\lefteqn{\ip{P_2 P_1 \Psi, \hat{Z} -Z^{*}}
 = \tr((\hat{Z} -Z^{*}) P_2 P_1 \Psi) }\\
& = &
\tr((\hat{Z} -Z^{*}) \frac{u_2 u_2^T}{n} \inv{n}\vecone_n 
\vecone_n^T  \Psi)  = 
\tr((\hat{Z} -Z^{*}) \frac{u_2}{n} \big( u_2^T \inv{n}\vecone_n \big) 
\vecone_n^T  \Psi) \\
& = &
\frac{\abs{C_1} - \abs{C_2}}{n}
\tr((\hat{Z} -Z^{*}) \frac{u_2}{n} \vecone_n^T  \Psi) =
(w_1 - w_2) \vecone_n^T \Psi (\hat{Z} -Z^{*}) \frac{u_2}{n}.
\eens
Thus
\bens
\abs{  \ip{P_2 P_1 \offd(\Psi), \hat{Z} -Z^{*}}}
& = &
\abs{w_1 - w_2}
\abs{\vecone_n^T \offd(\Psi) (\hat{Z} -Z^{*}) \frac{u_2}{n}} \\ 
& \le &
2 \abs{w_1 - w_2}
  \sum_{k=1}^n
  \abs{\big(\sum_{j\not=k}^n \Psi_{k j}\big)} 
\abs{\inv{2n}\big(\sum_{j \in \MC_1}
  (\hat{Z} -Z^{*})_{k j} - \sum_{j \in \MC_2}  (\hat{Z}
  -Z^{*})_{k j} \big)  } \\
& = &
2 \abs{w_1 - w_2} \sum_{k} \abs{S_{k}} \hat{w}_k.
\eens
Therefore, we have shown~\eqref{eq::doublesum}.
Now for the diagonal component, denote by
\bens
\abs{S_{\diag} }&= & \abs{\ip{P_2 P_1 \diag(\Psi), \hat{Z} - Z^{*}}}  = 
\abs{w_1 - w_2}
\abs{\vecone_n^T \diag(\Psi) (\hat{Z} -Z^{*}) \frac{u_2}{n}} \\ 
& := &  2 \abs{w_1 - w_2} \abs{\sum_{k \in C_1} \Psi_{kk} (-\hat{w}_k) +\sum_{k \in C_2}  \Psi_{kk} \hat{w}_k} \\
&\le&2 \abs{w_1 - w_2} \sum_{k \in [n]} \abs{\Psi_{kk}} \hat{w}_k 
\le  \abs{w_1 - w_2} \max_{k} \abs{\Psi_{kk}} \onenorm{\hat{Z} - Z^{*}}/n,
\eens
where recall for $u_2$ as in~\eqref{eq::u2},
\bens
\frac{(\hat{Z} -Z^{*}) u_2}{n}  = (-\hat{w}_1, \ldots, -\hat{w}_{n_1},
\hat{w}_{n_1 +1}, \ldots, \hat{w}_{n}).
\eens
Now
  \bens
\lefteqn{ \ip{P_2 \Psi, \hat{Z} -Z^{*}}    =
    \tr((\hat{Z} -Z^{*}) P_2 \Psi)    =
    \tr(P_2 \Psi ( \hat{Z} -Z^{*})) =
 \ip{\Psi, (\hat{Z} -Z^{*}) P_2}}\\
    &   = &
    \inv{n}(u_2^T \Psi ( \hat{Z} -Z^{*}) u_2)
    =  \ip{\Psi u_2, ( \hat{Z} -Z^{*}) u_2/n} \\
      &  = &
      2\sum_{k=1}^n
      \left(\sum_{t  \in \MC_1}  \Psi_{k t} - \sum_{t  \in \MC_2}
  \Psi_{k t} \right) \cdot 
\inv{2n}\big(  \sum_{j \in \MC_1}  (\hat{Z} -Z^{*})_{k j} - \sum_{j \in \MC_2}
(\hat{Z} -Z^{*})_{k j} \big) \\
& = &
2 \big(\sum_{k \in \MC_1} \Psi_{k k} (-\hat{w}_k) 
-\sum_{t  \in \MC_2} \Psi_{tt} \hat{w}_t \big) +
2 \big(\sum_{j \in \MC_1} Q_{j} (-\hat{w}_j) 
+\sum_{j  \in \MC_2} Q_{j} \hat{w}_j \big) = Q_{\diag} + Q_{\offd}.
\eens
For the diagonal component,  we have
\bens
\abs{Q_{\diag}} & := &
\abs{\ip{P_2 \diag(\Psi), \hat{Z} -Z^{*}}  }  =
2 \abs{\sum_{k \in C_1} \Psi_{kk} (-\hat{w}_k) +\sum_{k \in C_2} (- \Psi_{kk}) \hat{w}_k }\\
&\le&
2 \sum_{k \in [n]} \abs{\Psi_{kk}} \hat{w}_k \le \max_{k}
\abs{\Psi_{kk}} \onenorm{\hat{Z} - Z^{*}}/n.
\eens
Moreover, for the off-diagonal component,  we have
\bens
\abs{Q_{\offd}}
& = &
2\abs{\sum_{j \in \MC_1} Q_{j} (-\hat{w} _j) 
  +\sum_{j  \in \MC_2} Q_{j} \hat{w}_j}
\le 2 \sum_{j \in [n]} \abs{Q_{j}}{\hat{w} _j}.
\eens
Hence~\eqref{eq::P2ZZ} holds, given that 
\bens 
\abs{\ip{P_2 \offd(\Psi), \hat{Z} -Z^{*}}}
  &  \le &
2\sum_{k \in [n]} \hat{w}_k \abs{Q_k}. 
   \eens 
 \end{proofof2}

\subsection{Proof of  Proposition~\ref{prop::orderQ} }
\label{sec::prooforderQ}
\begin{proofof2}
  Denote by $\Psi =  \Z \Z^T -\E (\Z \Z^T)$.
  Recall
 \bens
 \forall j \in [n] \quad Q_j :=
 \sum_{t  \in \MC_1, t\not=j}  \Psi_{tj} -\sum_{t  \in \MC_2, t\not=j}
 \Psi_{tj} \quad \text{ and } \quad S_j := \sum_{i\not=j}^n \Psi_{j i}.
 \eens
It remains to show~\eqref{eq::Q2} and~\eqref{eq::Q3}. 
 The two bounds follow identical steps, and hence we will only 
 prove~\eqref{eq::Q2}. 
First, we bound for any positive integer $1\le q < n$ and any $\tau >0$,
  \ben
   \nonumber
\prob{\sum_{j=1}^q Q_j^* > \tau}
& = &
\prob{\exists J \in [n], \abs{J} = q, \sum_{j \in J} \abs{Q_j} > \tau} \\
 \nonumber
& = &
\prob{\max_{ J \in [n], \abs{J} = q} \max_{y =(y_1, \ldots, y_q) \in
    \{-1, 1\}^q}  \sum_{j \in J}   y_j Q_j > \tau} \\
\label{eq::q0}
& \le &
\sum_{J: \abs{J} = q} \sum_{y \in \{-1, 1\}^q}
\prob{\sum_{j \in J}   y_j Q_j > \tau}  = q_0(\tau).
\een

We now bound $q_0(\tau)$.
Denote by $w \in \{0, -1, 1\}^n$ the $0$-extended vector of $y \in \{-1, 1\}^{\abs{J}}$, that is,
for the fixed index set $J \subset [n]$, we have $w_j =y_j, \forall j \in J$ and
$w_j= 0$ otherwise. Now we can write for $A =(a_{tj}) = w \otimes
u_2$, where $u_2$ is the group  membership vector:
\bens
\forall t \in \MC_1, j \in J, \quad a_{jt} & =&  u_{2,t} y_j = y_j,\\
\forall t \in \MC_2,  j \in J, \quad  a_{jt} & = & u_{2,t} y_j = -y_j,
\eens
and $0$ otherwise. In other words, we have $n q$ number of non-zero
entries in $A$ and each entry has absolute value 1.
Thus $\twonorm{A} \le  \fnorm{A} = \sqrt{q n}$ and
\bens
&&
\sum_{j \in J}   y_j \left(\sum_{t  \in \MC_1,  t\not=j}  \ip{\Z_t, \Z_j} -
  \sum_{t  \in \MC_2,  t\not=j}  \ip{\Z_t, \Z_j}\right)
=  \sum_{j \in [n]}   \sum_{t  \not=j}   a_{jt}\ip{\Z_t, \Z_j},
\eens
since each row of $\offd(A)$ in $J$ has $n -1$ nonzero entries and
$\abs{J} = q$.
Set   $$\tau_q = C_4 (C_0 \max_{i} \twonorm{H_i})^2 q  (\sqrt{np \log(2e n/q)} + \sqrt{nq} \log(2e n/q)).$$
Now for each fixed index set $J \in [n]$ with $\abs{J} =q$,
and a fixed vector $y \in \{-1, 1\}^q$,
\ben
   \nonumber
\prob{\sum_{j \in J}   y_j Q_j > \tau_q}
& = &
\prob{\sum_{j \in J}   y_j \left(\sum_{t  \in \MC_1, t\not=j}  \Psi_{tj} -\sum_{t  \in \MC_2, t\not=j}    \Psi_{tj} \right) > \tau_q}\\
   \nonumber
& = &   \prob{\sum_{j} \sum_{t\not=j}  a_{jt}  \ip{\Z_t, \Z_j} > \tau_q}\\
& \le &
\label{eq::q1}
2 \exp \left(- c (C_4^2 \wedge C_4) q \log(2e n/q)\right).
\een
It remains to show the inequality in \eqref{eq::q1}.
By Theorem~\ref{thm::ZHW}, we have for $q \in [n]$, $A = (a_{ij}) \in
\R^{n \times n}$, and $\tau_q = C_4 (C_0 \max_{i} \twonorm{H_i})^2 q(
\sqrt{np \log(2e n/q)} + \sqrt{nq} \log(2e n/q) )$,
\bens
\nonumber
\lefteqn{\prob{\abs{\sum_{i=1}^n  \sum_{j \not=i}^n \ip{\Z_{i}, \Z_{j}} a_{ij}} >  \tau_q}}\\
& \le &
2 \exp \left(- c\min\left(\frac{\tau_q^2}{(C_0 \max_{i} \twonorm{H_i})^4  p
      \fnorm{A}^2}, \frac{\tau_q}{(C_0 \max_{i} \twonorm{H_i})^2
      \twonorm{A}} \right)\right)  \\
& \le &
2 \exp \left(- c\min\left(\frac{(C_4 q \sqrt{np \log(2e n/q)})^2}{p n q}, 
    \frac{C_4 q \sqrt{nq} \log(2e n/q)}{\sqrt{qn}} \right)\right)  \\
& \le &
2 \exp \left(- c (C_4^2 \wedge C_4)  q \log(2e n/q)\right),
\eens
where $\max_{i} \norm{Z_i}_{\psi_2} \le C C_0 \max_{i} \twonorm{H_i}$
in the sense of \eqref{eq::Zpsi2}. Thus~\eqref{eq::q1} holds.
Similarly, we can prove the same bound holds for~\eqref{eq::duet},
following an identical sequence of arguments.
Now we have for by~\eqref{eq::q0} and~\eqref{eq::q1},
\bens
q_0(\tau_q)
& \le &
 {n \choose q} 2^{q}   \exp\left(- c (C_4^2 \wedge C_4) q \log(2e 
  n/q)\right)  \\
& = & 
 e^{q \log (2 en/q)}  \exp\left(- c (C_4^2 \wedge C_4) q \log(2e 
  n/q)\right)  \\
& \le &
2 \exp \left(- c' (C_4^2 \wedge C_4) q \log(2e n/q)\right),
\eens
where ${n \choose q} 2^{q}   \le ( en/q)^q 2^q = (2 en/q)^q = 
\exp(q \log (2en/q))$. Hence following the calculation as done in~\cite{GV19},  we have
\bens 
\sum_{q=1}^n q_0(\tau_q) \le
\sum_{q=1}^n  {n \choose q} 2^{q} \exp\left(- c (C_4^2 \wedge C_4) q
  \log(2e   n/q)\right)
\le \sum_{q=1}^n e^{- c' q \log (2e n/q)}  \le \frac{c'}{n^2}.
\eens
The proposition thus holds.
\end{proofof2}

\section{The oracle estimator and its properties}
\label{sec::oracleSDP}
The presentation of this section follows from~\cite{Zhou23a}, which
we 
include for self-containment.
Both $A$ and $B$ are defined to bridge the gap between $YY^T$ and $R =
\E(Y)\E(Y)^T$, and the bias term 
\bens 
\E B -R & = & \E Y Y^T - \E (Y) \E(Y)^T- \E \lambda (E_n - I_n) - \E \tau 
I_n 
\eens 
is substantially smaller than $\E A-R$ in the operator norm, under
assumption (A2).
Proposition~\ref{prop::optsol} shows that
the optimization goals of SDP1 and SDP~\eqref{eq::sdpmain} are equivalent to 
that of {\bf Oracle SDP}.
\begin{proposition}{\textnormal{\citep{Zhou23a}}}
  \label{prop::optsol}
  The optimal solutions $\hat{Z}$ as in SDP~\eqref{eq::sdpmain}
  must have their  diagonals set to $I_n$.
Thus, the set of optimal solutions $\hat{Z}$ in~\eqref{eq::sdpmain}
coincide with
those on the convex subset $\M_{\opt}$ as in \eqref{eq::moptintro},
\ben
\label{eq::Aquiv}
\argmax_{Z \in  \M^{+}_{G} }  \ip{A , Z}
& = &  \argmax_{Z \in \M_{\opt}}   \ip{A , Z}
  =  \argmax_{Z \in \M_{\opt}} \ip{\tilde{A}, Z}\\
\label{eq::optsolAB}
  & = &
  \argmax_{Z \in \M_{\opt}}  (  \ip{A , Z} - \E \tau \ip{I_n, Z} ).
   \een 
\end{proposition}

Theorem~\ref{thm::reading} follows immediately from
Lemma~\ref{lemma::EBRtilt} and Theorem~\ref{thm::YYaniso}.
\begin{theorem}{\textnormal{{\bf{($R$ is the leading term)}}~\citep{Zhou23a}}}
  \label{thm::reading}
Suppose the conditions in Theorem~\ref{thm::exprate} hold.
Then with probability at least $1-2\exp(-cn)$, for the oracle
$B$ as in~\eqref{eq::defineBintro},
\bens 
\twonorm{B - R} & \le &  \xi n p \gamma,
\quad \text{ where the reference $R$ is as defined in~\eqref{eq::Rtilt}.}
\eens
\end{theorem}

\subsection{Proof of Lemma~\ref{lemma::EBRtilt}}
\label{sec::proofofEBR}
The proof follows from~\cite{Zhou23a}, Section H.
We include it here for self-containment.
We have the following facts about $R$ as defined in~\eqref{eq::Rtilt}.
\begin{fact}
  \label{fact::Rtrace}
First, $\vecone^T_n R \vecone_n = 0$. Hence
  \ben
  \label{eq::negTR}
\vecone^T_n\offd( R) \vecone_n & = &\vecone^T_n R 
\vecone_n- \tr(R) = - \tr(R) = -n p \gamma w_2 w_1, \;\text{ where} \\
 \nonumber
 \tr(R)/(p \gamma)
 & = & n w_2^2 w_1 +  n w_1^2 w_2=w_1 w_2 n \quad \text{ and }  \quad \twonorm{R} = \tr(R) =w_1 w_2 n p \gamma.
\een
Thus ${\tr(R)}/{n}  =  p \gamma w_1 w_2$ and for $n \ge 4$,
\ben
\label{eq::Rop}
\twonorm{\frac{\tr(R)}{n-1} \left( I_n-
    {E_n}/{n}\right)} & \le & \frac{n}{n-1}
p \gamma w_1 w_2 \le p \gamma /3,
\een
where $\twonorm{I_n-{E_n}/{n}} = 1$ since $I_n- {E_n}/{n}$ is a projection matrix.
\end{fact}

\begin{proofof}{Lemma~\ref{lemma::EBRtilt}}
We have by  Proposition~\ref{prop::biasfinal},
\ben
\nonumber                                   
\norm{\E B - R}
  & =&  \norm{W_0 - \mathbb{W} -\frac{\tr(R)}{(n-1)} (I_n - \frac{E_n}{n})}
  \\
 \label{eq::normEBR}
& \le &
  \norm{W_0}                  
 +\norm{ \mathbb{W}} +\norm{\frac{\tr(R)}{(n-1)} (I_n - \frac{E_n}{n}) },
 \een
 where the $\norm{\cdot}$ is understood to be either the operator or
 the cut norm. Then
\ben 
\label{eq::W0opnorm}
\twonorm{W_0}  & \le & \abs{V_1 -  V_2}  (w_2 \vee w_1).
\een
Moreover,
\ben
\label{eq::W3main}
\twonorm{\mathbb{W}}
& \le& \abs{V_1 - V_2} (w_1 \vee w_2).
\een
Combining~\eqref{eq::Rop}, \eqref{eq::normEBR},~\eqref{eq::W0opnorm},
and~\eqref{eq::W3main}, we have for
$w_{\min} := w_1 \wedge w_2$ and $ n  \xi \ge \inv{2 w_{\min}}$,
\bens
\twonorm{\E B - R} & \le  &
\twonorm{W_0} + \twonorm{\mathbb{W}} + 
\twonorm{\frac{\tr(R)}{(n-1)} (I_n - \frac{E_n}{n})} \\
& \le & 2 \abs{V_1 -  V_2} (w_1 \vee w_2) + p \gamma/3 \\
& = &
\frac{2}{3} \xi n p \gamma (1- w_{\min}) + p \gamma/3 
\le \frac{2}{3} \xi n p \gamma,
\eens
where we use the fact that
\bens
\frac{2}{3} \xi p \gamma n w_{\min} \ge  p \gamma/3  \; \; \text{ since
} \; \; 2 \xi n w_{\min} \ge  1.
\eens
The lemma thus holds for the general setting; when $V_1 = V_2$, we 
show the improved bounds in Lemma~H.7.~\cite{Zhou23a}.
\end{proofof}

\subsection{Proof sketch of Theorem~\ref{thm::YYaniso}}
\label{sec::reduction}
Let $Y$ be as in Definition~\ref{def::estimators}. 
We provide a proof outline for Theorem~\ref{thm::YYaniso} in this 
section. 
The idea of decomposition and reduction appears in the local analysis
as well; See Proofs for Lemma~\ref{lemma::S1} in
Section~\ref{sec::proofofS1main}, as well as
Corollary~\ref{coro::sumY}.
We have by~\eqref{eq::projection} and \eqref{eq::covproj}, and the triangle inequality,
\bens
\twonorm{\Lambda} =
\twonorm{YY^T -\E(YY^T)} & \le &
\twonorm{\hat{\Sigma}_Y - \Sigma_Y} +
\twonorm{\E(Y)(Y-\E(Y))^T +(Y-\E(Y)) \E(Y)^T},
\eens
where $\hat\Sigma_Y =  (Y-\E(Y))(Y-\E(Y))^T$ and
\ben
\label{eq::ZZproj}
\twonorm{\hat{\Sigma}_Y - \Sigma_Y}
& \le &   \twonorm{I-P_1}\twonorm{\Z \Z^T -\E (\Z \Z^T)}\twonorm{I-P_1}.
\een
To control $\norm{M_Y} =\norm{\E(Y)(Y-\E(Y))^T + (Y-\E(Y))(\E(Y))^T},$
we need to bound the projection of each mean-zero random vector $\Z_j,
\forall j \in [n]$, along the direction of $v := \mu^{(1)}-\mu^{(2)}$.
In other words,  a particular direction for which we compute the 
one-dimensional marginals is the direction between the two centers, namely,
$\mu^{(1)}$ and $\mu^{(2)}$.
Let $c, c', c_1, c_5, C_3, C_4, \ldots $ be absolute constants.
  \begin{lemma}{\textnormal{{\bf (Projection: probabilistic view)}~\cite{Zhou23a}}}
  \label{lemma::projerr}
  Suppose conditions in Theorem~\ref{thm::YYaniso} hold.
Then we have with probability at least $1 -2 \exp(-c n)$, 
\bens
\label{eq::Yprop}
\twonorm{\E(Y)(Y-\E(Y))^T + (Y-\E(Y))(\E(Y))^T }
& \le & 4 \sqrt{n}  \sqrt{w_1 w_2} \sup_{q \in S^{n-1}} \abs{\sum_{i}
q_i \ip{\Z_i, \mu^{(1)} -\mu^{(2)}}} \\
& \le & 2 C_3 C_0 (\max_{i} \twonorm{H^T_i \mu}) n \sqrt{p \gamma}.
\eens 
\end{lemma}
Hence Theorem~\ref{thm::YYaniso} follows immediately from 
Lemma~\ref{lemma::projerr} and Theorem~\ref{thm::YYcovcorr}. 
Lemma~\ref{lemma::projerr} follows from  
Lemma~\ref{lemma::tiltproject} and the sub-gaussian concentration of 
measure bound in~\eqref{eq::proanis}.
On the other hand, controlling the second component in 
\eqref{eq::projection} amounts to the problem of covariance 
estimation given the mean matrix $\E(Y)$, in view
of~\eqref{eq::ZZproj} and Theorem~\ref{thm::YYcovcorr}.
See~\cite{Zhou23a} for proofs; cf. Lemma 7.1. and Theorem 7.2 therein.

\bibliography{final,subgaussian,clustering}

\end{document}